\documentclass[leqno]{amsart}

\usepackage[foot]{amsaddr}

\usepackage{color}
\usepackage{geometry}
\usepackage{subfiles}
\usepackage{graphicx}
\usepackage[section]{placeins}
\usepackage{hyperref}
\usepackage{stmaryrd}
\usepackage{amsmath,amssymb,amsfonts,amsthm,textcomp}
\usepackage{mathtools}
\usepackage{tensor}
\usepackage{multicol}
\usepackage{subfig}
\usepackage{csquotes}
\usepackage{stackrel}
\usepackage{tikz}
\usepackage{siunitx}
\usepackage{soul}
\usepackage{pdfpages}
\usepackage{dsfont}
\usepackage{booktabs}
\usepackage[ruled,lined]{algorithm2e}
\graphicspath{
              {Images/}
             }

\newfont{\tenbfsl}{cmbxti9 scaled 1200}
\newfont{\tenbbb}{msbm10}
\newfont{\svnbbb}{msbm8}

\DeclareMathOperator*{\argmin}{arg\,min}

\newcommand{\bs}[1]{\boldsymbol{#1}}
\newcommand{\cl}[1]{\mathcal{#1}}

\newcommand{\bb}[1]{\mathbb{#1}}

\newcommand{\red}[1]{{\color{red} #1}}

\newcommand{\di}[1]{\,\mathrm{d}{#1}}

\newcommand{\trans}{\scriptscriptstyle\mskip-1mu\top\mskip-2mu}

\newcommand{\tr}{\mathrm{tr}\mskip2mu}

\newtheorem{cor}{Corollary}
\newtheorem{lem}{Lemma}
\newtheorem{prop}{Proposition}
\newtheorem{rmk}{Remark}
\newtheorem{set}{Definition}

\newtheorem{asu}{Assumption}

\begin{document}

\title[Randomized quasi-Monte Carlo estimator for nested integration]{Double-loop randomized quasi-Monte Carlo estimator for nested integration}
\author{Arved Bartuska$^{\#,1,2}$, Andr\'{e} Gustavo Carlon$^{1}$, Luis Espath$^3$, Sebastian Krumscheid$^5$, \& Ra\'{u}l Tempone$^{1,2,4}$}
\address{$^1$Department of Mathematics, RWTH Aachen University, Geb\"{a}ude-1953 1.OG, Pontdriesch 14-16, 161, 52062 Aachen, Germany}
\address{$^2$King Abdullah University of Science \& Technology (KAUST), Computer, Electrical and Mathematical Sciences \& Engineering Division (CEMSE), Thuwal 23955-6900, Saudi Arabia}
\address{$^3$School of Mathematical Sciences, University of Nottingham, Nottingham, NG7 2RD, United Kingdom}
\address{$^4$Alexander von Humboldt Professor in Mathematics for Uncertainty Quantification, RWTH Aachen University, Germany}
\address{$^5$Scientific Computing Center, and Institute for Applied and Numerical Mathematics, Karlsruhe Institute of Technology, Germany}
\email{$^\#$arved.bartuska@gmail.com}

\date{\today}

\begin{abstract}
\noindent
Nested integration of the form $\int f\left(\int g(\bs{y},\bs{x})\di{}\bs{x}\right)\di{}\bs{y}$, characterized by an outer integral connected to an inner integral through a nonlinear function $f$, is a challenging problem in various fields, such as engineering and mathematical finance.
The available numerical methods for nested integration based on Monte Carlo (MC) methods can be prohibitively expensive owing to the error propagating from the inner to the outer integral. Attempts to enhance the efficiency of these approximations using the quasi-MC (QMC) or randomized QMC (rQMC) method have focused on either the inner or outer integral approximation.
This work introduces a novel nested rQMC method that simultaneously addresses the approximation of the inner and outer integrals. The method leverages the unique nested integral structure to offer a more efficient approximation mechanism. As the primary contribution, we derive asymptotic error bounds for the bias and variance of our estimator, along with the regularity conditions under which these bounds can be attained. Incorporating Owen's scrambling techniques, we address integrands exhibiting infinite variation in the Hardy--Krause sense, enabling theoretically sound error estimates. Moreover, we derive a truncation scheme for applications in the context of expected information gain estimation. We verify the estimator quality through numerical experiments by comparing the computational efficiency of the nested rQMC method against standard nested MC estimation to highlight the computational savings and enhanced applicability of the proposed approach.
$\cdot$
65C05 
$\cdot$
65D30 
$\cdot$
65D32 

\end{abstract}

\maketitle

\section*{Acknowledgments}

This publication is based upon work supported by the King Abdullah University of Science and Technology (KAUST) Office of Sponsored Research (OSR) under Award No.~OSR-2019-CRG8-4033, the Alexander von Humboldt Foundation, the Deutsche Forschungsgemeinschaft (DFG, German Research Foundation) -- 333849990/GRK2379 (IRTG Hierarchical and Hybrid Approaches in Modern Inverse Problems), and was partially supported by the Flexible Interdisciplinary Research Collaboration Fund at the University of Nottingham Project ID 7466664. The authors thank Prof.~Robert Scheichl, Prof.~Fabio Nobile, and Dr.~Yang Liu for their helpful and constructive comments. 

\tableofcontents                        


\section{Introduction}
Nested integrals are ubiquitous in numerous fields (e.g., geology \cite{God18}, mathematical finance \cite{Xu20}, medical decision-making \cite{Fan22}, and optimal experimental design (OED) \cite{Rya03}). The Monte Carlo (MC) method is a commonly applied approximation technique for general integrals, especially in high dimensions. For nested integrals, both integration steps can be approximated using the MC method, resulting in the double-loop MC (DLMC) estimator \cite{Rya03}. To control the statistical error of an MC approximation of a single integral via the central limit theorem (CLT) \cite{Dur19}, a sample size of $\cl{O}(TOL^{-2})$ is required for an error tolerance $TOL>0$, related to an estimator's variance (to be specified later). In contrast, using DLMC to approximate nested integrals results in worse complexity, with the overall required number of samples of $\cl{O}(TOL^{-3})$ to achieve the same error tolerance~\cite{Rya03}.
The two MC estimators appearing in the DLMC estimator are connected through a nonlinear function; thus, the statistical error of the inner MC estimator induces bias in the DLMC estimator. The sample sizes of both MC estimators must be precisely controlled to ensure that the error in the DLMC estimator is below $TOL$ with a certain confidence, controlling the corresponding statistical and bias errors. Improving the performance of the DLMC estimator is the goal of ongoing research, with some approaches proposing Laplace approximation \cite{Lon13}, importance sampling \cite{Bec18}, and multilevel MC (MLMC) methods \cite{Fan22}.

The randomized quasi-MC (rQMC) method \cite{Caf98, Hic98, Nie92, Owe03, Dic10, Lec18, Cra76} is a promising technique to improve the efficiency of the standard MC method (i.e., the required computational work to meet a tolerance) while maintaining a nonintrusive sampling structure. The rQMC estimator uses deterministic points from a low-discrepancy sequence and randomizes the entire sequence while maintaining a low-discrepancy structure. Randomization removes bias in deterministic QMC estimators, leaving only the statistical error, which can be controlled using the Chebyshev inequality or the CLT under specific conditions \cite{Tuf04, Lec10, Lec18}. Given appropriate regularity properties of the integrand, the rQMC method can improve the order of convergence of the approximation error compared to the MC estimator. For sufficiently regular integrands, the computational work required by the rQMC estimator to achieve a tolerance $TOL$ is almost $\cl{O}(TOL^{-1})$, where the multiplicative term can increase with the integrand's dimension. Certain constructions of rQMC estimators based on scrambled digital nets even achieve comparable accuracy requiring computational work of almost $\cl{O}(TOL^{-\frac{2}{3}})$ if the integrand is regular enough. Owen \cite{Owe06} demonstrated that rQMC can yield improvements over the MC method, even for integrands with singularities at the boundaries of the integration domain (i.e., integrands with classically unfavorable regularity properties) for a randomization scheme introduced in \cite{Owe95}. Introducing weighted function spaces permits rQMC methods to remain competitive for high, or even infinite, dimensional integration domains under suitable regularity conditions~\cite{Dic13}. Finally, higher-order rQMC methods enable even faster convergence at the detriment of exponential dependence on the dimension \cite{Dic13}.

The estimation of nested integrals via MC was studied in \cite{Rya03} along with the theoretical error bounds in the context of the expected information gain (EIG) for OED. A more general case of a possibly discontinuous outer integrand was considered in \cite{Hon09}, while the study \cite{Rai18} investigated conditions under which nested MC estimators converge. Attempts to improve upon the efficiency of MC estimators were undertaken, e.g., in \cite{Gob17}, where the Markov chain MC method was combined with nonparametric regression for the rare event case, and in \cite{Fos17}, where variational inference was employed to create several estimators for OED applications. 

The evaluation of the integrand itself often requires an additional approximation at a specific level of accuracy. Examples of these accuracy levels include the number of time steps in temporal approximation schemes, the mesh size in finite element method (FEM) schemes, and the number of samples of inner integral estimators. A strategy for improving the effectiveness of standard MC estimators is MLMC \cite{Gil08}, where MC estimations for integrands of various accuracy levels are combined into one highly efficient estimator. Under certain conditions, MLMC estimation can be performed without a loss of accuracy from the added approximation (i.e., at the cost of standard MC estimation, $\cl{O}(TOL^{-2})$). An MLMC estimator for OED was derived in \cite{God20}, where the number of inner MC samples defined the levels, and \cite{Bec20} derived another MLMC estimator, where the number of inner MC samples and a FEM approximation defined the levels. The randomized MLMC method is an extension of the MLMC method, where the number of levels is randomized, thus removing the bias from the inner approximation. This method was applied in the context of OED in \cite{God22}; however, randomization is not always possible and places considerable restrictions on the integrands.

Researchers \cite{Gob22} have investigated the optimal error tolerance that can be reached using rQMC, given a fixed number of samples and a specific confidence level. They demonstrated that combining rQMC with robust estimation improves error tolerance values. The rQMC concepts were applied in the context of OED in \cite{Dro18} but only to the outer integral of a nested integration problem. In \cite{Dro18}, the inner integral was approximated using the MC method. For several numerical experiments, the sample standard deviation was reduced for this scheme compared to using the MC method for both integrals, which was demonstrated for a fixed number of outer and inner samples. In the exploratory work \cite{Kaa24}, Kaarnioja and Schillings proposed two nested rQMC estimators for OED governed by random PDEs. The randomness in the PDE coefficients is parametric in terms of a sequence of random variables. For a fixed discretization of the PDE, based on rank-1 lattice rules with random shifts, a computational cost of almost $\cl{O}(TOL^{-2})$ is derived for their first estimator to control the rQMC error. Regularity conditions on the integrands are derived to achieve independence of the parameter dimensions. Their second estimator, moreover, employs a sparse tensor structure to achieve a computational cost of almost $\cl{O}(TOL^{-1})$ to control the rQMC error. This last estimator relies on the very particular nested expectation structure of the EIG. A truncation of the Gaussian observation noise depending on $TOL$ is used in these estimators to address singularities that are ubiquitous in this problem. This truncation introduces a factor in their root-mean-square-error bound that grows exponentially with the noise dimension, as in\cite[Theorem 6.2]{Kaa24}. In the present work, we demonstrate that such behavior does not occur for bounds based on the milder Hardy--Krause variation or Owen's boundary growth condition, 
see Corollary~\ref{cor:total.error}.

In the work \cite{Fan22}, an rQMC method was employed to approximate the outer integral of a nested estimator in medical decision-making. These authors compared this rQMC method with an MLMC approach and standard nested (i.e., double-loop) MC by specifying a target mean squared error and observing the number of samples necessary to reach the target. The interplay between rQMC and MC estimators in the nested setting was not theoretically analyzed. Both MLMC and rQMC have similar performance results in practice, depending on the number of parameters and other model complexity  measures.

In \cite{God18}, nested integrals were approximated using MLMC techniques. The number of inner samples was increased at each level to reduce the bias induced in the outer approximation by the variance of the inner approximation. The rQMC estimator estimated the inner integral and reduced the variance, requiring a smaller sample size in the inner loop to achieve the error tolerance in the MLMC setting in \cite[Thm.~3.1]{Gil08}. This outcome was presented theoretically and via examples.
A similar approach was followed in \cite{Xu20}, in which a discontinuous inner integrand was approximated using a sigmoid (i.e., a smooth function), improving the performance of the rQMC method. The methods presented in \cite{God18} and \cite{Xu20} can achieve an error tolerance $TOL$ at a computational cost of $\cl{O}(TOL^{-2})$. 

This work applies the rQMC method to both integrals to build a double-loop rQMC (rDLQMC) estimator, reducing the number of samples necessary to estimate the nested integrals up to a specific error tolerance $TOL$ at a cost close to $\cl{O}(TOL^{-10/9})$. Moreover, we derive sufficient regularity conditions for both the outer and inner integrands and demonstrate the advantage of replacing either of the MC approximations with rQMC approximations. We also consider the case in which the inner integrand is only approximately provided in terms of a FEM discretization, resulting in additional bias in the outer approximation. We provide approximate error bounds for the bias and statistical error of the proposed estimator and derive a nearly optimal setting, verified via numerical examples from Bayesian OED. To make the derived error bounds applicable in this context, we verify the previously derived conditions for a simple EIG setting with truncated observation noise. The additional truncation error introduced by this scheme is shown to affect the total error only by an insignificant contribution. Although truncation makes the outer integrand sufficiently regular for each fixed tolerance, the corresponding multiplicative constant term grows as the truncation region expands. Thus the improved rate provided by rQMC estimators based on scrambled digital nets is retained only for the inner estimator, whereas the truncation mainly reduces constants for the outer integrand for the tolerances considered, resulting in an overall cost of almost $\cl{O}(TOL^{-4/3})$.

This work is structured as follows. Section~\ref{sec:MC.QMC} provides a brief overview of the MC and rQMC methods, including bounds on the $L^1$ and $L^2$ errors via the Koksma--Hlawka inequality \cite{Hla61, Nie92} and more elaborate error bounds proposed in \cite{Owe06} and later generalized in \cite{He23}. The subsequent Section~\ref{sec:nested} introduces the proposed novel nested rQMC estimator. As the main contributions of this work, asymptotic error bounds on the number of inner and outer samples are derived in Propositions~\ref{prop:B.DLQ} and \ref{prop:V.DLQ} along with conditions under which these rates can be attained. Based on these bounds, the nearly optimal setting for the rDLQMC estimator is derived in Proposition~\ref{prop:W.DLQ}. A simple example where both the inner and outer integrands are polynomials demonstrates the convergence rates achievable in the optimal setting for this estimator and presents the effect of the integrand dimensions on multiplicative terms in the estimation error. In Section~\ref{sec:EIG.estimation}, an overview of the EIG setting is provided, where the outer integrand is the logarithm and the inner integrand is the likelihood of data observations. Moreover, a truncation of the observation noise is introduced, and the precise conditions derived in Section~\ref{sec:nested} are verified for the EIG with truncation in Lemma~\ref{lem:boundary.condition} and~\ref{lem:cond}. However, verifying the Assumptions on the particular experiment model derived in Appendix~\ref{app:cond} and \ref{app:Condition.42} is not the subject of this work to maintain generality. In Lemma~\ref{lem:trunc}, we bound the truncation error introduced in Appendix~\ref{app:Condition.42}; and in Corollary~\ref{cor:total.error}, we present the total error of the EIG estimation for our estimator. This section also presents two Bayesian OED examples, where nested integrals frequently arise. The superior performance of the proposed estimators compared to other estimators is observed numerically in accordance with our theoretical results. A first example demonstrates the effect of the truncation of Gaussian observation noise and prior for the parameters of interest. This example allows for a closed-form solution and serves to showcase the effects of the noise level. The second example explores an application from solid mechanics involving the solution to a partial differential equation (PDE) with favorable regularity properties, demonstrating the practical applicability of the rDLQMC estimator. 
Finally, Section~\ref{sec:conclusion} provides concluding remarks and an outlook for future research directions.

























\section{Brief overview of Monte Carlo and randomized quasi-Monte Carlo integration}\label{sec:MC.QMC}
Before addressing the nested integration case, which is the focus of this work, we discuss the basic concepts for approximating integrals using MC, QMC, and rQMC for the reader's convenience.

\subsection{Monte Carlo method}\label{sec:MC}
We let $\varphi:[0,1]^d\to\bb{R}$, where $d$ is a positive integer, be square integrable. Then, the integral
\begin{equation}\label{eq:integral}
  I=\int_{[0,1]^d}\varphi(\bs{x})\di{}\bs{x},
\end{equation}
can be approximated using the MC estimator as follows:
\begin{equation}\label{eq:monte.carlo}
I_{\rm{MC}}\coloneqq\frac{1}{M}\sum_{m=1}^{M}\varphi(\bs{x}^{(m)}).
\end{equation}
The MC method uses random points $\bs{x}^{(1)},\ldots,\bs{x}^{(M)}$, which are independent and identically distributed (iid) samples from the uniform distribution $\cl{U}\left([0,1]^d\right)$, to approximate $I$ in \eqref{eq:integral}. Using the CLT to analyze the error of the MC estimator reveals that
\begin{equation}\label{eq:CLT.prob}
    \bb{P}\left(|I - I_{\rm{MC}}| \le \frac{C_{\alpha} \sqrt{\bb{V}[\varphi]}}{\sqrt{M}}\right) \ge 1 - \alpha,
\end{equation}
as $M\to\infty$. Here, $\alpha\in(0,1)$ is a parameter, $C_{\alpha} = \Phi^{-1}(1-\alpha/2)$, $\Phi^{-1}$ denotes the inverse cumulative distribution function (CDF) of the standard normal distribution, and $\bb{V}[g]$ represents the variance of the integrand. Alternatively, Chebyshev's inequality could be applied to obtain an error estimate similar to \eqref{eq:CLT.prob}, 
although with a larger constant $C_{\alpha}=1/\sqrt{\alpha}$.

\subsection{Quasi-Monte Carlo method}\label{sec:QMC}
The MC method converges with probability one to the true value of the integral by the strong law of large numbers, but the rate of $M^{-(1/2)}$ can be improved for certain integrands~\cite{Caf98}. We can employ the QMC method, which achieves a better convergence rate by using specific points to exploit the regularity properties of the integrand. The QMC estimator to approximate the integral \eqref{eq:integral} is given by
\begin{equation}\label{eq:det.quasi.monte.carlo}
 I_{\rm{Q}} \coloneqq \frac{1}{M}\sum_{m=1}^M\varphi(\bs{t}^{(m)}),
\end{equation}
where $\bs{t}^{(1)},\ldots,\bs{t}^{(M)}\in[0,1]^d$ are selected from a deterministic sequence of points with low discrepancy, \cite{Nie92, Hic98}. To analyze the error of this estimator, we introduce two useful concepts. The star discrepancy~\cite[Def.~(5.2)]{Caf98}, defined as
\begin{equation}
    D^{\ast}(\bs{t}^{(1)},\ldots,\bs{t}^{(M)})\coloneqq \sup_{(s_1,\ldots,s_d)\in[0,1]^d}\left|\frac{1}{M}\sum_{m=1}^M\mathds{1}_{\{\bs{t}^{(m)}\in\prod_{i=1}^d [0,s_i)\}}-\prod_{i=1}^ds_i\right|\;,
\end{equation}
indicates how evenly the points $\bs{t}^{(1)},\ldots,\bs{t}^{(M)}$ are distributed in $[0,1]^d$. This concept is independent of the particular integrand $\varphi$. The next concept we recall is the Hardy--Krause variation~\cite[Def.~(5.8)]{Caf98}, which is recursively defined for continuously differentiable integrands $\varphi$ as follows:
\begin{equation}\label{eq:VHK}
    V_{\rm{HK}}(\varphi)\coloneqq \int_{[0,1]^d}\left|\left(\prod_{i=1}^d\frac{\partial^d}{\partial z_i}\right)\varphi(\bs{z})\right|\di{}\bs{z}+\sum_{i=1}^dV_{\rm{HK}}(\varphi(\bs{z})|_{z_i=1}),
\end{equation}
where $\varphi(\bs{z})|_{z_i=1}:[0,1]^{d-1}\to\bb{R}$ represents the restriction of $\varphi(\bs{z})$ to $\varphi(z_1,\ldots,z_i=1,\ldots,z_d)$. For $\varphi:[0,1]\to\bb{R}$,
\begin{equation}
    V_{\rm{HK}}(\varphi)\coloneqq \int_{[0,1]}\left|\frac{\di{}}{\di{}z}\varphi(z)\right|\di{}z,
\end{equation}
which depends on the regularity of the integrand $\varphi$ and is independent of the points $\bs{t}^{(1)},\ldots,\bs{t}^{(M)}$. The mixed first-order derivatives are assumed to be continuous to guarantee the uniqueness of \eqref{eq:VHK}.
 The Koksma--Hlawka inequality \cite{Hla61, Nie92, Caf98} provides the following bound for the error of the QMC estimator \eqref{eq:det.quasi.monte.carlo}
\begin{equation}\label{eq:QMC.bound}
    |I-I_{\rm{Q}}|\leq D^{\ast}(\bs{t}^{(1)},\ldots,\bs{t}^{(M)})V_{\rm{HK}}(\varphi).
\end{equation}
Low-discrepancy sequences satisfying 
\begin{equation}\label{eq:star.discrepancy}
D^{\ast}(\bs{t}^{(1)},\ldots,\bs{t}^{(M)})=\cl{O}(M^{-1}\log(M)^{d-1})
\end{equation}
are provided in \cite{Nie92}. To suppress log terms in the notation, the bound \eqref{eq:star.discrepancy} can be written as follows:
\begin{equation}
    D^{\ast}(\bs{t}^{(1)},\ldots,\bs{t}^{(M)})\leq C_{\epsilon,d}M^{-1+\epsilon},
\end{equation}
for any $\epsilon>0$, where $C_{\epsilon,d}>0$ is independent of $M$, and $C_{\epsilon,d}\to\infty$ as $\epsilon\to 0$. Two common low-discrepancy sequences are digital sequences \cite{Nie92, Owe95, Owe03, Dic10} and rank-1 lattice rules \cite{Hic98}. The concept behind digital sequences is to split $[0,1]^d$ into equally spaced subintervals as follows:
\begin{equation}
    \prod_{i=1}^d\left[\frac{a_{i}}{2^{b_{i}}},\frac{a_{i}+1}{2^{b_{i}}}\right),
\end{equation}
where $a_i,b_i$ are integers, with $b_i\geq 0$, $0\leq a_i\leq 2^{b_i}-1$, and $1\leq i\leq d$ such that each subinterval contains the same fraction of each nonoverlapping block of $M$ points in the sequence. A number of points $M$, such that $\log_2(M)\in\bb{N}$, must be used to achieve the best results for this sequence type (for details on the construction of such sequences, refer to \cite{Dic10}). 

The deterministic QMC method introduced above has several practical disadvantages. The estimator \eqref{eq:det.quasi.monte.carlo} is biased, and the first point $\bs{t}^{(1)}$ is the origin for many common low-discrepancy sequences, which is problematic for standard classes of integrands (see Remark~\ref{rmk:general.domains}). The QMC method suffers from the curse of dimensionality for a finite number of points $M$ but an increasing dimension $d$. Furthermore, the Hardy--Krause variation $V_{\rm{HK}}$ is challenging to compute. Most importantly, the QMC error in \eqref{eq:QMC.bound} is finite for some integrands for which $V_{\rm{HK}}$ is infinite \cite{Lec18}. A critical class of such integrands is the set of functions with singularities at one or more corners of the domain $[0,1]^d$ (e.g., the log or inverse CDF of the normal distribution; Remark~\ref{rmk:general.domains}). A more accurate bound is provided in Owen \cite{Owe06}, depending on how fast the mixed derivatives of the integrand increase when approaching the boundary. 

\subsection{Randomized Quasi-Monte Carlo method}\label{sec:rQMC}
A popular strategy for obtaining a computable error estimate for QMC methods is using points from randomized low-discrepancy sequences. Randomization also removes bias and ensures that the origin is included with only negligible probability. The deterministic points $\bs{t}^{(1)},\ldots,\bs{t}^{(M)}$ are mapped via randomization, denoted by $\bs{\rho}$, to obtain randomized low-discrepancy points $\bs{x}^{(m)}=\bs{\rho}(\bs{t}^{(m)})$, where $1\leq m\leq M$ and the same randomization $\bs{\rho}$ is applied to each point. Thus, the points $\bs{x}^{(1)},\ldots,\bs{x}^{(M)}$ are uniformly distributed on $[0,1]^d$ but are not independent. We applied Owen's scrambling for digital sequences \cite{Nie92, Owe95, Owe03} throughout this work, in which a random permutation $\bs{\rho}_{\rm{sc}}$ is applied to the base two digits of $\bs{t}^{(m)}$, where $1\leq m\leq M$, to build on the error bounds derived in \cite{Owe06, He23} for these low-discrepancy sequences. This randomization scheme is only guaranteed to preserve the low discrepancy with probability one (refer to \cite{Owe95, Dic10} for details). More specifically, we applied the Sobol' sequence \cite{Sob67} throughout this work.
Using a particular randomization $\bs{\rho}$, the rQMC estimator is given by
\begin{equation}\label{eq:quasi.monte.carlo}
 I_{\rm{rQ}} \coloneqq \frac{1}{R}\sum_{r=1}^R\frac{1}{M}\sum_{m=1}^M\varphi(\bs{x}^{(r,m)}),
\end{equation}
and the shortened notation is
\begin{equation}
    \bs{x}^{(r,m)}\coloneqq\bs{\rho}^{(r)}(\bs{t}^{(m)}), \quad 1\leq r\leq R, \, 1\leq m\leq M.
\end{equation}
 The applications of rQMC methods in practice often involve integrands with an unbounded Hardy--Krause variation, e.g., the logarithm. Thus, a modified version of \cite[Thm.~5.7]{Owe06} is presented, providing an error bound for integrands that do not increase too rapidly as their argument approaches the boundary of the integration domain. A recent work~\cite{Liu23} derived a nonasymptotic error bound to account for the fact that optimal rQMC rates are typically not observed for a finite $M$ because of the log term in \eqref{eq:star.discrepancy} and also provides a minor correction to~\cite{Owe06}. Here, we carefully track constant terms appearing in these error estimates, which allows for the analysis of integral approximations nested in other integral approximations \eqref{eq:dlqmc}.
\begin{asu}[Boundary growth condition]\label{eq:boundary.growth}
Let $\bs{x}^{(m)}\sim\cl{U}\left([0,1]^d\right)$, where $1\leq m\leq M$, be from a randomized low-discrepancy sequence with one randomization such that $\bb{E}[D^{\ast}(\bs{x}^{(1)},\ldots,\bs{x}^{(M)})]\leq C_{\epsilon,d}M^{-1+\epsilon}$ for all $\epsilon>0$, where $C_{\epsilon,d}>0$ is independent of $M$, and $C_{\epsilon,d}\to\infty$ as $\epsilon\to 0$. Furthermore, let $\varphi:[0,1]^d\to\bb{R}$ and assume that there exists $0<b<\infty$ and $A_i>0$ for $1\leq i\leq d$, where $\max_{i}A_i<1/2$, such that
\begin{equation}
    \left|\left(\prod_{j\in u}\frac{\partial^u}{\partial z_j}\right)\varphi(\bs{z})\right|\leq b\prod_{i=1}^d\min (z_i,1-z_i)^{-A_i-\mathds{1}_{\{i\in u\}}}
\end{equation}
for all $\bs{z}=(z_1,\ldots,z_d)\in(0,1)^d$ and all $u\subseteq\{1,\ldots,d\}$ with the convention that $(\prod_{j\in u}\partial^u/\partial z_j)\varphi\equiv \varphi$ for $u=\emptyset$.
\end{asu}
Observe that Assumption~\ref{eq:boundary.growth} implies that $\varphi\in L^2([0,1]^{d})$.

\begin{prop}[Thm.~5.7 in \cite{Owe06}]\label{prop:owen}
    Given Assumption~\ref{eq:boundary.growth}, the rQMC estimator \ref{eq:quasi.monte.carlo} satisfies the following $L^1([0,1]^{d})$ error bound:
    \begin{equation}
        \bb{E}[|I-I_{\rm{rQ}}|]\le  b B_A C_{\epsilon,d} M^{-1+\epsilon+\max_iA_i},
    \end{equation}
for all $\epsilon>0$, where $B_A\to\infty$ as $\min_{i}A_i\to 0$, and $C_{\epsilon,d}\to\infty$ as $\epsilon\to 0$.
\end{prop}
The following extension of this result for the $L^2([0,1]^{d})$ error was provided in \cite{He23},
\begin{prop}[Thm.~3.1 in \cite{He23}]\label{cor:he}
Given Assumption~\ref{eq:boundary.growth}, the rQMC estimator \ref{eq:quasi.monte.carlo} satisfies the following $L^2([0,1]^{d})$ error bound:
\begin{equation}\label{eq:RMSE}
        \bb{E}[|I-I_{\rm{rQ}}|^2]^{\frac{1}{2}}\leq b B_A C_{\epsilon,d} M^{-1+\epsilon+\max_iA_i},
    \end{equation}
for all $\epsilon>0$, where $B_A\to\infty$ as $\min_{i}A_i\to 0$, and $C_{\epsilon,d}\to\infty$ as $\epsilon\to 0$.
    \end{prop}
The difference between the estimators \eqref{eq:monte.carlo} and \eqref{eq:quasi.monte.carlo} is in the points used to evaluate the function to be integrated. For practical purposes, we employed $R$ independent randomizations of the same low-discrepancy sequence, $\bs{x}^{(r,1)},\ldots,\bs{x}^{(r,M)}$, where $1\leq r\leq R$, to acquire an error bound using the Chebyshev inequality; that is,
\begin{equation}\label{eq:CLT.bound}
    \bb{P}\left(|I - I_{\rm{rQ}}| \le C_{\alpha} \sqrt{\bb{V}[I_{\rm{rQ}}]}\right) \ge 1 - \alpha
\end{equation}
for $\alpha\in(0,1)$, where $C_{\alpha}=1/\sqrt{\alpha}$. For practical reasons, the variance of the rQMC estimator \eqref{eq:quasi.monte.carlo} can be estimated using the sample variance:
\begin{equation}\label{eq:qmc.var}
\tilde{\bb{V}}[I_{\rm{rQ}}] \coloneqq
\frac{1}{R(R-1)}\sum_{r=1}^R\left(\frac{1}{M}\sum_{m=1}^M\varphi(\bs{x}^{(r,m)})-I_{\rm{rQ}}\right)^2.
\end{equation}
For certain functions $\varphi$ with desirable properties, such as smoothness and boundedness, more precise statements are possible \cite{Gob22, Owe08}. A CLT-based error estimate similar to \eqref{eq:CLT.bound} only holds asymptotically as $R\to\infty$. The estimate can still be used in practice to obtain a confidence interval of the error in \eqref{eq:CLT.bound}; however, keeping $R$ fixed and letting $M\to\infty$ is occasionally problematic, as noted in \cite{Tuf04, Lec10}. Specifically, the convergence of the distribution of the estimator \eqref{eq:quasi.monte.carlo} to a normal distribution is not guaranteed.

\begin{rmk}[Integration over general domains]\label{rmk:general.domains}
The rQMC method is commonly defined for integration over the unit cube and uniform random variables. For integrals over general domains (e.g., normally distributed random variables), the corresponding inverse CDF can be applied to maintain the general shape of the estimators \eqref{eq:monte.carlo} and \eqref{eq:quasi.monte.carlo}. This approach introduces singularities at the boundaries, leading to an unbounded Hardy--Krause variation \eqref{eq:VHK}. Thus, Proposition~\ref{prop:owen} or Proposition~\ref{cor:he} is required to obtain meaningful error bounds rather than relying on the Koksma--Hlawka inequality \eqref{eq:QMC.bound}.
\end{rmk}
Finally, for certain integrands, the rQMC method using scrambled Sobol' points can achieve convergence rates of almost three in the variance error (see~\cite[Theorem 6.25]{Dic13}). This is the case if the integrand has bounded generalized Hardy--Krause variation of order one. The key difference between the Hardy--Krause variation introduced above and the generalized Hardy--Krause variation of order one is that the mixed first partial derivatives must be in $L^2$ rather than in $L^1$ to satisfy this stricter regularity criterion. See \cite[Chapter 6]{Dic13} or \cite[Chapter 13]{Dic10} for the full definition. This applies, for example, to the inner integrand in the EIG setting (that is, the likelihood of data observations). We discuss the application of this criterion to the outer integrand and provide a short definition in Appendix~\ref{app:GHK}.

\section{Nested integration}\label{sec:nested}
The previous section discusses the fundamentals of the rQMC estimator \eqref{eq:quasi.monte.carlo}. Next, this section resumes the primary focus of this work, establishing the rDLQMC estimator \eqref{eq:dlqmc.S.R} for nested integration problems with a discretized integrand, deriving asymptotic error bounds in the number of samples and discretization parameters, and analyzing the optimal work for this estimator to meet a tolerance goal.
\subsection{Estimator and notation}
\begin{set}[Nested integral]
A nested integral is defined as follows:
\begin{align}\label{eq:double.loop}
  I {}&= \int_{[0,1]^{d_1}}f\left(\int_{[0,1]^{d_2}}g(\bs{y},\bs{x})\di{}\bs{x}\right)\di{}\bs{y},\nonumber\\
{}&=\int_{[0,1]^{d_1}}\tilde{f}\left(\bs{y}\right)\di{}\bs{y},
\end{align}
where
\begin{equation}
    \tilde{f}(\bs{y})\coloneqq f\left(\int_{[0,1]^{d_2}}g(\bs{y},\bs{x})\di{}\bs{x}\right),
\end{equation}
and $f:\bb{R}\to\bb{R}$ is nonlinear. Furthermore, $g:[0,1]^{d_1}\times[0,1]^{d_2}\to\bb{R}$ defines a nonlinear relation between $\bs{x}\in[0,1]^{d_2}$ and $\bs{y}\in[0,1]^{d_1}$, where $d_1, d_2$ are positive integers.
\end{set}

\begin{rmk}[Application: expected information gain]\label{ex:EIG}
The following integral is nested:
\begin{equation}\label{exa:EIG}
    I=\int_{[0,1]^{d_1}}\log\left(\frac{1}{\det(2\pi\bs{\Sigma})^{\frac{1}{2}}}\int_{[0,1]^{d_2}}e^{-\frac{1}{2}\left\lVert \bs{G}(\bs{y}_1)+\bs{\Sigma}^{\frac{1}{2}}\Phi^{-1}(\bs{y}_2)-\bs{G}(\bs{x})\right\rVert^2_{\bs{\Sigma}^{-1}}}\di{}\bs{x}\right)\di{}\bs{y},
\end{equation}
where
\begin{equation}
    g(\bs{y},\bs{x})= \frac{1}{\det(2\pi\bs{\Sigma})^{\frac{1}{2}}}e^{-\frac{1}{2}\left\lVert \bs{G}(\bs{y}_1)+\bs{\Sigma}^{\frac{1}{2}}\Phi^{-1}(\bs{y}_2)-\bs{G}(\bs{x})\right\rVert^2_{\bs{\Sigma}^{-1}}},
\end{equation}
and $f\equiv \log$. In these equations, $\bs{G}$ denotes a nonlinear function, $d_1>d_2$, $\bs{y}=(\bs{y}_1,\bs{y}_2)\in[0,1]^{d_1}$, where $\bs{y}_1$, $\bs{x}\in[0,1]^{d_2}$, $\bs{y}_2\in[0,1]^{d_1-d_2}$, and $\bs{\Sigma}\in\bb{R}^{(d_1-d_2)\times(d_1-d_2)}$ is a positive definite matrix. Furthermore, $\Phi^{-1}$ is the inverse CDF of the standard normal in dimension $d_1-d_2$. The integral \eqref{exa:EIG} is typically not solvable in closed form, motivating the application of numerical integration techniques to approximate $I$. This integral appears in the EIG, which is applied as a utility function in OED. In this context, $\bs{G}$ represents a model of the experiment and regularly necessitates the solution of a PDE, motivating a FEM discretization:
\begin{align}
    g(\bs{y},\bs{x}){}&\approx g_h(\bs{y},\bs{x}),\nonumber\\
    {}&\coloneqq \frac{1}{\det(2\pi\bs{\Sigma})^{\frac{1}{2}}}e^{-\frac{1}{2}\left\lVert \bs{G}_h(\bs{y}_1)+\bs{\Sigma}^{\frac{1}{2}}\Phi^{-1}(\bs{y}_2)-\bs{G}_h(\bs{x})\right\rVert^2_{\bs{\Sigma}^{-1}}},
\end{align}
where $\bs{G}_h$, $h>0$, denotes the FEM approximation of $\bs{G}$ (see Assumption~\ref{asu:FEM}). Furthermore, $\bs{x}$ and $\bs{y}_1$ are parameters of interest, and $\bs{\Sigma}^{\frac{1}{2}}\Phi^{-1}(\bs{y}_2)$ denotes the observation noise.
\end{rmk}
\begin{rmk}[Nested integral]
If $f$ in \eqref{eq:double.loop} is linear, the nested integral \eqref{eq:double.loop} simplifies to a nonnested integral in dimension $d_1+d_2$. If $g$ is separable (i.e., $g(\bs{y},\bs{x})\equiv g_{1}(\bs{y})g_{2}(\bs{x})$, for some functions $g_1$ and $g_2$), then the inner and outer integrals in \eqref{eq:double.loop} can be separately approximated. Numerical approximation schemes are substantially more efficient in these cases than for truly nested integrals.
\end{rmk}
The following DLMC estimator \cite{Rya03} results from a standard double-loop MC method to approximate a nested integral \eqref{eq:double.loop}:
\begin{equation}\label{eq:dlmc}
I_{\rm{DLMC}} \coloneqq \frac{1}{N}\sum_{n=1}^Nf\left(\frac{1}{M}\sum_{m=1}^Mg(\bs{y}^{(n)},\bs{x}^{(n,m)})\right),
\end{equation}
where the points $\bs{y}^{(n)}$, where $1\leq n\leq N$, are iid and sampled from $\cl{U}\left([0,1]^{d_1}\right)$, and $\bs{x}^{(n,m)}$, where $1\leq n\leq N$, and $1\leq m\leq M$, are iid and sampled from $\cl{U}\left([0,1]^{d_2}\right)$. The standard MC estimator \eqref{eq:monte.carlo} for a single integral is unbiased and has a variance that decreases with the number of samples, which holds for the inner MC estimator in \eqref{eq:dlmc}, where the variance decreases with the number of inner samples $M$. The outer MC estimator in \eqref{eq:dlmc} has a variance that decreases with $N$ and has a bias relative to the size of the variance of the inner integral estimate. Thus, many inner and outer samples are typically required to control the bias and variance of this estimator, significantly limiting its practical usefulness, predominantly for computationally demanding problems.

Directly evaluating the function $g$ is often not feasible. For example, if evaluating $g$ requires solving a PDE, only a FEM approximation $g_h$ with a discretization parameter $h$ may be accessible.
Thus, the following assumption is made.
\begin{asu}[Discretization rate]\label{asu:FEM}
    Let $g_h$ be an approximation of $g$ such that
\begin{equation}
\left|\int_{[0,1]^{d_1}}f\left(\int_{[0,1]^{d_2}}g_h(\bs{y},\bs{x})\di{}\bs{x}\right)\di{}\bs{y}-\int_{[0,1]^{d_1}}f\left(\int_{[0,1]^{d_2}}g(\bs{y},\bs{x})\di{}\bs{x}\right)\di{}\bs{y}\right| \leq C_{\rm{disc}}h^{\eta},
\end{equation}
where $\eta>0$ denotes the $h$-convergence rate, and $C_{\rm{disc}}>0$ is independent of $h$. The work of evaluating $g_h$ is assumed to be of order $\cl{O}(h^{-\gamma})$, for some $\gamma > 0$.
\end{asu}
Moreover, the variance of the DLMC estimator \eqref{eq:dlmc} for discretized integrand $g_h$ is assumed to be finite as $h\to 0$. Unless stated otherwise, a discretized version $g_h$ is applied in the numerical estimators, and the subscript was omitted for simplicity thereafter. For square-integrable and thrice-differentiable $f$ and square-integrable $g$, the bias of the DLMC estimator has the following upper bound:
\begin{equation}\label{eq:bias.bound}
  |\mathbb{E}[I_{\rm{DLMC}}]-I|\leq C_{\rm{disc}}h^{\eta}+\frac{C_{\mathrm{MC}}^{(3)}}{M},
\end{equation}
where $C_{\rm{MC}}^{(3)}>0$ represents a constant related to the variance of the inner MC estimation in \eqref{eq:dlmc}, and $C_{\rm{disc}}>0$ could be different from that introduced in Assumption~\ref{asu:FEM}. Furthermore, the variance of the DLMC estimator has the following upper bound:
\begin{equation}\label{eq:variance.bound.MC}
  \bb{V}[I_{\rm{DLMC}}]\leq \frac{C_{\mathrm{MC}}^{(1)}}{N}+\frac{C_{\mathrm{MC}}^{(2)}}{NM},
\end{equation}
where $C_{\rm{MC}}^{(1)},C_{\rm{MC}}^{(2)} >0$ are constants \cite{Bec18}.
The optimal work of the DLMC estimator for a specific error tolerance $TOL>0$ is given by
\begin{equation}\label{eq:Opt.Work}
    W_{DLMC}^{\ast}=\cl{O}\left(TOL^{-(3+\frac{\gamma}{\eta}) }\right),
\end{equation}
as $TOL\to 0$. Proofs for the specific case of approximating the EIG in Bayesian OED are provided in \cite{Bec18}, and adapting them to the general case of a nested integral is straightforward. 
Proofs for general integrands without additional inner discretization (e.g., FEM) are provided in \cite{Hon09}, where nonsmooth integrands are also considered. These proofs assume that the outer integrand is thrice differentiable; however, this is not a requirement for the estimator \eqref{eq:dlmc} to converge. In \cite{Rai18}, similar error bounds were derived using more relaxed requirements.
\begin{rmk}[Assumptions on the discretized inner integrand]
    Assumptions similar to the ones stated above regarding the convergence of the estimator with discretized inner integrand must be verified on a case-by-case basis. It is beyond the scope of this work to analyze this discretization error further for particular integrands or discretization methods. We refer to~\cite{Bec18} for a discussion of the DLMC estimator for the EIG application with FEM discretization.
\end{rmk}
To attain smaller error bounds and optimal work, in this work we replace both MC approximations in \eqref{eq:dlmc} with rQMC approximations and arrive at the rDLQMC estimator, defined below.
\begin{set}[rDLQMC estimator]
The rDLQMC estimator of a nested integral \eqref{eq:double.loop} is defined as follows:
\begin{equation}\label{eq:dlqmc.S.R}
I_{\rm{rDLQ}}^{(S,R)} \coloneqq \frac{1}{S}\sum_{s=1}^S\frac{1}{N}\sum_{n=1}^Nf\left(\frac{1}{R}\sum_{r=1}^R\frac{1}{M}\sum_{m=1}^Mg(\bs{y}^{(s,n)},\bs{x}^{(s,n,r,m)})\right),
\end{equation}
where $f:\bb{R}\to\bb{R}$ is nonlinear and $g:[0,1]^{d_1}\times[0,1]^{d_2}\to\bb{R}$ defines a nonlinear relation between $\bs{x}$ and $\bs{y}$. The samples take the following shape
\begin{equation}
    \bs{y}^{(s,n)}\coloneqq\bs{\rho}^{(s)}(\bs{u}^{(n)}), \quad 1\leq s\leq S, \, 1\leq n\leq N,
\end{equation}
where $\bs{\rho}^{(s)}$, with $1\leq s\leq S$, is the independent randomization of points in a digital sequence $\bs{u}^{(1)},\ldots,\bs{u}^{(N)}$ in dimension $d_1$, using Owen's scrambling, and
\begin{equation}\label{eq:x.nrm}
    \bs{x}^{(s,n,r,m)}\coloneqq\bs{\rho}^{(s,n,r)}(\bs{t}^{(m)}), \quad \, 1\leq s\leq S, \,1\leq n\leq N, \, 1\leq r\leq R, \, 1\leq m\leq M,
\end{equation}
where $\bs{\rho}^{(s,n,r)}$, with $1\leq s\leq S$, $1\leq n\leq N$, and $1\leq r\leq R$ are independent randomizations of points in a digital sequence $\bs{t}^{(1)},\ldots,\bs{t}^{(M)}$ in dimension $d_2$, using Owen's scrambling.
\end{set}
\begin{rmk}[Number of randomizations]
    Increasing the number of randomizations $R$ and $S$ in the estimator \eqref{eq:dlqmc.S.R} reduces the bias and statistical error at the MC rates in \eqref{eq:bias.bound} and \eqref{eq:variance.bound.MC}, respectively. Only one randomization (i.e., scramble) of the outer samples in the estimator in \eqref{eq:dlqmc} was applied, and numerous low-discrepancy points $N$ and $M$ were employed to maximize the effect of the low-discrepancy structure. Choosing independent randomizations of the inner samples for each outer sample enables the bounding of the variance in Proposition~\ref{prop:V.DLQ}. The rDLQMC estimator with one randomization is introduced below to avoid clutter and is applied throughout the following propositions. For practical applications, including the numerical experiments in Sections~\ref{sec:nested} and~\ref{sec:EIG.estimation}, the variance contribution of the inner and outer estimators was estimated separately, using pilot runs with a fixed number of randomizations $S,R>1$.
\end{rmk}
\begin{set}[rDLQMC estimator with $S,R=1$]
The rDLQMC estimator for a nested integral \eqref{eq:double.loop} with $S,R=1$ is defined as follows:
\begin{align}\label{eq:dlqmc}
I_{\rm{rDLQ}} {}&\coloneqq I_{\rm{rDLQ}}^{(1,1)},\nonumber\\
{}&=\frac{1}{N}\sum_{n=1}^Nf\left(\frac{1}{M}\sum_{m=1}^Mg(\bs{y}^{(n)},\bs{x}^{(n,m)})\right),
\end{align}
where $f:\bb{R}\to\bb{R}$ is nonlinear, and $g:[0,1]^{d_1}\times[0,1]^{d_2}\to\bb{R}$ defines a nonlinear relation between $\bs{x}$ and $\bs{y}$. The samples take the following shape:
\begin{equation}
    \bs{y}^{(n)}\coloneqq\bs{\rho}(\bs{u}^{(n)}), \quad 1\leq n\leq N,
\end{equation}
where $\bs{\rho}$ is a randomization of points in a digital sequence $\bs{u}^{(1)},\ldots,\bs{u}^{(N)}$ in dimension $d_1$, applying Owen's scrambling, and
\begin{equation}
    \bs{x}^{(n,m)}\coloneqq\bs{\rho}^{(n)}(\bs{t}^{(m)}), \quad \, 1\leq n\leq N, \, 1\leq m\leq M,
\end{equation}
where $\bs{\rho}^{(n)}$, with $1\leq n\leq N$, are independent randomizations of points in a digital sequence $\bs{t}^{(1)},\ldots,\bs{t}^{(M)}$ in dimension $d_2$, applying Owen's scrambling. 
\end{set}
Next, we analyze the error of the rDLQMC method for general integrands $f$ and $g$ based on certain regularity assumptions. These assumptions are verified in Appendix~\ref{app:cond} and~\ref{app:Condition.42} for the EIG application, where $f$ is the logarithm and $g$ is the likelihood of data observations based on auxiliary assumptions on the experiment model and a truncation scheme for the observation noise. These bounds will be revisited in Corollary~\ref{cor:optimal.work.s}, which only applies to certain rQMC constructions (for example, using the scrambled Sobol' sequence) under stricter regularity assumptions. First, the error is split into bias and statistical errors, respectively, and each is estimated individually, as follows:
\begin{equation}\label{eq:rDLQMC.error}
|I_{\rm{rDLQ}}-I|\leq\underbrace{|\mathbb{E}[I_{\rm{rDLQ}}]-I|}_{\text{bias error}}+\underbrace{|I_{\rm{rDLQ}}-\mathbb{E}[I_{\rm{rDLQ}}]|}_{\text{statistical error}}.
\end{equation}
The following notation is defined:
\begin{equation}
\hat{g}_M(\bs{y})\coloneqq\frac{1}{M}\sum_{m=1}^Mg(\bs{y},\bs{x}^{(m)}),
\end{equation}
\begin{align}\label{eq:bar.g}
\bar{g}(\bs{y}){}&\coloneqq\mathbb{E}[\hat{g}_M(\bs{y})|\bs{y}],\nonumber\\
{}&=\mathbb{E}[g(\bs{y},\bs{x})|\bs{y}],\nonumber\\
{}&=\int_{[0,1]^{d_2}} g(\bs{y},\bs{x})\di{}\bs{x},
\end{align}
where the conditioning is on the randomization $\bs{\rho}$ of the outer variable $\bs{y}=\bs{\rho}(\bs{u})$, and
\begin{equation}
\Delta g_M(\bs{y})\coloneqq\hat{g}_M(\bs{y})-\bar{g}(\bs{y}),
\end{equation}
where the dependence on $\bs{x}^{(m)}$,  where $1\leq m\leq M$, is suppressed, and the index $n$ is omitted for brevity. The notation $\bar{g}_h(\bs{y})$ is defined analogously to~\eqref{eq:bar.g} for the discretized inner integrand $g_h$ as presented in Assumption~\ref{asu:FEM}. Next, the assumptions for establishing the bound on the bias error of the rDLQMC estimator are presented.
\subsection{Bias error analysis}
\begin{asu}[Integrability assumption for the composition of the outer and inner integrand]\label{asu:f.3.diff}
There exists $0<\delta_{K} <1$ and $1\leq p\leq\infty$ such that $\tilde{f}(\cdot)\equiv f(K\bar{g}_{h}(\cdot)):[0,1]^{d_1}\to\bb{R}$ is an element of the Sobolev space $W^{3,p}([0,1]^{d_1})$ for all $K\in(1-\delta_{K},1+\delta_{K})$, where $\bar{g}_{h}$ is as in \eqref{eq:bar.g}.
\end{asu}
It is demonstrated in the proof of Proposition~\ref{prop:B.DLQ} below that the open interval $(1-\delta_{K},1+\delta_{K})$ in Assumption~\ref{asu:f.3.diff}, where $0< \delta_{K}<1$, can be arbitrarily small.
\begin{asu}[Inverse inequality for the inner integrand]\label{asu:g.L4}
    Let $g_{h}(\cdot,\cdot):[0,1]^{d_1}\times[0,1]^{d_2}\to\bb{R}$ and assume that there exists $0\leq k<\infty$ such that
    \begin{equation}\label{eq:reg:g}
    \sup_{\bs{y}\in[0,1]^{d_1}} \left\lvert \frac{V_{\rm{HK}}(g_{h}(\bs{y},\cdot))}{ \bar{g}_{h}(\bs{y})} \right\rvert \leq k\;,
\end{equation}
with $V_{\rm{HK}}$ as in \eqref{eq:VHK}, $g_h$ and $\bar{g}_h$ as in Assumption~\ref{asu:FEM} and \eqref{eq:bar.g} as $h\to 0$. Furthermore, if Assumption~\ref{asu:f.3.diff} holds, assume that there exists $1\leq q\leq\infty$ with $1/p+1/q=1$ for $p$ as in Assumption~\ref{asu:f.3.diff} such that $\bar{g}_{h}\in L^{3q}([0,1]^{d_1})$.
\end{asu}
\begin{prop}[Bias of the rDLQMC estimator]\label{prop:B.DLQ}
Given Assumptions~\ref{asu:FEM}--\ref{asu:g.L4}, if $|f'''|$ is monotonic, the bias of the rDLQMC estimator \eqref{eq:dlqmc} has the following upper bound:
\begin{equation}\label{eq:Bias.constraint}
  |\mathbb{E}[I_{\rm{rDLQ}}]-I|\leq C_{\mathrm{disc}}h^{\eta}+ \frac{\bb{E}[|\bar{g}_{h}|^2\left|f''(\bar{g}_{h})\right|]k^2C_{\epsilon,d_2}^2}{2M^{2-2\epsilon}}+\frac{\bb{E}[\left|\bar{g}_{h}\right|^3\left|f'''(K\bar{g}_{h})\right|]k^3C_{\epsilon,d_2}^3}{6M^{3-3\epsilon}},
\end{equation}
for any $\epsilon>0$, where $C_{\epsilon,d_2}\to\infty$ as $\epsilon\to 0$.
\end{prop}

\begin{proof}
We begin by splitting the bias error into a discretization bias and an inner sampling bias on the discretized level as follows:
\begin{align}
    {}&|\mathbb{E}[I_{\rm{rDLQ}}]-I|\nonumber\\{}&=\left|\bb{E}\left[\frac{1}{N}\sum_{n=1}^Nf\left(\frac{1}{M}\sum_{m=1}^Mg_h(\bs{y}^{(n)},\bs{x}^{(n,m)})\right)\right]-\int_{[0,1]^{d_1}}f\left(\int_{[0,1]^{d_2}}g(\bs{y},\bs{x})\di{}\bs{x}\right)\di{}\bs{y}\right|,\nonumber\\
    {}&= \left|\bb{E}\left[f\left(\frac{1}{M}\sum_{m=1}^Mg_h(\bs{y},\bs{x}^{(n,m)})\right)\right]-\int_{[0,1]^{d_1}}f\left(\int_{[0,1]^{d_2}}g(\bs{y},\bs{x})\di{}\bs{x}\right)\di{}\bs{y}\right|,\nonumber\\
    {}&\leq \underbrace{\left|\bb{E}\left[f\left(\frac{1}{M}\sum_{m=1}^Mg_h(\bs{y},\bs{x}^{(n,m)})\right)\right]-\int_{[0,1]^{d_1}}f\left(\int_{[0,1]^{d_2}}g_h(\bs{y},\bs{x})\di{}\bs{x}\right)\di{}\bs{y}\right|}_{\text{inner sampling bias}},\nonumber\\{}&\quad+\underbrace{\left|\int_{[0,1]^{d_1}}f\left(\int_{[0,1]^{d_2}}g_h(\bs{y},\bs{x})\di{}\bs{x}\right)\di{}\bs{y}-\int_{[0,1]^{d_1}}f\left(\int_{[0,1]^{d_2}}g(\bs{y},\bs{x})\di{}\bs{x}\right)\di{}\bs{y}\right|}_{\text{discretization bias}}.
\end{align}
It follows directly from Assumption~\ref{asu:FEM} that the discretization bias is bounded by $C_{\rm{disc}}h^{\eta}$ for some $C_{\rm{disc}}, \eta>0$.

For the bias from the inner sampling, we omit the index $h$ for simplicity. The following Taylor expansion for $f(\hat{g}_M(\bs{y}))$ around $\bar{g}(\bs{y})$ is considered:
\begin{align}\label{eq:Taylor.IM}
  f(\hat{g}_M(\bs{y}))-f(\bar{g}(\bs{y}))= f'(\bar{g}(\bs{y}))\Delta g_M(\bs{y})&{}+\frac{1}{2}f''(\bar{g}(\bs{y}))\left(\Delta g_M(\bs{y})\right)^2\nonumber\\
  {}&+ \frac{1}{2}\left(\Delta g_M(\bs{y})\right)^3\int_0^1f'''\left(\bar{g}(\bs{y})+s \Delta g_M(\bs{y})\right)(1-s)^2\di{}s.
\end{align}
The absolute value of the expectation of the left-hand side of \eqref{eq:Taylor.IM} is the inner sampling bias of the rDLQMC estimator. Thus, each term on the right-hand side must be estimated. For the first term in \eqref{eq:Taylor.IM},
\begin{align}
    |\bb{E}[f'(\bar{g}(\bs{y}))\Delta g_M(\bs{y})]|={}&|\bb{E}[\bb{E}[f'(\bar{g}(\bs{y}))\Delta g_M(\bs{y})|\bs{y}]]|,\nonumber\\
    ={}&|\bb{E}[f'(\bar{g}(\bs{y}))\underbrace{\bb{E}[\Delta g_M(\bs{y})|\bs{y}]}_{=0}]|,\nonumber\\
    ={}&0.
\end{align}
For the second term in \eqref{eq:Taylor.IM},
\begin{align}
    \left|\bb{E}\left[f''(\bar{g}(\bs{y}))\left(\Delta g_M(\bs{y})\right)^2\right]\right|={}&\left|\bb{E}\left[\bb{E}\left[f''(\bar{g}(\bs{y}))\left(\Delta g_M(\bs{y})\right)^2\Bigg|\bs{y}\right]\right]\right|,\nonumber\\
    ={}&\left|\bb{E}\left[f''(\bar{g}(\bs{y}))\bb{E}\left[\left(\Delta g_M(\bs{y})\right)^2\Bigg|\bs{y}\right]\right]\right|.
\end{align}
For the inner expectation, by Assumption~\ref{asu:g.L4} and the Koksma--Hlawka inequality \eqref{eq:QMC.bound},
\begin{align}\label{eq:assumption.4.KH}
    \bb{E}\left[\left(\Delta g_M(\bs{y})\right)^2\Bigg|\bs{y}\right]{}&=\bb{E}\left[\left|\int_{[0,1]^{d_2}}g(\bs{y},\bs{x})\di{}\bs{x}-\frac{1}{M}\sum_{m=1}^{M}g(\bs{y},\bs{x}^{(m)})\right|^2\Bigg|\bs{y}\right],\nonumber\\
    {}&\leq (V_{\rm{HK}}(g(\bs{y},\cdot)))^2C_{\epsilon, d_2}^2M^{-2+2\epsilon},\nonumber\\
    {}&\leq (\bar{g}(\bs{y}))^2k^2C_{\epsilon, d_2}^2M^{-2+2\epsilon},
\end{align}
for any $\epsilon>0$, where $C_{\epsilon, d_2}\to\infty$ as $\epsilon\to0$. By H\"{o}lder's inequality, combining the two previous expressions results in
\begin{align}\label{eq:Holder.finite}
    \left|\bb{E}\left[f''(\bar{g}(\bs{y}))\left(\Delta g_M(\bs{y})\right)^2\right]\right|
    {}&\leq \bb{E}\left[\left|f''(\bar{g}(\bs{y}))\right|(\bar{g}(\bs{y}))^2\right] k^2C_{\epsilon, d_2}^2M^{-2+2\epsilon},\nonumber\\
    {}&\leq \bb{E}\left[\left|f''(\bar{g}(\bs{y}))\right|^p\right]^{\frac{1}{p}}\bb{E}\left[|\bar{g}(\bs{y})|^{2q}\right]^{\frac{1}{q}} k^2C_{\epsilon, d_2}^2M^{-2+2\epsilon},\nonumber\\
    {}&<\infty,
\end{align}
per Assumption~\ref{asu:f.3.diff} and \ref{asu:g.L4}.
For the third term in \eqref{eq:Taylor.IM},
\begin{align}\label{eq:higher.order.term.bias}
    {}&\left|\bb{E}\left[\left(\Delta g_M(\bs{y})\right)^3\int_0^1f'''\left(\bar{g}(\bs{y})+s \Delta g_M(\bs{y})\right)(1-s)^2\di{}s\right]\right|,\nonumber\\
    {}&=\left|\bb{E}\left[\bb{E}\left[\left(\Delta g_M(\bs{y})\right)^3\int_0^1f'''\left(\bar{g}(\bs{y})+s \Delta g_M(\bs{y})\right)(1-s)^2\di{}s\Bigg|\bs{y}\right]\right]\right|,\nonumber\\
    {}&\leq\bb{E}\left[\bb{E}\left[\left|\Delta g_M(\bs{y})\right|^3\int_0^1\left|f'''\left(\bar{g}(\bs{y})+s \Delta g_M(\bs{y})\right)\right|(1-s)^2\di{}s\Bigg|\bs{y}\right]\right].
\end{align}
Owen's scrambling yields a low-discrepancy sequence with probability one. From Assumption~\ref{asu:g.L4}, assuming that $|f'''|$ is monotonically decreasing, with probability one,
\begin{align}\label{eq:monotonicity}
    \int_0^1\left|f'''\left(\bar{g}(\bs{y})+s \Delta g_M(\bs{y})\right)\right|(1-s)^2\di{}s{}&\leq \int_0^1\left|f'''\left(\bar{g}(\bs{y})-s |\Delta g_M(\bs{y})|\right)\right|(1-s)^2\di{}s,\nonumber\\
    {}&\leq \left|f'''\left(\bar{g}(\bs{y})- |\Delta g_M(\bs{y})|\right)\right|\int_0^1(1-s)^2\di{}s,\nonumber\\
    {}&= \frac{1}{3}\left|f'''\left(\bar{g}(\bs{y})- |\Delta g_M(\bs{y})|\right)\right|,\nonumber\\
    {}&\leq \frac{1}{3}\left|f'''\left(\bar{g}(\bs{y})- C_{\epsilon, d_2}M^{-1+\epsilon}V_{\rm{HK}}(g(\bs{y},\cdot))\right)\right|,\nonumber\\
    {}&\leq \frac{1}{3}\left|f'''\left(\bar{g}(\bs{y})- C_{\epsilon, d_2}M^{-1+\epsilon}k|\bar{g}(\bs{y})|\right)\right|.
\end{align}
We introduce the following notations:
\begin{equation}
    \tilde{k}\coloneqq {\rm{sgn}}(\bar{g}(\bs{y}))C_{\epsilon, d_2}M^{-1+\epsilon}k
\end{equation}
and
\begin{equation}
    K\coloneqq 1-\tilde{k},
\end{equation}
where $\rm{sgn(\cdot)}$ is the sign function, and note that $-1< \tilde{k}<1$ for $M$ greater than some $M_0(\tilde{k})$ and thus $0<K<2$. From this, it follows that
\begin{align}
    \frac{1}{3}\left|f'''\left(\bar{g}(\bs{y})- C_{\epsilon, d_2}M^{-1+\epsilon}k|\bar{g}(\bs{y})|\right)\right|
    {}&= \frac{1}{3}\left|f'''\left((1-\tilde{k})\bar{g}(\bs{y})\right)\right|,\nonumber\\
    {}&= \frac{1}{3}\left|f'''\left(K\bar{g}(\bs{y})\right)\right|.
\end{align}
Substituting it into \eqref{eq:higher.order.term.bias} results in
\begin{align}
    \left|\bb{E}\left[\left(\Delta g_M(\bs{y})\right)^3\int_0^1f'''\left(\bar{g}(\bs{y})+s \Delta g_M(\bs{y})\right)(1-s)^2\di{}s\right]\right|{}&\leq \frac{1}{3}\bb{E}\left[\bb{E}\left[\left|\Delta g_M(\bs{y})\right|^3\Bigg|\bs{y}\right]\left|f'''\left(K\bar{g}(\bs{y})\right)\right|\right].
\end{align}
By a similar derivation as in~\eqref{eq:assumption.4.KH} and~\eqref{eq:Holder.finite}, it follows that
\begin{align}
    \frac{1}{3}\bb{E}\left[\bb{E}\left[\left|\Delta g_M(\bs{y})\right|^3\Bigg|\bs{y}\right]\left|f'''\left(K\bar{g}(\bs{y})\right)\right|\right]{}&\leq \frac{1}{3}\bb{E}\left[\left|f'''\left(K\bar{g}(\bs{y})\right)\right|^p\right]^{\frac{1}{p}}\bb{E}\left[\left|\bar{g}(\bs{y})\right|^{3q}\right]^{\frac{1}{q}}k^3C_{\epsilon, d_2}^3M^{-3+3\epsilon},\nonumber\\
    {}&<\infty,
\end{align}
for any $\epsilon>0$. 
For an $|f'''|$ that is monotonically increasing, similar arguments hold 
for $K\coloneqq 1+\tilde{k}$.
\end{proof}
\begin{rmk}[Dependence on the dimension $d_2$]\label{rmk:d2}
    Under appropriate regularity conditions, the dependence on the dimension $d_2$ in the above result can be neglected for certain constructions of low-discrepancy sequences (that is, lattice rules); cf.,~\cite{Kaa24}. The study of these constructions is beyond the scope of this work.
\end{rmk}
\subsection{Statistical error analysis}
The statistical error is bounded in probability by the Chebyshev inequality using the estimator variance $\bb{V}[I_{\rm{rDLQ}}]$; that is, for
\begin{equation}
    \epsilon_{\rm{rDLQ}}\coloneqq C_{\alpha}\sqrt{\bb{V}[I_{\rm{rDLQ}}]},
\end{equation}
we have
\begin{equation}
    \bb{P}\left(|I_{\rm{rDLQ}}-\bb{E}[I_{\rm{rDLQ}}]|\leq \epsilon_{\rm{rDLQ}}\right)\geq 1-\alpha
\end{equation}
for $\alpha\in(0,1)$, where $C_{\alpha}=1/\sqrt{\alpha}$. Next, the assumptions for establishing the bound on the variance of the rDLQMC estimator are presented.
\begin{asu}[Integrability assumption for the composition of the outer and inner integrand]\label{asu:f.3.diff.variance}
There exists $0<\delta_{K} <1$ and $1\leq p\leq\infty$ such that $\tilde{f}(\cdot)\equiv f(K\bar{g}_h(\cdot)):[0,1]^{d_1}\to\bb{R}$ is an element of the Sobolev space $W^{2,2p}([0,1]^{d_1})$ for all $K\in(1-\delta_{K},1+\delta_{K})$, where $\bar{g}_h$ is as in Assumption~\ref{asu:FEM} and \eqref{eq:bar.g} as $h\to 0$.
\end{asu}
\begin{asu}[Inverse inequality for the inner integrand]\label{asu:g.L4.v}
Let $g_h(\cdot,\cdot):[0,1]^{d_1}\times[0,1]^{d_2}\to\bb{R}$ and assume that there exists $0\leq k<\infty$ such that
    \begin{equation}
    \sup_{\bs{y}\in[0,1]^{d_1}} \left\lvert \frac{V_{\rm{HK}}(g_h(\bs{y},\cdot))}{ \bar{g}_h(\bs{y})} \right\rvert \leq k\;,
\end{equation}
with $V_{\rm{HK}}$ as in \eqref{eq:VHK}, $g_h$ and $\bar{g}_h$ as in Assumption~\ref{asu:FEM} and \eqref{eq:bar.g} as $h\to 0$. Furthermore, if Assumption~\ref{asu:f.3.diff.variance} holds, assume that there exists $1\leq q\leq\infty$ with $1/p+1/q=1$ for $p$ as in Assumption~\ref{asu:f.3.diff.variance} such that $\bar{g}_h\in L^{4q}([0,1]^{d_1})$.
\end{asu}
\begin{prop}[Variance of the rDLQMC estimator]\label{prop:V.DLQ}
For $\tilde{f}\equiv f(\bar{g}_h)$ satisfying Assumptions~\ref{eq:boundary.growth} and \ref{asu:f.3.diff.variance}, $g_h$ satisfying Assumption~\ref{asu:g.L4.v}, if $|f''|$ is monotonic, the variance of the rDLQMC estimator \eqref{eq:dlqmc} has the following upper bound:
\begin{equation}\label{eq:variance}
  \mathbb{V}[I_{\rm{rDLQ}}]\leq \frac{2b^2B_A^2C_{\epsilon, d_1}^2}{N^{2-2\epsilon-2\max_{i}A_i}} +\frac{2\mathbb{E}[|\bar{g}_h|^2\left|f'(\bar{g}_h)\right|^2]k^2C_{\epsilon, d_2}^2}{NM^{2-2\epsilon}}+\frac{\mathbb{E}[|\bar g_h|^4\left|f''(K\bar{g}_h)\right|^2](\Gamma_{d_1}+1)k^4C_{\epsilon, d_2}^4}{2NM^{4-4\epsilon}},
\end{equation}
for any $\epsilon>0$, where $C_{\epsilon,d_1},C_{\epsilon,d_2}\to\infty$ as $\epsilon\to 0$, $B_A\to\infty$ as $\min_{i}A_i\to 0$, and $\Gamma_{d_1}\to\infty$ as $d_1\to\infty$.
\end{prop}
\begin{proof}
The index $h$ is omitted for simplicity, and the following Taylor expansion of $f(\hat{g}_M(\bs{y}))$ around $\bar{g}(\bs{y})$ is considered:
\begin{align}\label{eq:Taylor.IM.2}
  f(\hat{g}_M(\bs{y}))&=f(\bar{g}(\bs{y})) + f'(\bar{g}(\bs{y}))\Delta g_M(\bs{y})+ \left(\Delta g_M(\bs{y})\right)^2\int_0^1f''\left(\bar{g}(\bs{y})+s \Delta g_M(\bs{y})\right)(1-s)\di{}s.
\end{align}
By the law of total variance, we obtain
\begin{align}\label{eq:total.variance}
\mathbb{V}[I_{\rm{rDLQ}}]{}&=\mathbb{V}\left[\frac{1}{N}\sum_{n=1}^Nf(\hat{g}_M(\bs{y}^{(n)}))\right],\nonumber\\
{}&=\mathbb{V}\left[\mathbb{E}\left[\frac{1}{N}\sum_{n=1}^Nf(\hat{g}_M(\bs{y}^{(n)}))\Bigg|\{\bs{y}^{(n)}\}_{n=1}^N\right]\right] + \mathbb{E}\left[\mathbb{V}\left[\frac{1}{N}\sum_{n=1}^Nf(\hat{g}_M(\bs{y}^{(n)}))\Bigg|\{\bs{y}^{(n)}\}_{n=1}^N\right]\right],
\end{align}
where conditioning occurs on the randomization $\bs{\rho}$ of the outer random variables $\bs{y}^{(n)}=\bs{\rho}(\bs{u}^{(n)})$, where $1\leq n\leq N$. Moreover, using \eqref{eq:Taylor.IM.2} for the first term in \eqref{eq:total.variance} yields the following:
\begin{align}\label{eq:ve.3}
  &{}\mathbb{V}\left[\mathbb{E}\left[\frac{1}{N}\sum_{n=1}^Nf(\hat{g}_M(\bs{y}^{(n)}))\Bigg|\{\bs{y}^{(n)}\}_{n=1}^N\right]\right],\nonumber\\
  {}&=\mathbb{V}\left[\frac{1}{N}\sum_{n=1}^N \Bigg(f(\bar{g}(\bs{y}^{(n)})) +\mathbb{E}\left[\left(\Delta g_M(\bs{y}^{(n)})\right)^2\int_0^1f''\left(\bar{g}(\bs{y}^{(n)})+s \Delta g_M(\bs{y}^{(n)})\right)(1-s)\di{}s\Bigg|\bs{y}^{(n)}\right]\Bigg)\right]\nonumber\\
  &{}\leq 2\mathbb{V}\left[\frac{1}{N}\sum_{n=1}^N f(\bar{g}(\bs{y}^{(n)})) \right]+2\mathbb{V}\left[\frac{1}{N}\sum_{n=1}^N \mathbb{E}\left[\left(\Delta g_M(\bs{y}^{(n)})\right)^2\int_0^1f''\left(\bar{g}(\bs{y}^{(n)})+s \Delta g_M(\bs{y}^{(n)})\right)(1-s)\di{}s\Bigg|\bs{y}^{(n)}\right]\right],
\end{align}
where the last line follows from the Cauchy--Schwarz inequality (see Appendix~\ref{app:variance}). For the first term in \eqref{eq:ve.3}, we have
\begin{align}\label{eq:var.hk}
    2\mathbb{V}\left[\frac{1}{N}\sum_{n=1}^N f(\bar{g}(\bs{y}^{(n)})) \right]{}&=2\mathbb{E}\left[\left|\frac{1}{N}\sum_{n=1}^N f(\bar{g}(\bs{y}^{(n)}))-\int_{[0,1]^{d_1}}f(\bar{g}(\bs{y}))\di{}\bs{y}\right|^2 \right],\nonumber\\
    {}&\leq 2 b^2 B_A^2 C_{\epsilon,d_1}^2 N^{-2+2\epsilon+2\max_iA_i},
\end{align}
by Assumption~\ref{eq:boundary.growth} and Proposition~\ref{cor:he}.
For the second term in \eqref{eq:ve.3}, the variance of scrambled points from a digital sequence is bounded by the MC variance times a constant $\Gamma_{d_1}>0$, where $\Gamma_{d_1}\to\infty$ as $d_1\to\infty$ \cite{Owe97, Pan23}; thus
\begin{align}
    {}&2\mathbb{V}\left[\frac{1}{N}\sum_{n=1}^N \mathbb{E}\left[\left(\Delta g_M(\bs{y}^{(n)})\right)^2\int_0^1f''\left(\bar{g}(\bs{y}^{(n)})+s \Delta g_M(\bs{y}^{(n)})\right)(1-s)\di{}s\Bigg|\bs{y}^{(n)}\right]\right],\nonumber\\
    {}&\leq 2\Gamma_{d_1}\mathbb{V}\left[\frac{1}{N}\sum_{n=1}^N \mathbb{E}\left[\left(\Delta g_M(\bs{z}^{(n)})\right)^2\int_0^1f''\left(\bar{g}(\bs{z}^{(n)})+s \Delta g_M(\bs{z}^{(n)})\right)(1-s)\di{}s\Bigg|\bs{z}^{(n)}\right]\right],
\end{align}
where $\bs{z}^{(n)}$, with $1\leq n\leq N$, is iid in $[0,1]^{d_{1}}$. From this, we obtain
\begin{align}\label{eq:ve.h.o.t}
{}&2\Gamma_{d_1}\mathbb{V}\left[\frac{1}{N}\sum_{n=1}^N \mathbb{E}\left[\left(\Delta g_M(\bs{z}^{(n)})\right)^2\int_0^1f''\left(\bar{g}(\bs{z}^{(n)})+s \Delta g_M(\bs{z}^{(n)})\right)(1-s)\di{}s\Bigg|\bs{z}^{(n)}\right]\right],\nonumber\\
{}&=\frac{2\Gamma_{d_1}}{N}\mathbb{V}\left[\mathbb{E}\left[\left(\Delta g_M(\bs{z})\right)^2\int_0^1f''\left(\bar{g}(\bs{z})+s \Delta g_M(\bs{z})\right)(1-s)\di{}s\Bigg|\bs{z}\right]\right],\nonumber\\
{}&\leq \frac{2\Gamma_{d_1}}{N}\mathbb{E}\left[\left|\mathbb{E}\left[\left(\Delta g_M(\bs{z})\right)^2\int_0^1f''\left(\bar{g}(\bs{z})+s \Delta g_M(\bs{z})\right)(1-s)\di{}s\Bigg|\bs{z}\right]\right|^2\right],\nonumber\\
{}&\leq \frac{2\Gamma_{d_1}}{N}\mathbb{E}\left[\mathbb{E}\left[\left|\Delta g_M(\bs{z})\right|^4\left(\int_0^1\left|f''\left(\bar{g}(\bs{z})+s \Delta g_M(\bs{z})\right)\right|(1-s)\di{}s\right)^2\Bigg|\bs{z}\right]\right].
\end{align}
With an argument similar to that in \eqref{eq:monotonicity}, by Assumption~\ref{asu:g.L4.v} and because $|f''|$ is monotonic,
\begin{equation}
    \int_0^1\left|f''\left(\bar{g}(\bs{z})+s \Delta g_M(\bs{z})\right)\right|(1-s)\di{}s\leq\frac{1}{2}\left|f''\left(K\bar{g}(\bs{z})\right)\right|,
\end{equation}
where $0< K<2$ for a sufficiently large $M$. Substituting this into \eqref{eq:ve.h.o.t} results in the following:
\begin{align}
    \frac{2\Gamma_{d_1}}{N}{}&\mathbb{E}\left[\mathbb{E}\left[\left|\Delta g_M(\bs{z})\right|^4\left(\int_0^1\left|f''\left(\bar{g}(\bs{z})+s \Delta g_M(\bs{z})\right)\right|(1-s)\di{}s\right)^2\Bigg|\bs{z}\right]\right],\nonumber\\
    {}&\leq \frac{\Gamma_{d_1}}{2N}\mathbb{E}\left[\mathbb{E}\left[\left|\Delta g_M(\bs{z})\right|^4\Big|\bs{z}\right]\left|f''\left(K\bar{g}(\bs{z})\right)\right|^2\right].
\end{align}
By a similar derivation as in~\eqref{eq:assumption.4.KH} and~\eqref{eq:Holder.finite}, it follows that
\begin{align}
    \frac{\Gamma_{d_1}}{2N}\mathbb{E}\left[\mathbb{E}\left[\left|\Delta g_M(\bs{z})\right|^4\Big|\bs{z}\right]\left|f''\left(K\bar{g}(\bs{z})\right)\right|^2\right]{}&\leq \frac{\Gamma_{d_1}}{2N}\mathbb{E}\left[|\bar{g}(\bs{z})|^{4q}\right]^{\frac{1}{q}}\mathbb{E}\left[\left|f''\left(K\bar{g}(\bs{z})\right)\right|^{2p}\right]^{\frac{1}{p}}k^4C_{\epsilon,d_2}^4M^{-4+4\epsilon},\nonumber\\
    {}&<\infty,
\end{align}
for any $\epsilon>0$ by Assumption~\ref{asu:f.3.diff.variance} and \ref{asu:g.L4.v}.
For the second term in \eqref{eq:total.variance}, using \eqref{eq:Taylor.IM.2} yields the following:
\begin{align}\label{eq:ev.3}
  &{}\mathbb{E}\left[\mathbb{V}\left[\frac{1}{N}\sum_{n=1}^Nf(\hat{g}_M(\bs{y}^{(n)}))\Bigg|\{\bs{y}^{(n)}\}_{n=1}^N\right]\right]=\frac{1}{N^2}\sum_{n=1}^N\mathbb{E}\left[\mathbb{V}\left[f(\hat{g}_M(\bs{y}^{(n)}))\Bigg|\bs{y}^{(n)}\right]\right],\nonumber\\
  {}&=\frac{1}{N}\mathbb{E}\left[\mathbb{V}\left[f(\hat{g}_M(\bs{y}))\Bigg|\bs{y}\right]\right],\nonumber\\
  &{}=\frac{1}{N}\mathbb{E}\left[\mathbb{V}\left[f'(\bar{g}(\bs{y}))\Delta g_M(\bs{y})+ \left(\Delta g_M(\bs{y})\right)^2\int_0^1f''\left(\bar{g}(\bs{y})+s \Delta g_M(\bs{y})\right)(1-s)\di{}s\Bigg|\bs{y}\right]\right],\nonumber\\
  {}&\leq\frac{2}{N}\mathbb{E}\left[\mathbb{V}\left[f'(\bar{g}(\bs{y}))\Delta g_M(\bs{y})\Bigg|\bs{y}\right]\right]+ \frac{2}{N}\mathbb{E}\left[\mathbb{V}\left[\left(\Delta g_M(\bs{y})\right)^2\int_0^1f''\left(\bar{g}(\bs{y})+s \Delta g_M(\bs{y})\right)(1-s)\di{}s\Bigg|\bs{y}\right]\right],
\end{align}
because the inner randomizations are independent. For the first term in \eqref{eq:ev.3},
\begin{align}
    \frac{2}{N}\mathbb{E}\left[\mathbb{V}\left[f'(\bar{g}(\bs{y}))\Delta g_M(\bs{y})\Bigg|\bs{y}\right]\right]{}&=\frac{2}{N}\mathbb{E}\left[\left(f'(\bar{g}(\bs{y}))\right)^2\mathbb{V}\left[\Delta g_M(\bs{y})\Bigg|\bs{y}\right]\right],\nonumber\\
    {}&\leq \frac{2}{N}\mathbb{E}\left[\left(f'(\bar{g}(\bs{y}))\right)^2\mathbb{E}\left[\left|\Delta g_M(\bs{y})\right|^2\Bigg|\bs{y}\right]\right],
\end{align}
and by Assumptions~\ref{asu:f.3.diff.variance} and \ref{asu:g.L4.v}, it follows again that
\begin{align}\label{eq:var:HK.2}
    \frac{2}{N}\mathbb{E}\left[\left(f'(\bar{g}(\bs{y}))\right)^2\mathbb{E}\left[\left|\Delta g_M(\bs{y})\right|^2\Bigg|\bs{y}\right]\right]{}&\leq \frac{2}{N}\mathbb{E}\left[\left(f'(\bar{g}(\bs{y}))\right)^{2p}\right]^{\frac{1}{p}}\bb{E}\left[|\bar{g}(\bs{y})|^{2q}\right]^{\frac{1}{q}}k^2C_{\epsilon,d_2}^2M^{-2+2\epsilon},\nonumber\\
    {}&<\infty,
\end{align}
for any $\epsilon>0$.
For the second term in \eqref{eq:ev.3}, we have
\begin{align}
    \frac{2}{N}{}&\mathbb{E}\left[\mathbb{V}\left[\left(\Delta g_M(\bs{y})\right)^2\int_0^1f''\left(\bar{g}(\bs{y})+s \Delta g_M(\bs{y})\right)(1-s)\di{}s\Bigg|\bs{y}\right]\right],\nonumber\\
    {}&\leq\frac{2}{N}\mathbb{E}\left[\mathbb{E}\left[\left|\left(\Delta g_M(\bs{y})\right)^2\int_0^1f''\left(\bar{g}(\bs{y})+s \Delta g_M(\bs{y})\right)(1-s)\di{}s\right|^2\Bigg|\bs{y}\right]\right],\nonumber\\
    {}&\leq\frac{2}{N}\mathbb{E}\left[\mathbb{E}\left[\left|\Delta g_M(\bs{y})\right|^4\left(\int_0^1\left|f''\left(\bar{g}(\bs{y})+s \Delta g_M(\bs{y})\right)\right|(1-s)\di{}s\right)^2\Bigg|\bs{y}\right]\right],\nonumber\\
    {}&\leq\frac{1}{2N}\mathbb{E}\left[\mathbb{E}\left[\left|\Delta g_M(\bs{y})\right|^4\Bigg|\bs{y}\right]\left|f''\left(K\bar{g}(\bs{y})\right)\right|^2\right],\nonumber\\
    {}&\leq\frac{1}{2N}\mathbb{E}\left[\left|\bar g(\bs{y})\right|^{4q}\right]^{\frac{1}{q}}\bb{E}\left[\left|f''\left(K\bar{g}(\bs{y})\right)\right|^{2p}\right]^{\frac{1}{p}}k^4C_{\epsilon,d_2}^4M^{-4+4\epsilon},\nonumber\\
    {}&<\infty,
\end{align}
for any $\epsilon>0$ by Assumptions~\ref{asu:f.3.diff.variance} and \ref{asu:g.L4.v}. 
\end{proof}
\begin{rmk}[rQMC analysis for specific problems]
    Detailed analyses of the performance of rQMC estimators for integrands based on specific finite element formulations are available (e.g., \cite{Kuo12, Kuo15, Gra15}). However, for the generality of the bias and variance bounds of the rQMC method, the assumptions are restricted to the differentiability and integrability of the integrands. Moreover, the dependence on the dimension $d_2$ can be neglected under appropriate regularity conditions (see Remark~\ref{rmk:d2}).
\end{rmk}
\begin{rmk}[Stability with respect to the discretized inner integrand]
    The quantities appearing in Propositions~\ref{prop:B.DLQ} and~\ref{prop:V.DLQ} depend on the discretized inner integrand $g_h$. We assume that the relevant moments of $\bar g_h$ and of the compositions involving $f$ remain finite as $h \to 0$. Since a detailed analysis of the discretization error is beyond the scope of the present work, Assumption~\ref{asu:FEM} is used to quantify its effect on the resulting estimator.
\end{rmk}

\subsection{Optimal work allocation}
Using Proposition~\ref{prop:B.DLQ} and~\ref{prop:V.DLQ}, we analyzed the work required for the rDLQMC estimator in terms of the number of samples and the discretization parameter of the inner integrand to achieve a prescribed tolerance. First, we ensured that the rDLQMC error \eqref{eq:rDLQMC.error} is smaller than some tolerance $TOL>0$; that is,
\begin{align}\label{eq:total.error}
   |I_{\rm{rDLQ}}-I|{}&\leq |\mathbb{E}[I_{\rm{rDLQ}}]-I|+|I_{\rm{rDLQ}}-\mathbb{E}[I_{\rm{rDLQ}}]|,\nonumber\\
   {}&\leq |\mathbb{E}[I_{\rm{rDLQ}}]-I|+C_{\alpha}\sqrt{\bb{V}[I_{\rm{rDLQ}}]}, \quad \text{with probability $1-\alpha$}\nonumber\\
   {}&\leq TOL.
\end{align}
An error-splitting parameter $\kappa\in(0,1)$ is introduced to control the bias and statistical error separately to obtain the following:
\begin{equation}
    |\mathbb{E}[I_{\rm{rDLQ}}]-I|\leq (1-\kappa)TOL,
\end{equation}
and
\begin{equation}
    C_{\alpha}\sqrt{\bb{V}[I_{\rm{rDLQ}}]}\leq \kappa TOL,
\end{equation}
which is written as
\begin{equation}
    \bb{V}[I_{\rm{rDLQ}}]\leq \left(\frac{\kappa TOL}{C_{\alpha}}\right)^2.
\end{equation}
Optimizing the error-splitting parameter $\kappa$ enables balancing the number of inner and outer samples to control the total error \eqref{eq:total.error} efficiently.
Substituting the error bounds established in Propositions~\ref{prop:B.DLQ} and~\ref{prop:V.DLQ} yields the following sufficient conditions to meet the error criterion \eqref{eq:total.error}:
\begin{equation}\label{eq:Bias.Constraint}
    C_{\rm{disc}}h^{\eta}+\frac{C_{\rm{Q},\epsilon}^{(3)}}{M^{2-2\epsilon}}\leq(1-\kappa)TOL,
\end{equation}
and
\begin{equation}\label{eq:Variance.Constraint}
    \frac{C_{\rm{Q},\epsilon}^{(1)}}{N^{2-2\epsilon-2\max_{i}A_i}}+\frac{C_{\rm{Q},\epsilon}^{(2)}}{NM^{2-2\epsilon}}\leq\left(\frac{\kappa TOL}{C_{\alpha}}\right)^2,
\end{equation}
where the constant terms in the bounds for the bias and statistical error (\eqref{eq:Bias.constraint} and \eqref{eq:variance}) are combined and higher-order terms are neglected, where $C_{\rm{Q},\epsilon}^{(1)},C_{\rm{Q},\epsilon}^{(2)},C_{\rm{Q},\epsilon}^{(3)}\to\infty$ as $\epsilon\to 0$ and $C_{\rm{Q},\epsilon}^{(1)}\to\infty$ as $\max_{i}A_i\to 0$.

\begin{prop}[Optimal work of the rDLQMC estimator]\label{prop:W.DLQ}
Given Assumptions~\ref{eq:boundary.growth} to~\ref{asu:g.L4.v}, assuming that $|f''|$ and $|f'''|$ are monotonic, the total work of the optimized rDLQMC estimator \eqref{eq:dlqmc} for a specified error tolerance $TOL>0$ is given by
\begin{equation}\label{eq:work.dlq}
    W_{\rm{rDLQ}}^{\ast}=\cl{O}\left(TOL^{-
    \left(
        \frac{2}{2-2\epsilon-2\max_{i}A_i}+\frac{1}{2-2\epsilon}+\frac{\gamma}{\eta}
    \right) 
    }\right),
\end{equation}
for any $\epsilon>0$
as $TOL\to 0$, where the constant implied by the $\cl{O}$ notation approaches infinity as $\epsilon\to 0$.
\end{prop}

\begin{proof}
The computational work for the rDLQMC estimator with $R=S=1$ is
\begin{equation}\label{eq:Work}
  W_{\mathrm{rDLQ}}=\cl{O}\left(NMh^{-\gamma}\right).
\end{equation}
We obtain the optimal rDLQMC estimator setting by minimizing the computational work \eqref{eq:Work} subject to the constraints~\eqref{eq:Bias.Constraint} and \eqref{eq:Variance.Constraint}; that is,
\begin{equation}\label{eq:Lagrangian}
  (N^{\ast},M^{\ast},h^{\ast},\kappa^{\ast})=\argmin_{(N,M,h,\kappa)}NMh^{-\gamma} \quad \text{subject to}\quad \begin{cases}\frac{C_{\rm{Q},\epsilon}^{(1)}}{N^{2-2\epsilon-2\max_{i}A_i}}+\frac{C_{\rm{Q},\epsilon}^{(2)}}{NM^{2-2\epsilon}}\leq\left(\frac{\kappa TOL}{C_{\alpha}}\right)^2\\C_{\rm{disc}}h^{\eta}+\frac{C_{\rm{Q},\epsilon}^{(3)}}{M^{2-2\epsilon}}\leq(1-\kappa)TOL\end{cases}.
\end{equation}
This problem is solved using Lagrange multipliers to derive the optimal $M^{\ast}$ and $h^{\ast}$ in terms of $\kappa$ and $N$. The equation for $\kappa^{\ast}$ is cubic and has a closed-form solution, but it is unwieldy to state explicitly. The last remaining equation for $N^{\ast}$ has no closed-form solution; therefore, a simplified version must be solved, and the convergence of the resulting solution to the true solution as $TOL\to 0$ must be demonstrated. The optimal values for $M$ and $h$ satisfy the following:
\begin{equation}\label{eq:M.opt}
M^{\ast}=\left(\frac{C_{\rm{Q},\epsilon}^{(2)}}{N\left(\frac{\kappa TOL}{C_{\alpha}}\right)^2-\frac{C_{\rm{Q},\epsilon}^{(1)}}{N^{1-2\epsilon-2\max_{i}A_i}}}\right)^{\frac{1}{2-2\epsilon}},
\end{equation}
and
\begin{equation}\label{eq:h.opt}
h^{\ast}=\left(\frac{(1-\kappa)TOL-\frac{C_{\rm{Q},\epsilon}^{(3)}}{C_{\rm{Q},\epsilon}^{(2)}}\left(N\left(\frac{\kappa TOL}{C_{\alpha}}\right)^2-\frac{C_{\rm{Q},\epsilon}^{(1)}}{N^{1-2\epsilon-2\max_{i}A_i}}\right)}{C_{\rm{disc}}}\right)^{1/\eta}.
\end{equation}
Moreover, the optimal $\kappa^{\ast}$ is given by the real root of
\begin{align}
  &\left[\frac{C_{\rm{Q},\epsilon}^{(3)}}{C_{\rm{Q},\epsilon}^{(2)}}\left(\frac{NTOL}{C_{\alpha}}\right)^2\left(\eta + \gamma(2-2\epsilon)\right)\right]\kappa^{\ast 3} + \left[NTOL(\eta+\frac{\gamma(2-2\epsilon)}{2})\right]\kappa^{\ast 2}\nonumber\\
  -&\left[N\left(\eta TOL+\frac{C_{\rm{Q},\epsilon}^{(3)}}{C_{\rm{Q},\epsilon}^{(2)}}\frac{C_{\rm{Q},\epsilon}^{(1)}}{N^{1-2\epsilon-2\max_{i}A_i}}(\eta+\gamma(2-2\epsilon))\right)\right]\kappa^{\ast}-\frac{\gamma(2-2\epsilon)C_{\rm{Q},\epsilon}^{(1)}C_{\alpha}^2}{2N^{1-2\epsilon-2\max_{i}A_i}TOL}=0.
\end{align}
Finally, the solution to the following yields the optimal $N^{\ast}$:
\begin{align}\label{eq:N.num}
  0={}&\left[\frac{C_{\rm{Q},\epsilon}^{(3)}}{C_{\rm{Q},\epsilon}^{(2)}}\left(\frac{\kappa TOL}{C_{\alpha}}\right)^2\right]N^{3-2\epsilon-2\max_{i}A_i}-\left[TOL\left(1-\kappa\left(1+\frac{\gamma}{2\eta}\right)\right)\right]N^{2-2\epsilon-2\max_{i}A_i}-\frac{C_{\rm{Q},\epsilon}^{(3)}}{C_{\rm{Q},\epsilon}^{(2)}}C_{\rm{Q},\epsilon}^{(1)}N\nonumber\\
  {}&+\frac{\gamma C_{\alpha}^2C_{\rm{Q},\epsilon}^{(1)}(1-2\epsilon-2\max_{i}A_i)}{2\eta\kappa TOL}.
\end{align}
It immediately follows from \eqref{eq:Lagrangian} that such an $N^{\ast}$ exists. The optimal discretization parameter $h^{\ast}$ is assumed to be independent of the sampling method (e.g., MC or rQMC). The bias constraint \eqref{eq:Bias.Constraint} is split as follows:
\begin{align}
    C_{\rm{disc}}h^{\ast \eta}&\leq\frac{1}{2}(1-\kappa)TOL,\label{eq:bias.split.1}\\
    \frac{C_{\rm{Q},\epsilon}^{(3)}}{M^{\ast2-2\epsilon}}&\leq\frac{1}{2}(1-\kappa)TOL,\label{eq:bias.split.2}
\end{align}
to derive the asymptotic rates in terms of the tolerance ($TOL$), implying that
\begin{equation}\label{eq:h.rate}
h^{\ast}\propto TOL^{\frac{1}{\eta}},
\end{equation}
and
\begin{equation}\label{eq:M.rate}
M^{\ast}\propto TOL^{-\frac{1}{2-2\epsilon}}.
\end{equation}
Next, we demonstrate that $N\propto TOL^{-\frac{2}{2-2\epsilon-2\max_{i}A_i}}$. From the variance constraint \eqref{eq:Variance.Constraint}, we obtain the following:
\begin{equation}\label{eq:appr.variance}
    \left( \frac{\kappa TOL}{C_{\alpha}} \right)^2N^{2-2\epsilon-2\max_{i}A_i} - (1 - \kappa) TOL N^{1-2\epsilon-2\max_{i}A_i} = C_{\rm{Q},\epsilon}^{(1)}.
\end{equation}
The term on the right-hand side of \eqref{eq:appr.variance} is constant in $TOL$.
The equation in \eqref{eq:appr.variance} can be solved by ignoring the second term on the left-hand side, yielding the following approximate solution:
\begin{equation}\label{eq:app.sol.N}
    N\approx \left(\frac{C_{\alpha}^2C_{\rm{Q},\epsilon}^{(1)}}{\kappa^2}\right)^{\frac{1}{2-2\epsilon-2\max_{i}A_i}} TOL^{-\frac{2}{2-2\epsilon-2\max_{i}A_i}}.
\end{equation}
We assess whether the ignored term in \eqref{eq:appr.variance} approaches zero as $TOL\to 0$ by inserting \eqref{eq:app.sol.N} to determine whether this approximation converges to the true solution.
For this term, we have
\begin{align}
(1 - \kappa) TOL N^{1-2\epsilon-2\max_{i}A_i}\approx{}&(1-\kappa)\left(\frac{C_{\alpha}^2C_{\rm{Q},\epsilon}^{(1)}}{\kappa^2}\right)^{\frac{1-2\epsilon-2\max_{i}A_i}{2-2\epsilon-2\max_{i}A_i}}TOL^{1-\frac{2-4\epsilon-4\max_{i}A_i}{2-2\epsilon-2\max_{i}A_i}},\nonumber\\
={}&(1-\kappa)\left(\frac{C_{\alpha}^2C_{\rm{Q},\epsilon}^{(1)}}{\kappa^2}\right)^{\frac{1-2\epsilon-2\max_{i}A_i}{2-2\epsilon-2\max_{i}A_i}}TOL^{\frac{2\epsilon+2\max_{i}A_i}{2-2\epsilon-2\max_{i}A_i}},
\end{align}
where, for the exponent of $TOL$, $0<(2\epsilon+2\max_{i}A_i)/(2-2\epsilon-2\max_{i}A_i)<1$ as $\epsilon\to 0$ because of the assumption that $\max_{i}A_i< 1/2$. Thus, this term approaches zero as $TOL\to 0$. In contrast, ignoring the first term in \eqref{eq:appr.variance} results in the approximate solution $N\approx-(C_{\rm{Q},\epsilon}^{(1)}/(1-\kappa))^{1/(1-2\epsilon-2\max_{i}A_i)}TOL^{-1/(1-2\epsilon-2\max_{i}A_i)}$. With this solution, the first term in \eqref{eq:appr.variance} has an exponent of $TOL$ which is strictly negative; thus, it approaches negative infinity as $TOL\to 0$.
\end{proof}

The constants $C_{\rm{Q},\epsilon}^{(1)}$, $C_{\rm{Q},\epsilon}^{(2)}$, and $C_{\rm{Q},\epsilon}^{(3)}$ introduced in \eqref{eq:Bias.Constraint} and \eqref{eq:Variance.Constraint} can be estimated using $S,R>1$ randomizations. Rather than using the approximate solution \eqref{eq:app.sol.N}, we can also solve the equation in \eqref{eq:N.num} numerically. Splitting the bias constraint \eqref{eq:Bias.Constraint} more elaborately than in \eqref{eq:bias.split.1} and \eqref{eq:bias.split.2} (i.e., by a splitting parameter possibly different from $1/2$) might improve the near-optimality of the approximation \eqref{eq:app.sol.N}.

The above result in Proposition~\ref{prop:W.DLQ} is valid for all constructions of rQMC estimators, that is, combinations of low-discrepancy sequences and randomization methods. However, given certain additional assumptions of the integrands, rQMC estimators using Sobol' points with Owen's scrambling can achieve even faster convergence, and thus a lower computational cost. In particular, for integrands $\varphi$ with bounded generalized Hardy--Krause variation of order one, denoted as $V_{\rm{GHK}}$  (see~\cite[Definition 6.24]{Dic13}), the variance of an rQMC estimator with $R=1$ randomizations and $M$ Sobol' points randomized via Owen's scrambling is bounded as follows (see~\cite[Theorem 6.25]{Dic13}):
\begin{equation}\label{prop.BGHK.var}
    \bb{V}[I_{\rm{rQ}}]\leq \left(V_{\rm{GHK}}(\varphi)\right)^2C_{\epsilon, d}^2M^{-3+2\epsilon},
\end{equation}
where $C_{\epsilon, d}\to\infty$ as $\epsilon\to0$. We present the total computational work for the rDLQMC estimator under this setting below.
\begin{cor}[Optimal work for scrambled Sobol' sequence]\label{cor:optimal.work.s}
Let $f\circ \bar{g}$ have bounded generalized Hardy--Krause variation $V_{\rm{GHK}}$ of order one (see~\cite[Definition 6.24]{Dic13}). Moreover, let  the assumptions of Proposition~\ref{prop:W.DLQ} hold with $V_{\rm{HK}}$ replaced by $V_{\rm{GHK}}$; then, the total work of the optimized rDLQMC estimator \eqref{eq:dlqmc} based on the Sobol' sequence with Owen's scrambling for a specified error tolerance $TOL>0$ is given by
\begin{equation}\label{eq:work.dlq.scr}
    W^{\ast}_{\rm{rDLQ}}=\cl{O}\left(TOL^{-\left(\frac{2}{3-2\epsilon}+\frac{4-4\epsilon}{(3-2\epsilon)^2}+\frac{\gamma}{\eta}\right)}\right),
\end{equation}
for any $\epsilon>0$
as $TOL\to 0$, where the constant implied by the $\cl{O}$ notation approaches infinity as $\epsilon\to 0$.
\end{cor}
Before we proceed with the proof of Corollary~\ref{cor:optimal.work.s}, we provide some context on the results above. As $\epsilon$ approaches zero, the rate of increase of the optimized total work of the rDLQMC estimator reduces from almost $3/2+\gamma/\eta$ (Proposition~\ref{prop:W.DLQ}) to almost $10/9+\gamma/\eta$ (Corollary~\ref{cor:optimal.work.s}) when using the Sobol' sequence and Owen's scrambling for integrands with bounded generalized Hardy--Krause variation of order one. This rate is achieved for the polynomial example presented below where both the inner and outer integrands trivially satisfy this strict regularity assumption. We demonstrate in Appendix~\ref{app:GHK} how a truncation of the observation noise in EIG estimation results in bounded generalized Hardy--Krause variation of order one even for the outer integrand. However, the resulting bound grows as the region of truncation approaches the boundary, thus negating the improved convergence rate. Applying the truncation nonetheless resulted in an observed outer variance which was reduced by a constant factor compared to not using the truncation (see Figures~\ref{fig:ex0.5.pilot}-\ref{fig:ex0.5.evt}). The inner integrand in the EIG setting is sufficiently regular and the resulting optimal rate is thus almost $4/3+\gamma/\eta$ as a direct consequence of Corollary~\ref{cor:optimal.work.s}. This is because the inner variance and bias converge at rate almost $3$ compared to almost rate $2$ (for other constructions of rQMC estimators) or rate $1$ for MC estimators.
\begin{proof}[Proof of Corollary~\ref{cor:optimal.work.s}]
    For the bound on the bias error in~\eqref{eq:Holder.finite} in the proof of Proposition~\ref{prop:B.DLQ}, it follows from~\cite[Theorem 6.25]{Dic13} that 

\begin{align}
    \left|\bb{E}\left[f''(\bar{g}(\bs{y}))\left(\Delta g_M(\bs{y})\right)^2\right]\right|
    {}&\leq \bb{E}\left[\left|f''(\bar{g}(\bs{y}))\right|(\bar{g}(\bs{y}))^2\right] k^2C_{\epsilon, d_2}^2M^{-3+2\epsilon},\nonumber\\
    {}&\leq \bb{E}\left[\left|f''(\bar{g}(\bs{y}))\right|^p\right]^{\frac{1}{p}}\bb{E}\left[|\bar{g}(\bs{y})|^{2q}\right]^{\frac{1}{q}} k^2C_{\epsilon, d_2}^2M^{-3+2\epsilon},\nonumber\\
    {}&<\infty,
\end{align}
for any $\epsilon>0$, where $C_{\epsilon, d_2}\to\infty$ as $\epsilon\to0$. Similarly, for the error bound in~\eqref{eq:var.hk} in Proposition~\ref{prop:V.DLQ}, it follows that
\begin{align}\label{eq:variance.ghk}
    \mathbb{V}\left[\frac{1}{N}\sum_{n=1}^N f(\bar{g}(\bs{y}^{(n)})) \right]{}&\leq (V_{\rm{GHK}}( f(\bar{g})))^2C_{\epsilon,d_1}^2 N^{-3+2\epsilon},\nonumber\\
    {}&<\infty,
\end{align}
for any $\epsilon>0$, where $C_{\epsilon, d_1}\to\infty$ as $\epsilon\to0$, and for the error bound in~\eqref{eq:var:HK.2}, it follows that
\begin{align}\label{eq:variance.ghk.2}
    \frac{2}{N}\mathbb{E}\left[\left(f'(\bar{g}(\bs{y}))\right)^2\mathbb{E}\left[\left|\Delta g_M(\bs{y})\right|^2\Bigg|\bs{y}\right]\right]{}&\leq \frac{2}{N}\mathbb{E}\left[\left(f'(\bar{g}(\bs{y}))\right)^{2p}\right]^{\frac{1}{p}}\bb{E}\left[|\bar{g}(\bs{y})|^{2q}\right]^{\frac{1}{q}}k^2C_{\epsilon,d_2}^2M^{-3+2\epsilon},\nonumber\\
    {}&<\infty,
\end{align}
for any $\epsilon>0$, where $C_{\epsilon, d_2}\to\infty$ as $\epsilon\to0$. The higher-order terms in the proofs of Propositions~\ref{prop:B.DLQ} and~\ref{prop:V.DLQ} are still of higher or comparable order even under the original milder assumptions. It remains to solve the updated optimization problem
\begin{equation}\label{eq:Lagrangian.ghk}
  (N^{\ast},M^{\ast},h^{\ast},\kappa^{\ast})=\argmin_{(N,M,h,\kappa)}NMh^{-\gamma} \quad \text{subject to}\quad \begin{cases}\frac{C_{\rm{Q},\epsilon}^{(1)}}{N^{3-2\epsilon}}+\frac{C_{\rm{Q},\epsilon}^{(2)}}{NM^{3-2\epsilon}}\leq\left(\frac{\kappa TOL}{C_{\alpha}}\right)^2\\C_{\rm{disc}}h^{\eta}+\frac{C_{\rm{Q},\epsilon}^{(3)}}{M^{3-2\epsilon}}\leq(1-\kappa)TOL\end{cases}.
\end{equation}
Retracing the steps in the proof of Proposition~\ref{prop:W.DLQ}, we find that the following choice of optimal samples
\begin{equation}
    N^{\ast}\propto TOL^{-\frac{2}{3-2\epsilon}},
\end{equation}
\begin{equation}
    M^{\ast}\propto TOL^{-\frac{4-4\epsilon}{(3-2\epsilon)^2}},
\end{equation}
and
\begin{equation}
    h^{\ast}\propto TOL^{\frac{1}{\eta}},
\end{equation}
yields optimal total work of
\begin{equation}
    W^{\ast}_{\rm{rDLQ}}=\cl{O}\left(TOL^{-\left(\frac{2}{3-2\epsilon}+\frac{4-4\epsilon}{(3-2\epsilon)^2}+\frac{\gamma}{\eta}\right)}\right),
\end{equation}
for any $\epsilon\to 0$ as $TOL\to 0$, where the constant implied by the $\cl{O}$ notation approaches infinity as $\epsilon\to 0$. This is demonstrated by separating the two constraints in~\eqref{eq:Lagrangian.ghk} into four constraints as follows:
\begin{equation}\label{eq:ghk.N}
    \frac{C_{\rm{Q},\epsilon}^{(1)}}{N^{3-2\epsilon}}\leq\frac{1}{2}\left(\frac{\kappa TOL}{C_{\alpha}}\right)^2,
\end{equation}
\begin{equation}\label{eq:ghk.N.M}
    \frac{C_{\rm{Q},\epsilon}^{(2)}}{NM^{3-2\epsilon}}\leq\frac{1}{2}\left(\frac{\kappa TOL}{C_{\alpha}}\right)^2,
\end{equation}
\begin{equation}\label{eq:ghk.h}
    C_{\rm{disc}}h^{\eta}\leq\frac{1}{2}(1-\kappa)TOL,
\end{equation}
and
\begin{equation}\label{eq:ghk.M}
    \frac{C_{\rm{Q},\epsilon}^{(3)}}{M^{3-2\epsilon}}\leq\frac{1}{2}(1-\kappa)TOL.
\end{equation}
It follows from~\eqref{eq:ghk.h} that $h^{\ast}\propto TOL^{\frac{1}{\eta}}$. For the optimal $N^{\ast}$ and $M^{\ast}$, we make the following ansatz:
\begin{equation}
    N^{\ast}\propto TOL^{-a},
\end{equation}
and
\begin{equation}
    M^{\ast}\propto TOL^{-b}.
\end{equation}
Thus, we obtain the optimization problem
\begin{equation}\label{eq:ghk.Lagrangian}
  \min a+b \quad \text{subject to}\quad \begin{cases}\frac{2}{3-2\epsilon}&\leq a\\\frac{1}{3-2\epsilon}&\leq b\\2&\leq a+(3-2\epsilon)b
  \end{cases},
\end{equation}
with the solution $a=2/(3-2\epsilon)$ and $b=(4-4\epsilon)/(3-2\epsilon)^2$.
\end{proof}

\subsection{Pilot-based practical calibration}
The rQMC rates in \eqref{eq:Bias.Constraint} and \eqref{eq:Variance.Constraint} are typically not observed in practice for finite $N$ and $M$. Furthermore, obtaining closed-form expressions for the constants $C_{\rm{Q},\epsilon}^{(1)}$, $C_{\rm{Q},\epsilon}^{(2)}$, and $C_{\rm{Q},\epsilon}^{(3)}$ is challenging; hence, the constraints on the bias and statistical error are restated as follows:
\begin{equation}\label{eq:Bias.Constraint.num}
    \frac{C_{\rm{Q}}^{(3)}}{M^{1+\delta}}\leq(1-\kappa)TOL,
\end{equation}
and
\begin{equation}\label{eq:Variance.Constraint.num}
    \frac{C_{\rm{Q}}^{(1)}}{N^{1+\beta}}+\frac{C_{\rm{Q}}^{(2)}}{NM^{1+\delta}}\leq\left(\frac{\kappa TOL}{C_{\alpha}}\right)^2,
\end{equation}
respectively, where $0<\beta,\delta<2$, correspond to the observed improvement of the rQMC convergence over the MC convergence and the dependence on $\epsilon$ was dropped to lighten the notation. Using pilot runs, we estimate the bias bound $\frac{1}{2}|\bb{E}[f''(\bar{g}(\bs{y}))\left(\Delta g_M(\bs{y})\right)^2]|$ of the rDLQMC estimator as a function of the inner number of samples $M$ and fit the constant $C_{\rm{Q}}^{(3)}$ and rate $1+\delta$. Similarly, we estimate the outer variance $\mathbb{V}[\mathbb{E}[\frac{1}{N}\sum_{n=1}^Nf(\hat{g}_M(\bs{y}^{(n)}))|\{\bs{y}^{(n)}\}_{n=1}^N]]$ and inner variance $\mathbb{E}[\mathbb{V}[\frac{1}{N}\sum_{n=1}^Nf(\hat{g}_M(\bs{y}^{(n)}))|\{\bs{y}^{(n)}\}_{n=1}^N]]$ and fit the constants $C_{\rm{Q}}^{(1)}$ and rate $1+\beta$ as well as the constant $C_{\rm{Q}}^{(2)}$, respectively. The algorithm for running the pilot to estimate the rates and constants is provided in Appendix~\ref{app:pilot}. The total work for the pilot is
\begin{equation}
W_{\rm{pilot}}\propto W_{\rm{pilot}}^{\rm{out}}+W_{\rm{pilot}}^{\rm{in}},
\end{equation}
where
\begin{equation}
      W_{\rm{pilot}}^{\rm{out}}
  =S_{\rm{out}}N_{\rm{out}}R_{\rm{out}}M_{\rm{out}}\,W_g,
  \end{equation}
  and
  \begin{equation}
  W_{\rm{pilot}}^{\rm{in}}
  =S_{\rm{in}}N_{\rm{in}}R_{\rm{in}}M_{\rm{in}}\, W_g,
\end{equation}
and $W_g$ denotes the cost of evaluating the inner integrand $g$. The outer pilot uses many outer randomizations $S_{\rm{out}}$ and low-discrepancy points $N_{\rm{out}}$ and very few inner samples to estimate $C_{\rm{Q}}^{(1)}$ and $1+\beta$. The inner pilot uses very few outer samples and many inner randomizations $R_{\rm{in}}$ and low-discrepancy points $M_{\rm{in}}$ to estimate $C_{\rm{Q}}^{(2)}$, $C_{\rm{Q}}^{(3)}$, and $1+\delta$.
We applied the confidence constant $C_{\alpha} = \Phi^{-1}(1-\alpha/2)$ implied by the CLT in the numerical examples. A more cautious approach could be taken, where $C_{\alpha}=1/\sqrt{\alpha}$, as implied by the Chebyshev inequality. The investigation of the applicability of the CLT for rQMC is the subject of current research \cite{Lec23}.
\begin{rmk}[Effects of inner randomization on the outer rQMC convergence]
    The effectiveness of the rQMC estimators compared to the MC estimators relies on the specific low-discrepancy structure of the samples used for the evaluation of the integrands. Introducing randomness in these samples jeopardizes this structure. Nesting rQMC estimators with independent randomizations of the inner samples for each outer sample effectively transforms the outer samples into MC samples. However, the proof of Proposition~\ref{prop:W.DLQ} demonstrates that the effectiveness of the outer rQMC estimation can be preserved asymptotically if the number of outer and inner samples is appropriately increased. To the best of our knowledge, this insight is novel and constitutes a significant deviation from the analysis of other estimators discussed in the literature.
\end{rmk}
\begin{rmk}[Rounding optimal parameters for practical purposes]
Because of the properties of low-discrepancy sequences, rounding the optimal number of samples $N^{\ast}$ and $M^{\ast}$ up to the nearest power of two is beneficial. 
\end{rmk}



\subsection{Example 1: Polynomial example with exact sampling}\label{sec:pol.example}
\begin{figure}[ht]
	\subfloat[Outer variance]{%
		\includegraphics[width=0.45\textwidth]{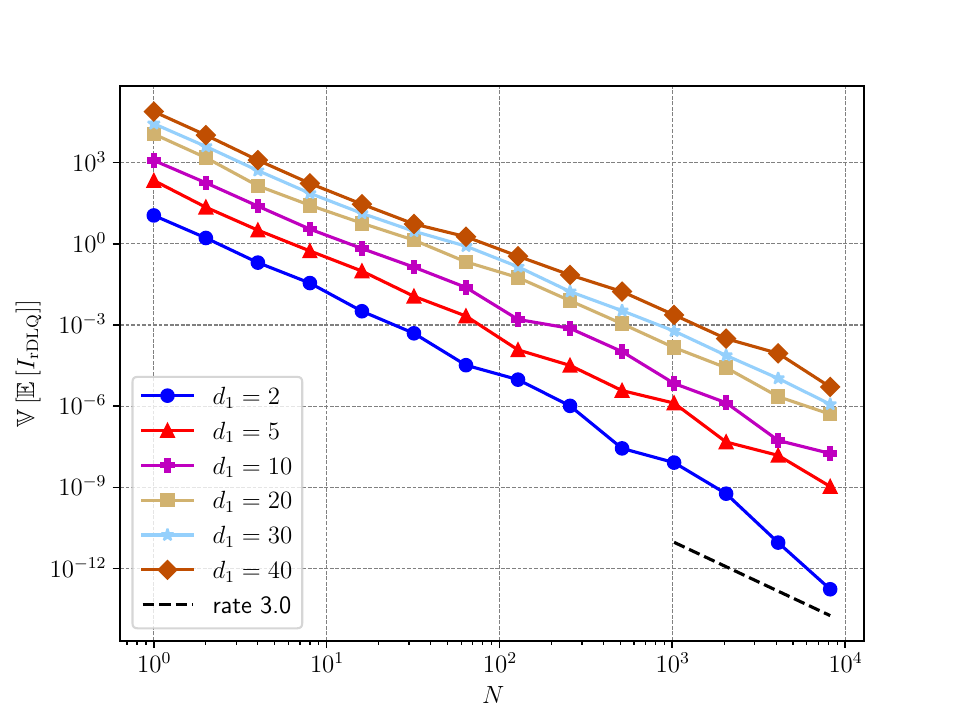}
	}
    \subfloat[Inner variance]{%
		\includegraphics[width=0.45\textwidth]{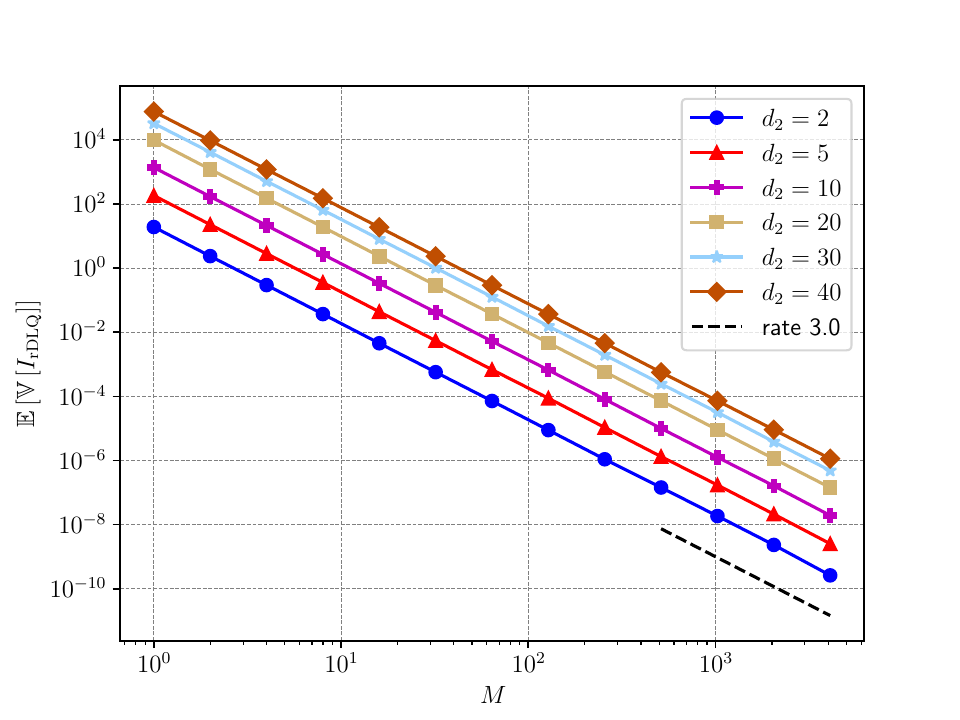}
	}\\
    \subfloat[Bias]{%
		\includegraphics[width=0.45\textwidth]{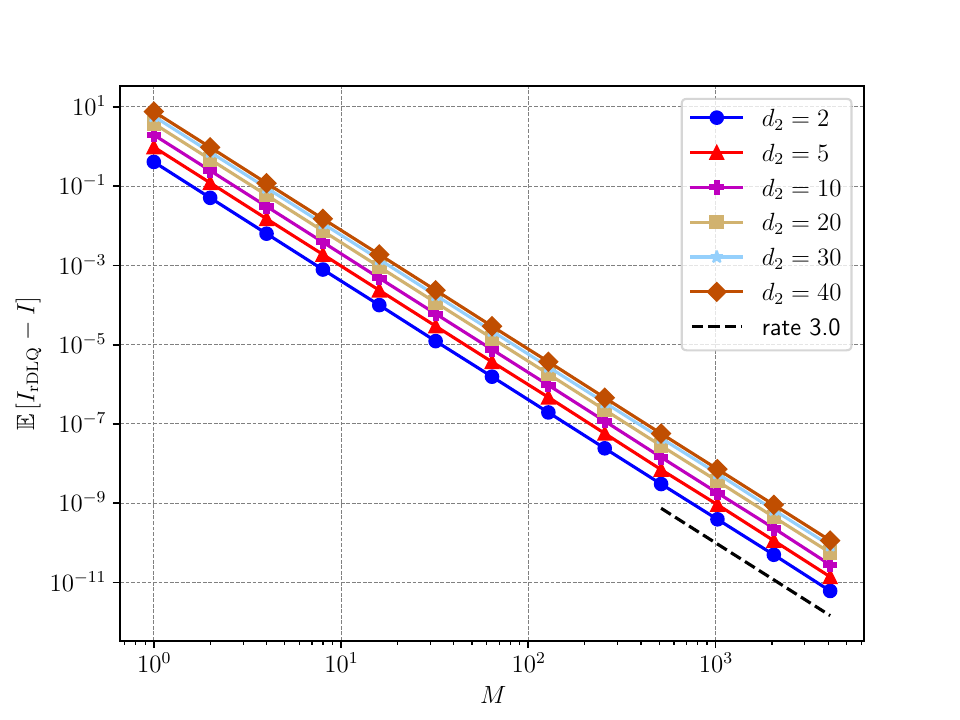}
	}
	\subfloat[Constants]{%
		\includegraphics[width=0.45\textwidth]{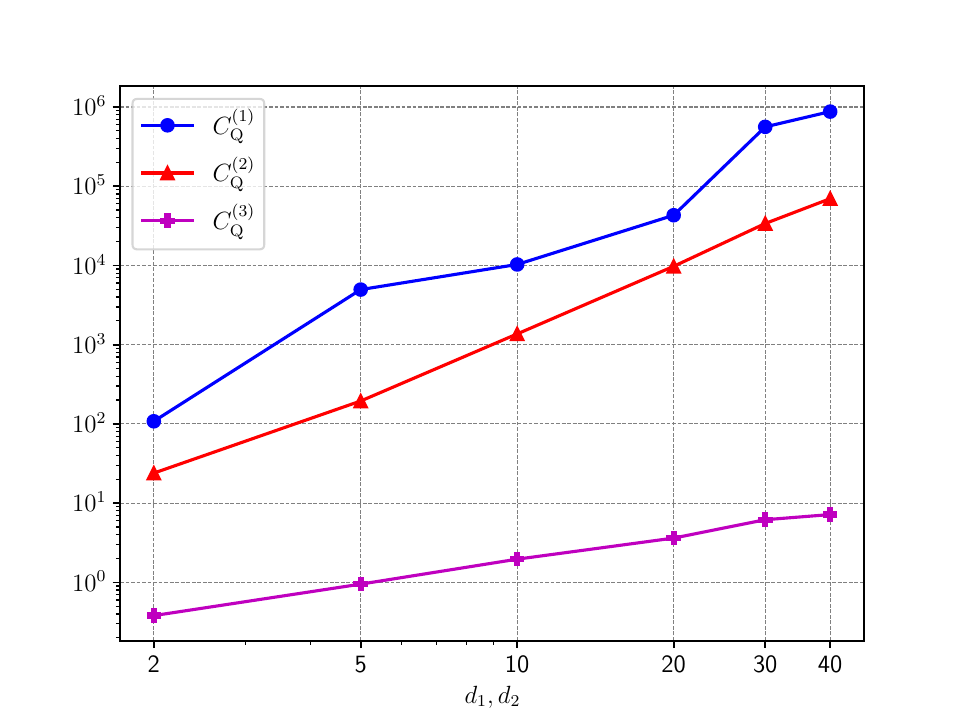}
	}
	\caption{Example 1 (polynomial example): results of the pilot runs to estimate the outer variance as a function of the outer number of points $N$, the inner variance and bias as functions of the inner number of points $M$, as well as the corresponding constants $C_{\rm{Q}}^{(1)}, C_{\rm{Q}}^{(2)}, C_{\rm{Q}}^{(3)}$ against the integrand dimensions. The same random seeds were used for all pilot runs. The bias differs from the inner variance by a factor independent of the inner samples.}
	\label{fig:ex0.dims}
\end{figure}

Before introducing the EIG formulation and demonstrating how to apply the rDLQMC estimator to the resulting nested integral, which is the primary application considered in this work, we present a simple polynomial example to verify our derived error bounds and efficiency analysis under the optimal regularity assumptions.
For this example, both the inner integrand $g$ and the outer integrand $f$ are polynomials over $[0,1]^{d_1}$ and $[0,1]^{d_2}$, respectively. Moreover, we analyze the effect of the integrand dimensions $d_1$ and $d_2$ on the multiplicative constants $C_{\rm{Q}}^{(1)}$, $C_{\rm{Q}}^{(2)}$, and $C_{\rm{Q}}^{(3)}$. We choose $d_1=d_2=d$ and consider the following outer integrand:
\begin{equation}\label{exa:1.f}
    f(z)\coloneqq az^2+bz, \quad \text{where }a,b,z\in\mathbb{R}.
\end{equation}
Moreover, we consider the following inner integrand:
\begin{equation}\label{exa:1.g}
    g(\bs{y},\bs{x})\coloneqq \bs{c}^{\trans}\bs{y} + \bs{d}^{\trans}\bs{x} + e\bs{y}^{\trans}\bs{x}, \quad \text{where } \bs{c},\bs{d}\in\mathbb{R}^{d},\, e\in\mathbb{R},\,\bs{y},\bs{x}\in[0,1]^{d}.
\end{equation}
The analytical solution to the nested integration problem $\int f\left(\int g(\bs{y},\bs{x})\di{}\bs{x}\right)\di{}\bs{y}$ for $f$ as in~\eqref{exa:1.f} and $g$ as in~\ref{exa:1.g} is given as follows:
\begin{align}\label{eq:ex0.I}
    I={}&\sum_{i=1}^d\left[a\left(\frac{c_i^2}{12}+e\frac{c_i}{3}+de\frac{d_i}{4}+(d-1)e\frac{c_i}{4}+\frac{c_i\sum_{j=1}^dc_j}{4}+e\frac{d_i\sum_{j=1}^dd_j}{4}+\frac{c_i\sum_{j=1}^dd_j}{2}\right)+b\left(\frac{c_i}{2}+\frac{d_i}{2}\right)\right]\nonumber\\
    {}&+d\frac{ae^2}{12}+d(d-1)\frac{ae^2}{16}+d\frac{be}{4}.
\end{align}
For this example, we choose all coefficients to be one, that is, $a=b=e=1$, and $\bs{c}=\bs{d}=(1,\ldots,1)^{\trans}$, resulting in the following outer and inner integrands:
\begin{equation}\label{def:pol.f}
    f(z)\coloneqq z^2+z, \quad \text{where } z\in\mathbb{R}
\end{equation}
and
\begin{equation}\label{def:pol.g}
    g(\bs{y},\bs{x})\coloneqq \sum_{i=1}^{d}\left[y_i + x_i +  y_ix_i\right], \quad \text{where } \bs{y},\bs{x}\in[0,1]^{d}.
\end{equation}
It follows immediately that both $f$ and $g$ as defined in~\eqref{def:pol.f} and~\eqref{def:pol.g} have bounded generalized Hardy--Krause variation of order one. Figure~\ref{fig:ex0.dims} shows the results of the pilot runs described in Algorithm~\ref{alg.pil}. We ran the pilot twice, once with $S=2^7$, $N=2^{15}$, $R=1$, and $M=4$ to estimate the outer variance and once with $S=1$, $N=4$, $R=2^{15}$, and $M=2^{12}$ to estimate the inner variance and bias. This procedure was repeated for dimensions $d_1,d_2=d\in\{2,5,10,20,30,40\}$. We observed convergence rates $\beta$ and $\delta$ close to two, implying convergence of the outer variance (Panel (A)) and inner variance (Panel (B)) as well as the bias (Panel (C)) at a rate close to three in all dimensions considered. The associated multiplicative constants display exponential growth in the dimension (Panel (C)).

\begin{figure}[ht]
	\subfloat[Computational work]{%
		\includegraphics[width=0.45\textwidth]{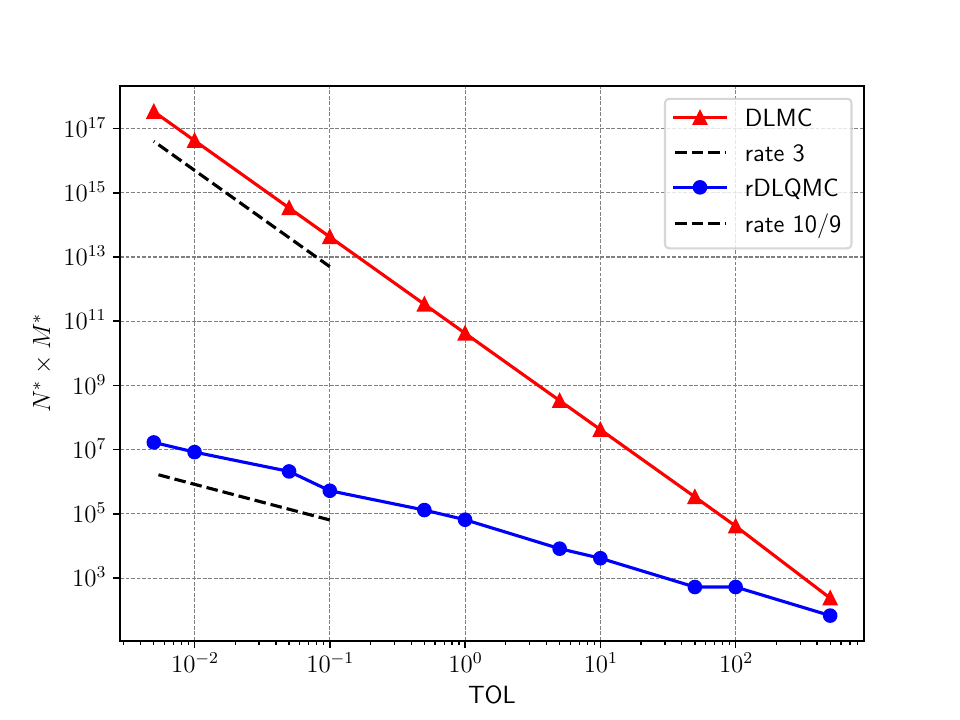}
	}
    	\hfill
	\subfloat[Error vs.~tolerance]{%
	\begin{tikzpicture}
	\node[inner sep=0pt] (comp) at (0,0)
		{\includegraphics[width=0.45\textwidth]{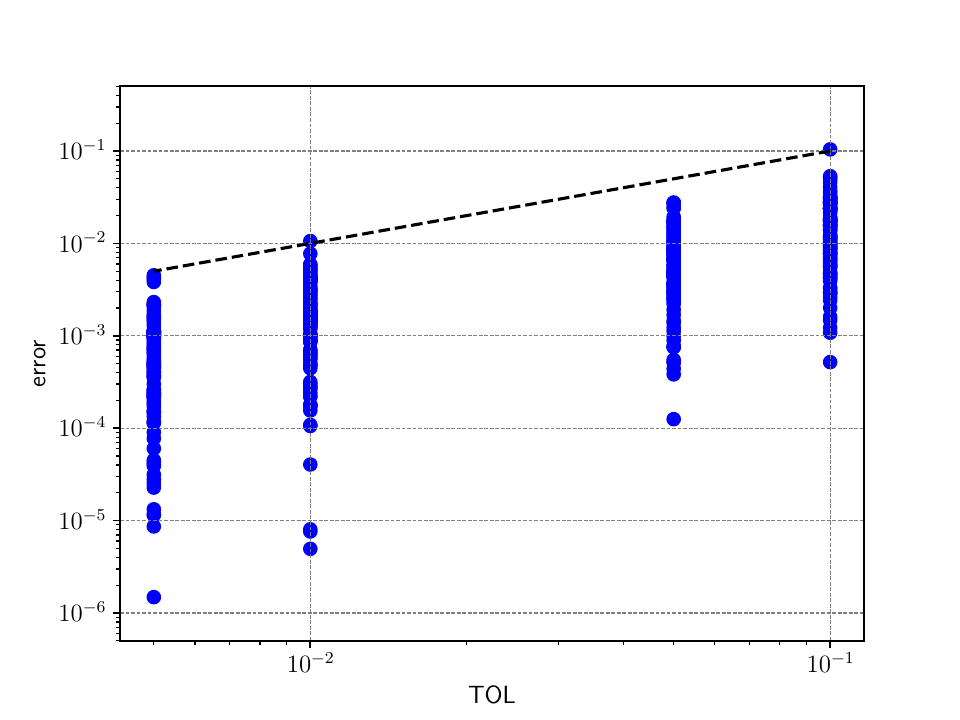}};
		\node at (-2.54,0.9) {\small \red{0}};
    \node at (-1.3,1.15) {\small \red{1}};
      \node at (1.5,1.65) {\small \red{0}};
       \node at (2.7,1.9) {\small \red{1}};
		\end{tikzpicture}
	}
	\caption{Example 1 (polynomial example): optimal computational work $N^{\ast}\times M^{\ast}$ and error vs.~tolerance consistency plot for confidence parameter $C_{\alpha}=1.96$ ($\alpha=0.05$) in dimension $d_1=d_2=30$ for the DLMC and the rDLQMC estimators. The samples $N^{\ast}$ and $M^{\ast}$ are rounded up to a power of two for the rDLQMC estimator to ensure optimal behavior.}
	\label{fig:ex0.EvT}
\end{figure}

Based on the results of the pilot runs, we obtained the optimal numbers of samples $N^{\ast}$ and $M^{\ast}$ as functions of the error tolerance $TOL$ in dimension $d_1=d_2=30$. The total computational work for the rDLQMC estimator without discretization of the inner integrand is proportional to $N^{\ast}\times M^{\ast}$. Figure~\ref{fig:ex0.EvT} (A) presents the rate of increase of the optimal work as a function of $TOL$ at rate $10/9$ as implied by Corollary~\ref{cor:optimal.work.s}. This rate is possible because of the combination of Sobol' points with Owen's scrambling. For comparison, the DLMC estimator requires computational work increasing at rate three as a function of the error tolerance. We used $N=2^{22}$ and $M=4$ MC samples to estimate the outer variance and $N=4$ and $M=2^{17}$ MC samples to estimate the inner variance and bias for the DLMC estimator. To verify that the rDLQMC estimator with $S=R=1$ randomizations and $N^{\ast}$ and $M^{\ast}$ derived from the pilot runs actually achieves the required accuracy, we compared 100 runs for various tolerances with the analytical solution in~\eqref{eq:ex0.I} as a reference. The results in Panel (B) indicate that for confidence parameter $C_{\alpha}=1.96$, less than five out of 100 runs exceeded the error tolerance; thus, the rates and constants obtained from the pilots appear to be conservative estimates.



\FloatBarrier

\section{Expected information gain estimation}\label{sec:EIG.estimation}
In Bayesian OED \cite{Cha95}, the aim is to maximize the EIG of an experiment, given by the expected Kullback--Leibler divergence \cite{Sha48, Kul59, Kul51, Lin56} of the posterior distribution of the parameters of interest $\bs{\theta}$ with respect to their prior distribution. We assumed the following data model:
\begin{equation}\label{eq:data.model}
    \bs{Y}_i(\bs{\xi}) = \bs{G}(\bs{\theta}_t,\bs{\xi}) + \bs{\varepsilon}_i, \quad 1\leq i\leq N_e,
\end{equation}
where $\bs{Y}=(\bs{Y}_1,\ldots,\bs{Y}_{N_e})\in\bb{R}^{d_y\times N_e}$ represents noisy data generated from the deterministic model $\bs{G}:\bb{R}^{d_{\theta}}\times\bb{R}^{d_{\xi}}\to\bb{R}^{d_y}$ (evaluated at the true parameter vector $\bs{\theta}_t\in\bb{R}^{d_{\theta}}$). The design of the experiment is denoted as $\bs{\xi}\in\bb{R}^{d_{\xi}}$. Random observation noise is denoted by $\bs{\varepsilon}_i\in\bb{R}^{d_y}$ for $N_e$ available independent observations under a consistent experimental setup, where $d_y, d_{\theta}, d_{\xi}$ are positive integers. The noise $\bs{\varepsilon}_i$ is assumed to follow a centered normal distribution with a known covariance matrix $\bs{\Sigma}_{\bs{\varepsilon}}\in\bb{R}^{d_y\times d_y}$ independent of $\bs{\theta}$ and $\bs{\xi}$. The knowledge of the parameters of interest $\bs{\theta}$ before experimenting is encompassed in the prior probability density function (PDF) $\pi(\bs{\theta})$. After the experiment, the knowledge is described by the posterior PDF of $\bs{\theta}$, given by Bayes' theorem as follows:
\begin{equation}\label{eq:Bayes}
  \pi(\bs{\theta}|\bs{Y},\bs{\xi})=\frac{p(\bs{Y}|\bs{\theta},\bs{\xi})\pi(\bs{\theta})}{p(\bs{Y}|\bs{\xi})},
\end{equation}
where
\begin{equation}\label{eq:likelihood}
 p(\bs{Y}|\bs{\theta},\bs{\xi})\coloneqq \det(2\pi\bs{\Sigma_{\bs{\varepsilon}}})^{-\frac{N_e}{2}}\exp\left(-\frac{1}{2}\sum_{i=1}^{N_e}\bs{r}(\bs{Y}_i,\bs{\theta},\bs{\xi})\cdot\bs{\Sigma}_{\bs{\varepsilon}}^{-1}\bs{r}(\bs{Y}_i,\bs{\theta},\bs{\xi})\right)
\end{equation}
represents the likelihood function. The data residual is given by \eqref{eq:data.model} as follows:
\begin{equation}\label{eq:residual}
\bs{r}(\bs{Y}_i, \bs{\theta},\bs{\xi})\coloneqq \bs{Y}_i(\bs{\xi})-\bs{G}(\bs{\theta},\bs{\xi}), \quad 1\leq i\leq N_e.
\end{equation}
We omitted the design parameter $\bs{\xi}$ for simplicity and distinguished between $\pi$ for the PDFs of the parameters of interest and $p$ for the PDFs of the data.
The amount of information regarding $\bs{\theta}$ gained from the experiment is given by
\begin{equation}\label{eq:EIG.posterior}
    EIG\coloneqq \int_{\cl{Y}}\int_{\bs{\Theta}}\log\left(\frac{\pi(\bs{\theta}|\bs{Y})}{\pi(\bs{\theta})}\right)\pi(\bs{\theta}|\bs{Y})\di{}\bs{\theta}p(\bs{Y})\di{}\bs{Y},
\end{equation}
which is rewritten in terms of the likelihood function and the prior using Bayes' theorem and marginalization:
\begin{align}\label{eq:EIG}
    EIG{}&= \int_{\bs{\Theta}}\int_{\cl{Y}}\log\left(\frac{p(\bs{Y}|\bs{\theta})}{\int_{\bs{\Theta}}p(\bs{Y}|\bs{\vartheta})\pi(\bs{\vartheta})\di{}\bs{\vartheta}}\right)p(\bs{Y}|\bs{\theta})\di{}\bs{Y}\pi(\bs{\theta})\di{}\bs{\theta},\nonumber\\
    {}&=\int_{\bs{\Theta}}\int_{\cl{Y}}\log\left(p(\bs{Y}|\bs{\theta})\right)p(\bs{Y}|\bs{\theta})\di{}\bs{Y}\pi(\bs{\theta})\di{}\bs{\theta}-\int_{\bs{\Theta}}\int_{\cl{Y}}\log\left(\int_{\bs{\Theta}}p(\bs{Y}|\bs{\vartheta})\pi(\bs{\vartheta})\di{}\bs{\vartheta}\right)p(\bs{Y}|\bs{\theta})\di{}\bs{Y}\pi(\bs{\theta})\di{}\bs{\theta},
\end{align}
where $\bs{\vartheta}$ indicates a dummy variable for integration. The second term in \eqref{eq:EIG} forms a nested integration problem to which the rDLQMC estimator is applied. Specifically, in this setting, the nonlinear outer function $f$ in \eqref{eq:double.loop} is the logarithm, and the inner function $g$ is the likelihood function.

\begin{rmk}[Closed-form expression for the first term in the EIG]\label{rmk:num}
From \eqref{eq:residual}, the likelihood in the logarithm in the first term of \eqref{eq:EIG} depends only on $\bs{\varepsilon}_i$, where $1\leq i\leq N_e$, and can be computed in closed form for $\bs{\varepsilon}_i\sim\cl{N}\left(\bs{0},\bs{\Sigma}_{\bs{\varepsilon}}\right)$, where $\bs{\Sigma}_{\bs{\varepsilon}}$ is a positive definite diagonal matrix in $\bb{R}^{d_y\times d_y}$ with entry $\sigma_{\varepsilon\{j,j\}}^2$, where $1\leq j\leq d_y$. The solution to this integration problem in \eqref{eq:EIG} is given by
\begin{equation}\label{eq:numerator}
   -\frac{N_e}{2}\sum_{j=1}^{d_y} \left(\log\left(2\pi\sigma_{\varepsilon\{j,j\}}^2\right)+1\right).
\end{equation}
A similar result is presented in \cite{Tsi17}. Appendix~\ref{ap:Rmk} provides a derivation for completeness.
\end{rmk}

The first term in~\eqref{eq:numerator} relates to the scaling factor $\det(2\pi\bs{\Sigma_{\bs{\varepsilon}}})^{-\frac{N_e}{2}}$ for diagonal $\bs{\Sigma_{\bs{\varepsilon}}}$. As this factor appears in both terms in~\eqref{eq:EIG}, it cancels out and need not be computed. To apply the rDLQMC estimator to the EIG, we introduce $\bs{y}_1\coloneqq \bs{\theta}$, $\bs{\Sigma}_{\bs{\varepsilon}}^{\frac{1}{2}}\Phi^{-1}(\bs{y}_{2,i})\coloneqq \bs{\varepsilon}_i$ for $1\leq i\leq N_e$, where $\bs{y}_2\coloneqq(\bs{y}_{2,1},\ldots,\bs{y}_{2,N_e})$, and $\bs{x}\coloneqq \bs{\vartheta}$. The rDLQMC estimator for the EIG with uniform prior is thus given as follows:
\begin{align}\label{eq:dlqmc.S.R.EIG}
I_{\rm{rDLQ}}^{(S,R)} ={}&-\frac{N_ed_{y}}{2} 
-\frac{1}{S}\sum_{s=1}^S\frac{1}{N}\sum_{n=1}^N\log\left(\frac{1}{R}\sum_{r=1}^R\frac{1}{M}\sum_{m=1}^M e^{-\frac{1}{2}\sum_{i=1}^{N_e}\left\lVert \bs{G}(\bs{y}_1^{(s,n)})+\bs{\Sigma}_{\bs{\varepsilon}}^{\frac{1}{2}}\Phi^{-1}(\bs{y}_{2,i}^{(s,n)})-\bs{G}(\bs{x}^{(s,n,r,m)})\right\rVert^2_{\bs{\Sigma}_{\bs{\varepsilon}}^{-1}}}\right),
\end{align}
where $\bs{y}\coloneqq(\bs{y}_1,\bs{y}_{2})\in[0,1]^{d_{\theta}+N_e\times d_{y}}$ and $\bs{x}\in[0,1]^{d_{\theta}}$ (see Remark~\ref{ex:EIG}). We thus have that $d_1=d_{\theta}+N_e\times d_{y}$ and $d_2=d_{\theta}$. For a Gaussian prior with given mean and covariance, an appropriate inverse CDF and a shift must be applied to $\bs{y}_1$ and $\bs{x}$ as well. Appendices~\ref{app:cond} to~\ref{app:GHK} demonstrate how to apply the error bounds derived in the previous section to EIG estimation. 
A truncation of the Gaussian observation noise is introduced to bound higher-order terms in the Taylor expansions in the proofs of Propositions~\ref{prop:B.DLQ} and~\ref{prop:V.DLQ}. Numerical experiments show that in practice, the convergence rates stated in those propositions are still attained, even without truncating the observation noise. The inner integrand moreover has bounded generalized Hardy--Krause variation of order one. Although the outer integrand is bounded once the truncation is applied, the resulting error bound for the outer variance grows at a rate of essentially $\cl{O}(TOL^{-d_y})$ as the region of truncation approaches the boundary of the integration domain; thus negating the advantage incurred from the improved convergence described in Corollary~\ref{cor:optimal.work.s}. These results are detailed in Appendix~\ref{app:GHK}. Nonetheless, numerical experiments demonstrate that implementing the truncation improves multiplicative factors in the observed outer variance for the considered tolerances. In Appendix~\ref{app:cond}, we derived the condition that the data model $\bs{G}$ applied to the inverse CDF of the parameters of interest must be uniformly Lipschitz continuous and have bounded mixed derivatives for Proposition~\ref{prop:owen} to guarantee good rQMC convergence rates for the EIG. In Appendix~\ref{app:Condition.42}, we verified Condition~\eqref{eq:reg:g} for the observation noise distributed according to a truncated normal distribution and demonstrated weak dependence of the constant
\begin{equation}\label{eq:k}
    k=\cl{O}\left(\log(TOL^{-1})^{\frac{
    d_{\theta}}{2}}\right)
\end{equation}
in Condition~\eqref{eq:reg:g} on $TOL$. The additional error introduced by this truncation was also bounded in Appendix~\ref{app:Truncation.EIG} as $o(TOL)$ for observation noise distributed according to a normal distribution. In Corollary~\ref{cor:total.error}, the total error of the rDLQMC estimator applied to estimate the EIG for observation noise from a standard normal distribution is presented, specifying the dependence of the bounds on the bias and statistical error on the error tolerance $TOL$ introduced via $k$ in \eqref{eq:k}. The verification of the Lipschitz condition for specific data models is left for future research. A similar truncation may be used also for parameters of interest with a normal prior.

Furthermore, an importance sampling distribution for $\bs{\vartheta}$ based on the Laplace approximation \cite{Lon13, Lon21, Sch20, Bec18, Bec20, Wac17, Sti86, Tie86, Tie89, Kas90, Lon15, Hel22, Ouy23, Bor11} can reduce the variance of the (rQ)MC methods even further (for a thorough discussion of the Bayesian OED formulation, see
\cite{Lon13, Bec18, Car20}). Without importance sampling, the inner integral approximation can produce numerical underflow \cite{Bec18} unless a vast number of inner samples is used, rendering these methods impractical in those cases.
Alternatively, a single-loop (rQ)-MC method can be used to estimate the EIG by approximating the posterior in \eqref{eq:EIG.posterior} by the Laplace approximation. In this case, no nested integration is necessary. Such estimators were introduced in the Ph.D.~                thesis~\cite{Bar25}. We did not encounter numerical underflow in our numerical experiments and the analysis of such estimators is left for future research.
We now demonstrate the effectiveness and applicability of the derived double-loop methods for two EIG examples numerically. The first example is linear in the parameters of interest. Moreover, a normal prior is assumed for the parameters of interest; thus, there exists a closed-form solution for the EIG. We use this setting to verify that the error convergence rates derived in previous sections can be obtained in practice with and without truncation of the observation noise. We also observe the effect of the size of the noise covariance on the approximation error.

The second example involves the solution to a PDE and showcases the effectiveness of the proposed method for a more applied example. The outer integrand in this latter example is in dimension 15, which is moderately high-dimensional in this context.
\subsection{Example 2: Linear EIG example with exact sampling}\label{sec:L.EIG}
We start with a simple experimental model adapted from \cite{Fen19}, given as follows:
\begin{equation}
    \bs{G}(\bs{\theta},\xi)\coloneqq \bs{A}(\xi)\bs{\theta},
\end{equation}
where $\bs{\theta}\sim\cl{N}(\bs{0},\bs{I}_{2\times 2})$ are the parameters of interest and
\begin{equation}
    \bs{A}(\xi)\coloneqq \begin{pmatrix}
        \xi & 0\\
        0 & 1-\xi
    \end{pmatrix},
\end{equation}
where $\xi\in[0,1]$ indicates the design. For this experiment, we consider additive Gaussian noise, that is, 
\begin{equation}
    \bs{Y}=\bs{G}(\bs{\theta},\xi)+\bs{\varepsilon},
\end{equation}
where $\bs{\varepsilon}\sim\cl{N}(\bs{0},\bs{\Sigma}_{\bs{\varepsilon}})$, and $\bs{\Sigma}_{\bs{\varepsilon}}\coloneqq\sigma_{\varepsilon}^{2}\bs{I}_{2\times 2}$. We only consider one experiment, that is, $N_e=1$, and thus the outer integration is $N_e\times d_y+d_{\theta}=4$ dimensional, and the inner integration is $d_{\theta}=2$ dimensional. We focus on the case $\xi=0.5$ and obtain the following closed-form expression for the EIG:
\begin{equation}
    EIG = \log\left(1+\frac{0.25}{\sigma_{\varepsilon}^2}\right).
\end{equation}
For now, we consider the noise covariance fixed at $\sigma_{\varepsilon}\coloneqq10^{-1}$. This example provides several interesting features. There is a closed-form solution for the EIG, allowing us to test the error convergence rates against the true solution. Moreover, the Gaussian distribution of both the noise and parameters of interest make it challenging for rQMC estimators. We therefore use this example to showcase the effect of the truncation of the Gaussian noise and parameters of interest on the rQMC convergence rates in the inner and outer estimators. Finally, we can easily adjust the experiment noise to study its effect on the rQMC error; in particular on the multiplicative constants.

\begin{figure}[ht]
	\subfloat[Outer variance]{%
		\includegraphics[width=0.45\textwidth]{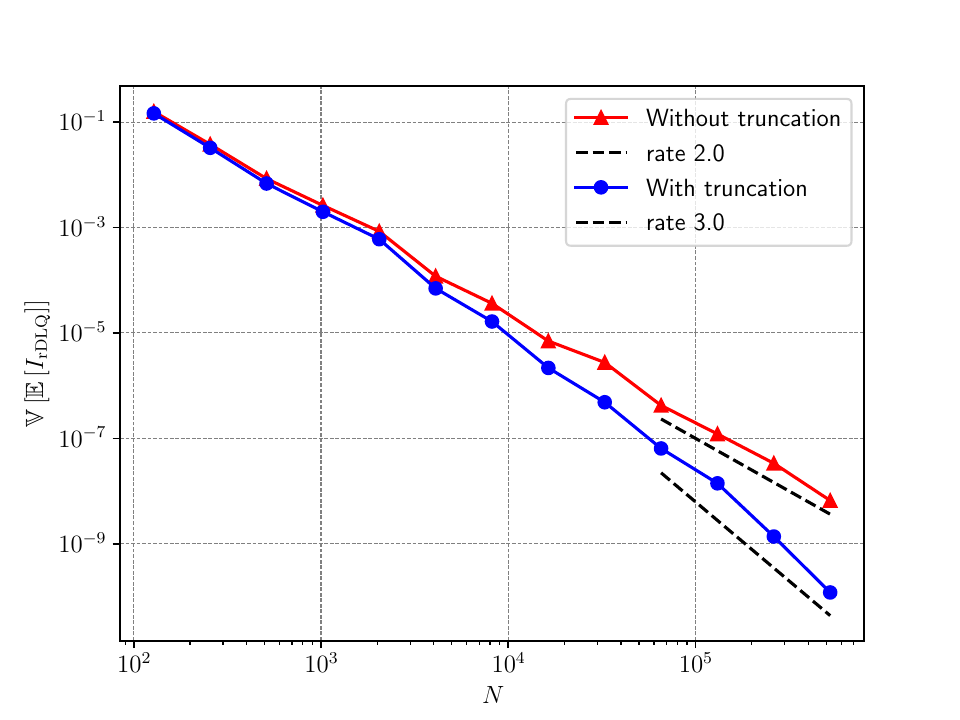}
	}
    \subfloat[Inner variance]{%
		\includegraphics[width=0.45\textwidth]{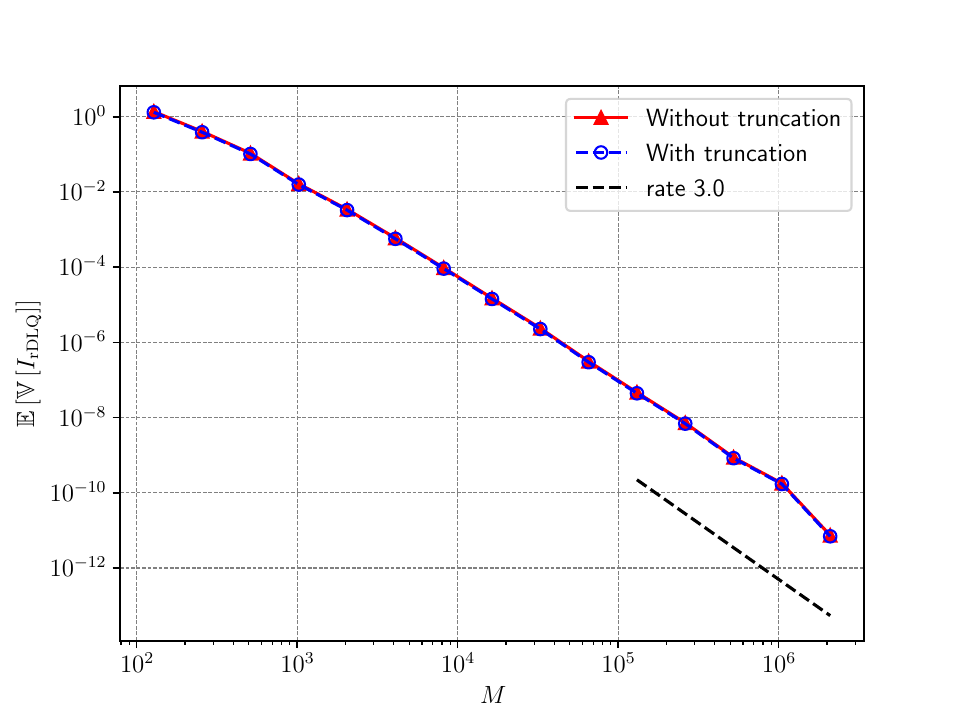}
	}\\
    \subfloat[Bias]{%
		\includegraphics[width=0.45\textwidth]{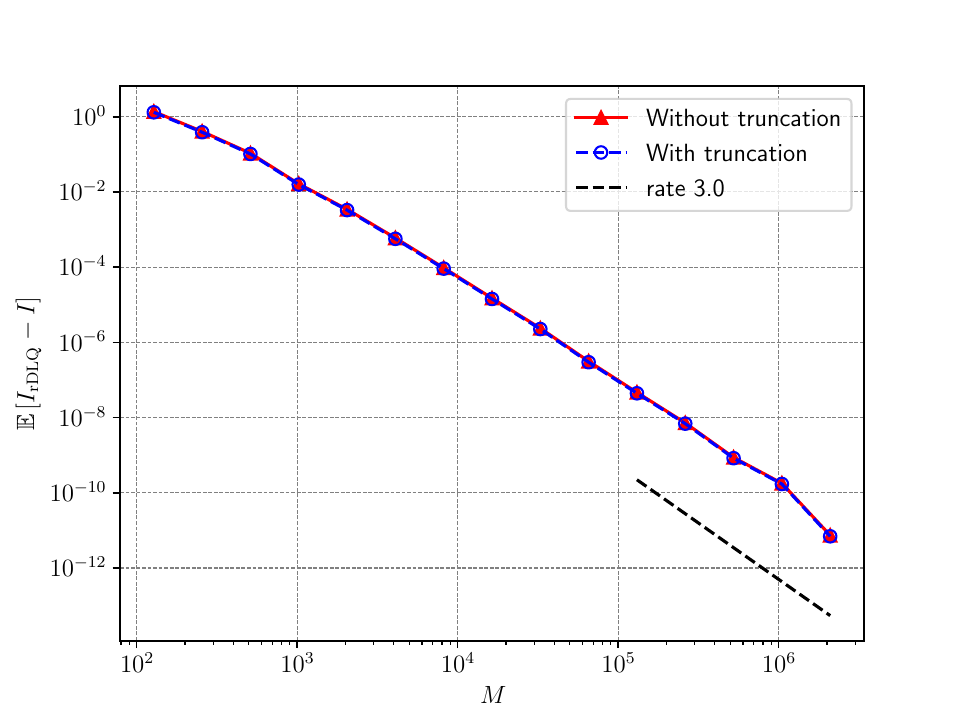}
	}
    \caption{Example 2 (linear EIG example): results of the pilot runs to estimate the outer variance, inner variance, and bias with and without truncation of the observation noise and parameters of interest. The bias and inner variance were unaffected by the truncation. Moreover, bias and inner variance only differ from each other by a multiplicative constant. The same random seeds were used for all pilot runs.}
	\label{fig:ex0.5.pilot}
\end{figure}

\begin{figure}[ht]
    \subfloat[Outer variance as a function of $TOL$]{%
		\includegraphics[width=0.45\textwidth]{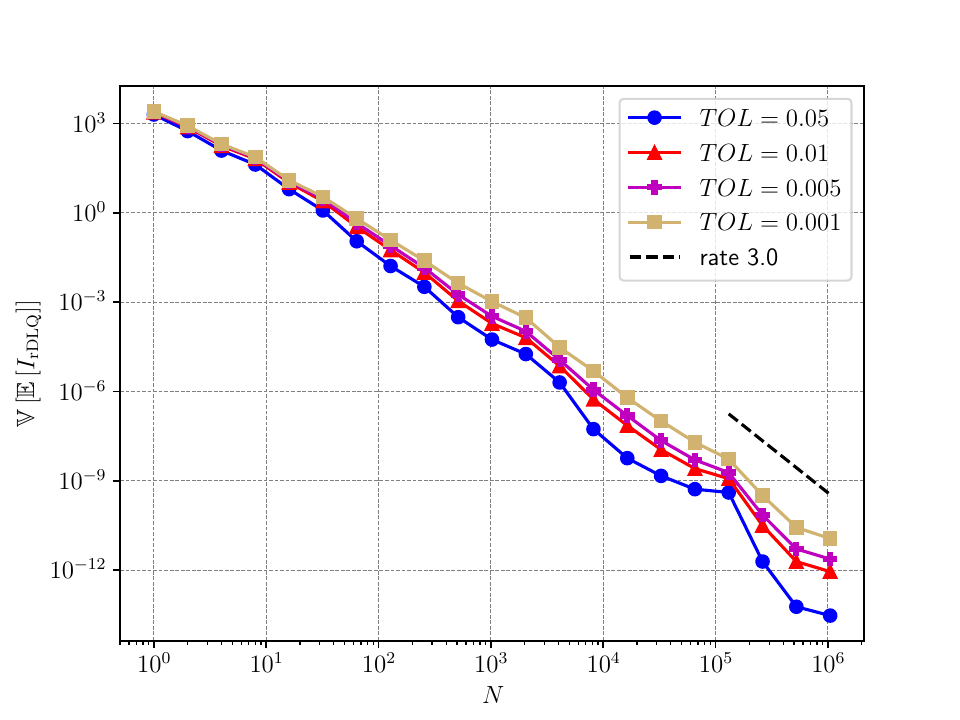}
	}
    \subfloat[Constant $C_{\rm{Q}}^{(1)}$ as a function of $TOL$]{%
		\includegraphics[width=0.45\textwidth]{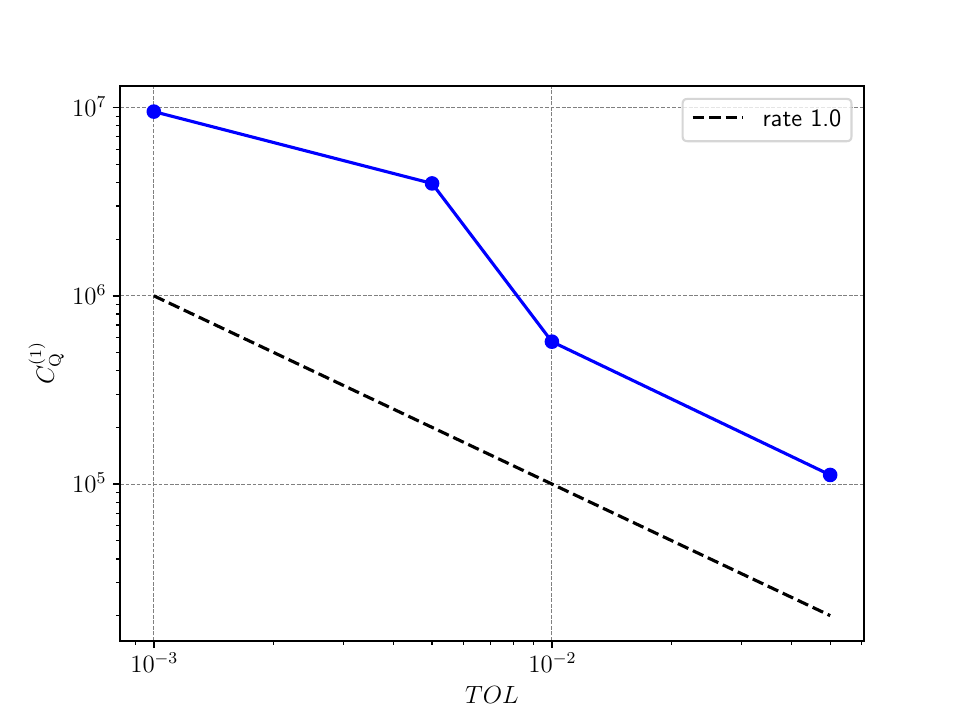}
	}
    \caption{Example 2 (linear EIG example): results of the pilot runs to estimate the outer variance with a truncation of the observation noise and parameters of interest that depends on the error tolerance $TOL$. As $TOL$ decreases, the multiplicative constant $C_{\rm{Q}}^{(1)}$ increases approximately at a rate of one.}
	\label{fig:ex0.5.pilot.truncation.c}
\end{figure}

Figure~\ref{fig:ex0.5.pilot} presents the results of the pilot runs. We ran the pilot with $S=2^6$, $N=2^{19}$, $R=1$, and $M=1$ to estimate the outer variance (Panel (A)) and with $S=1$, $N=1$, $R=2^9$, and $M=2^{21}$ to estimate the inner variance (Panel (B)) and bias (Panel (C)). Then we ran the pilots with the same settings with the truncation of the observation noise and prior. We observed a convergence rate of $1+\beta\approx3.07$ for the outer variance when using the truncation, whereas only a rate of rate $1+\beta\approx2.03$ was observed without truncation. Rate $1+\delta\approx3.06$ was observed for the inner variance and bias with and without truncation. For the truncation, we used 
\begin{equation}
    c(TOL)=(2(1+p))^{\frac{1}{2}}\log(TOL^{-1})^{\frac{1}{2}},
\end{equation}
for a fixed $TOL=10^{-3}$ and $p=0.5$. The truncation depends on the error tolerance $TOL$ to ensure consistent estimates. As $TOL$ decreases, the region of truncation extends towards the boundary as specified by $c(TOL)$. The multiplicative constant $C_{\rm{Q}}^{(1)}$ estimated in the pilot runs includes the squared generalized Hardy--Krause variation of order one of the truncated nested integrand (see~\eqref{eq:variance.ghk}). Theoretical results in Lemma~\ref{lem:GHK} and the subsequent discussion indicate that this factor exhibits unfavorable dependence on the dimension $d_1$, similarly to the result in~\cite[Theorem 6.2]{Kaa24}. For $d_1=4$, however, the analyses corresponding to the Hardy--Krause variation and Owen's boundary growth condition provide tighter error bounds for this example. In the numerical experiments, we observed an increase of approximately $\cl{O}(TOL^{-1})$ for $p=\epsilon=0$ in Figure~\ref{fig:ex0.5.pilot.truncation.c}. This suggests an overall dependence on $TOL$ comparable to that predicted by the analyses based on the Hardy--Krause variation and Owen's boundary growth condition. For practical purposes, $p>0$ must be chosen to ensure the required bound on the truncation error ($o(TOL)$, see Corollary~\ref{cor:total.error}). However, as $TOL$ approaches zero, smaller and smaller $p$ may be selected. It is therefore challenging to obtain reliable estimates of the required constants and rates from pilot runs. An adaptive procedure similar to the one discussed in~\cite{Col14} may be more suitable. Both methods (with and without truncation) eventually yield comparable overall increase in computational cost as $TOL$ approaches zero; however, lower outer variance, and thus total cost, can be observed for the truncated outer integrand for large or intermediate values of $TOL$.

\begin{figure}[ht]
    \subfloat[Optimal computational work $N^{\ast}\times M^{\ast}$ as a function of the error tolerance $TOL$. The theoretically optimal rate is displayed as a reference]{%
		\includegraphics[width=0.45\textwidth]{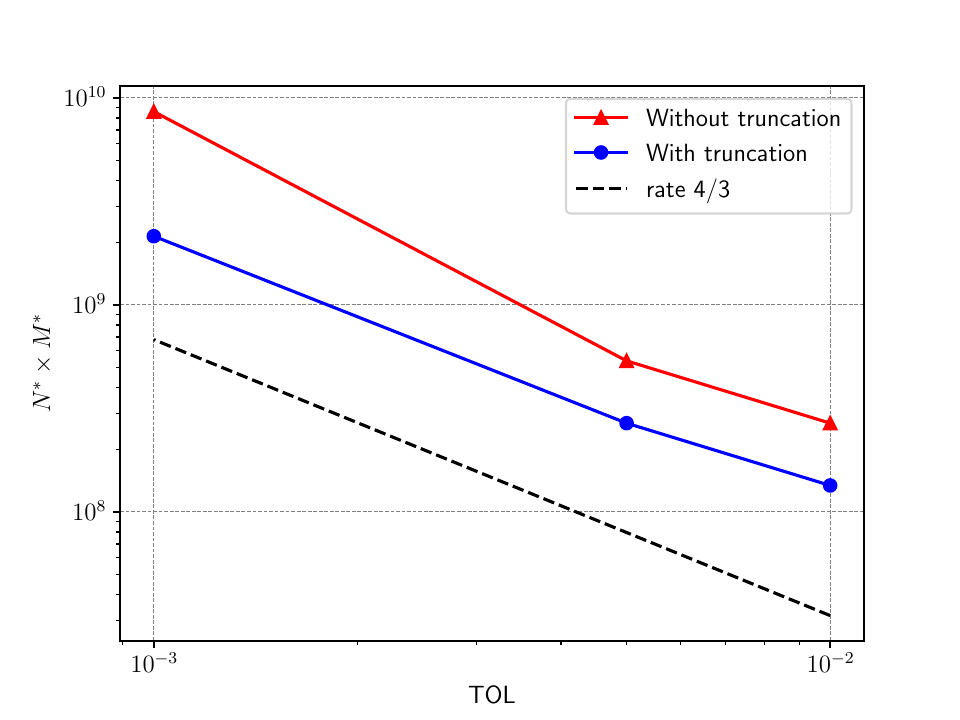}
	}\\
	\subfloat[Error vs.~tolerance without truncation]{%
	\begin{tikzpicture}
	\node[inner sep=0pt] (comp) at (0,0)
		{\includegraphics[width=0.45\textwidth]{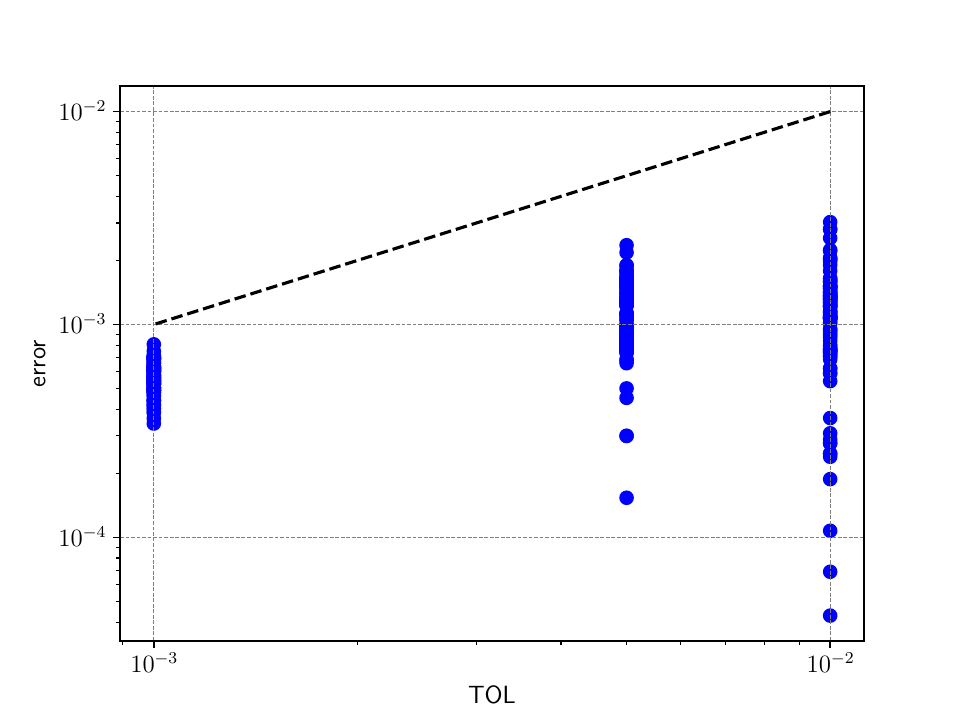}};
		\end{tikzpicture}
	}
	\hfill
	\subfloat[Error vs.~tolerance with truncation]{%
	\begin{tikzpicture}
	\node[inner sep=0pt] (comp) at (0,0)
		{\includegraphics[width=0.45\textwidth]{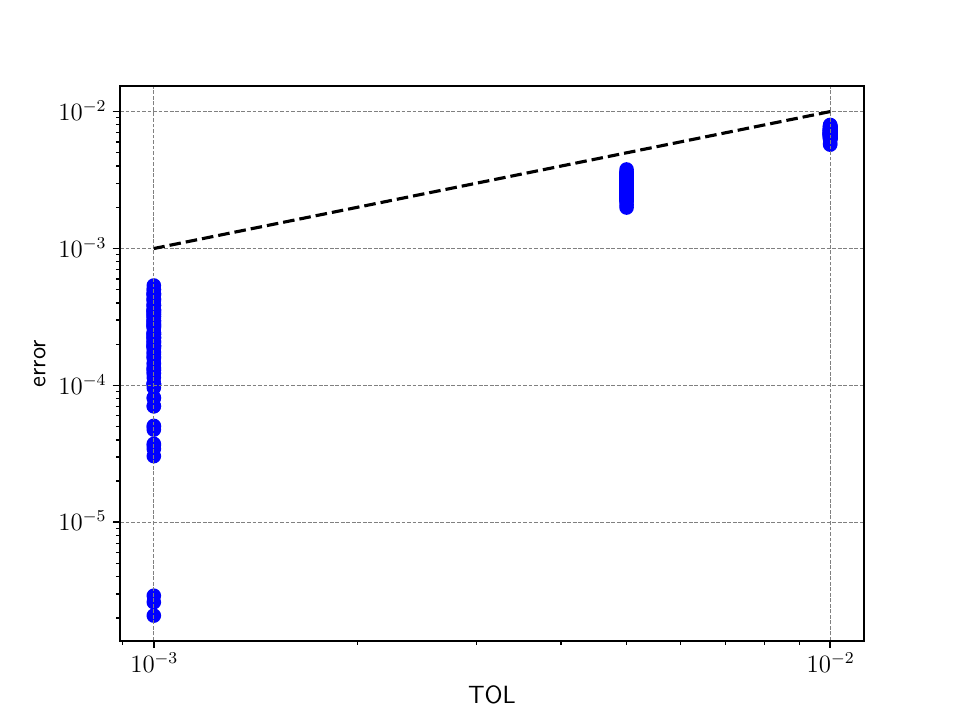}};
		\end{tikzpicture}
	}
	\caption{Example 2 (linear EIG example): optimal computational work $N^{\ast}\times M^{\ast}$ and error vs.~tolerance consistency plot for confidence parameter $C_\alpha=1.96$ ($\alpha=0.05$) with and without truncation of the observation noise and parameters of interest. The samples $N^{\ast}$ and $M^{\ast}$ are rounded up to a power of two to ensure optimal behavior of the rQMC estimators.}
	\label{fig:ex0.5.evt}
\end{figure}

From the pilot estimates, we obtained the optimal numbers of samples $N^{\ast}$ and $M^{\ast}$ as a function of $TOL$. Figure~\ref{fig:ex0.5.evt} (A) demonstrates that applying the truncation improves the computational cost compared to not using the truncation for the considered tolerances. However, as $TOL$ approaches zero, the cost for both methods increases at a rate of $4/3$. The accuracy of the rDLQMC estimator with and without truncation is demonstrated in Panel~(B) and Panel~(C), respectively for 100 runs and varying tolerances.

\begin{figure}[ht]
	\subfloat[Outer variance]{%
		\includegraphics[width=0.45\textwidth]{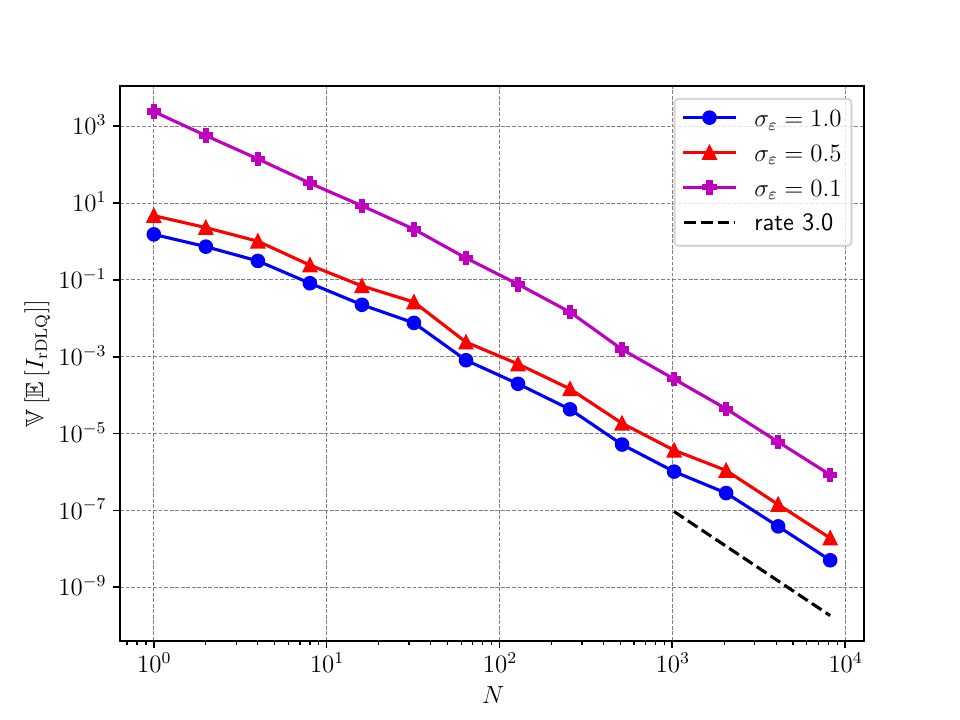}
	}
    \subfloat[Inner variance]{%
		\includegraphics[width=0.45\textwidth]{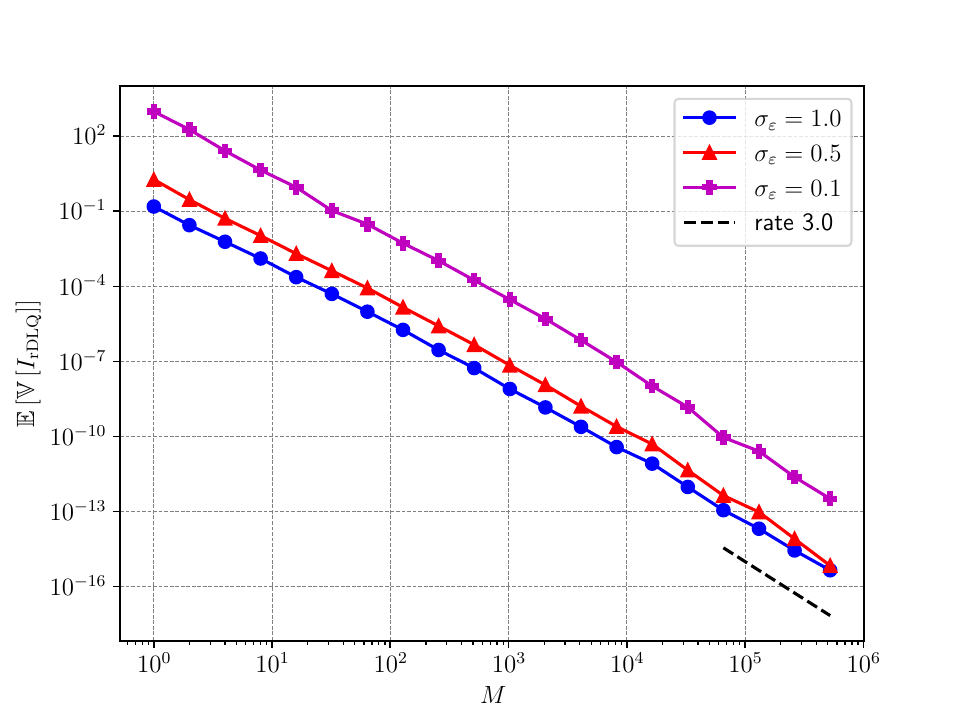}
	}\\
    \subfloat[Bias]{%
		\includegraphics[width=0.45\textwidth]{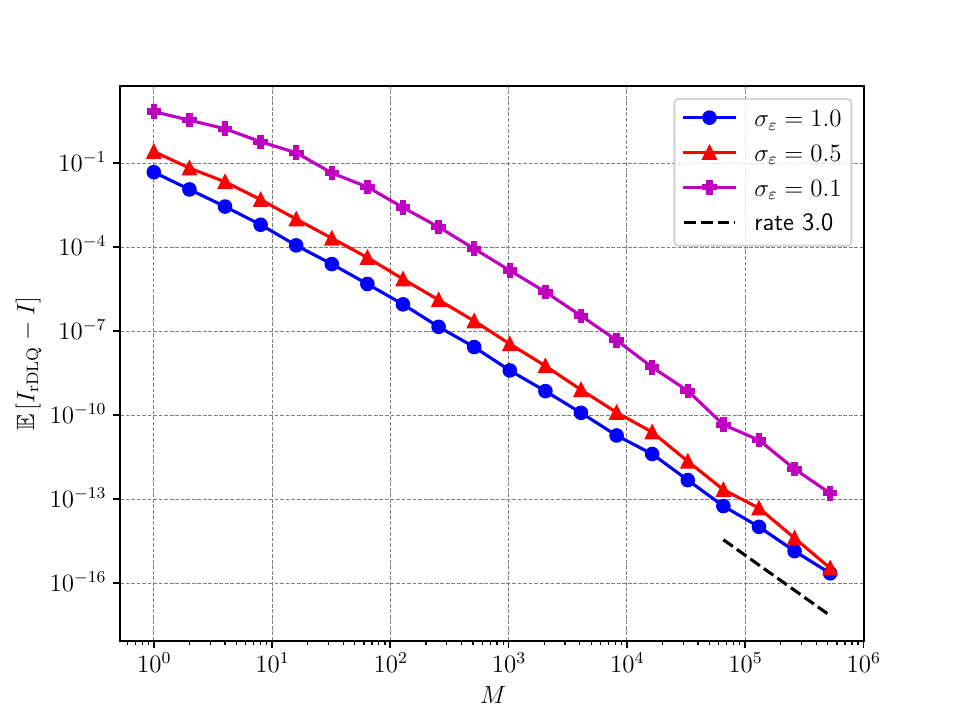}
	}
	\subfloat[Constants]{%
		\includegraphics[width=0.45\textwidth]{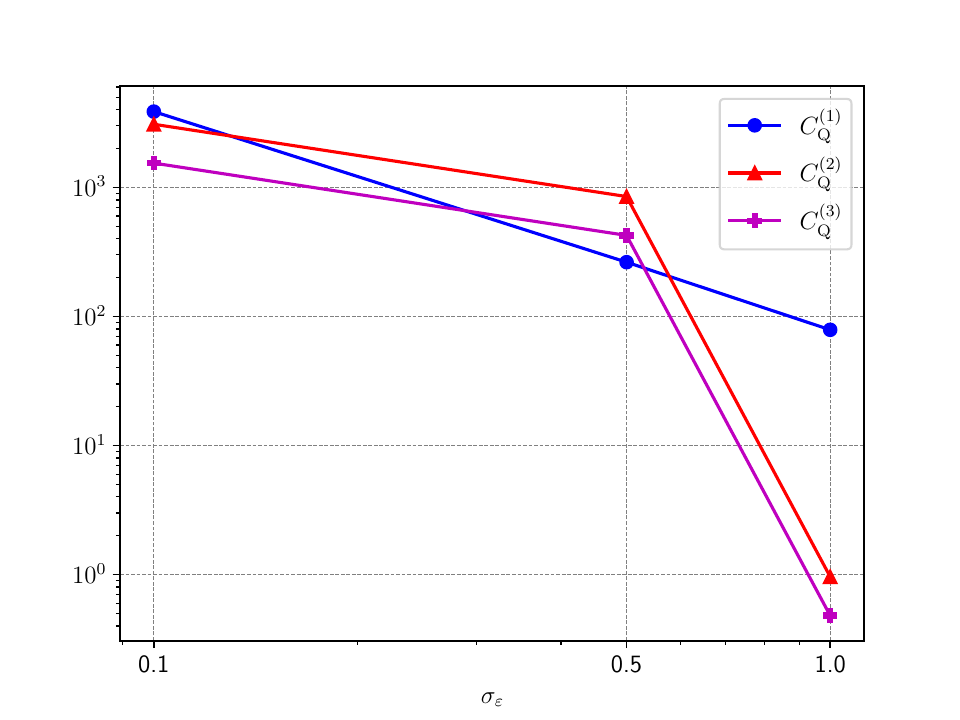}
	}
	\caption{Example 2 (linear EIG example): results of the pilot runs to estimate the outer variance, inner variance, and bias, as well as the corresponding constants $C_{\rm{Q}}^{(1)}, C_{\rm{Q}}^{(2)}, C_{\rm{Q}}^{(3)}$ as a function of the noise level $\sigma_{\varepsilon}$. Bias and inner variance only differ from each other by a multiplicative constant. The same random seeds were used for all pilot runs.}
	\label{fig:ex0.5.noise}
\end{figure}

To analyze the effect of the noise level on the rQMC error, we ran the pilots with truncation for $\sigma_{\varepsilon}\in\{1.0,0.5,0.1\}$. Figure~\ref{fig:ex0.5.noise} Panel (A) presents the outer variance, Panel (B) the inner variance, and Panel (C) the bias. The optimal convergence rates are observed for each setting. In Panel (D), the dependence of the multiplicative constants $C_{\rm{Q}}^{(1)}$, $C_{\rm{Q}}^{(2)}$, and $C_{\rm{Q}}^{(3)}$ on $\sigma_{\varepsilon}$ is demonstrated. These constants increase as $\sigma_{\varepsilon}$ approaches zero. For $\sigma_{\varepsilon}<0.1$, we observed numerical underflow, indicating that importance sampling based on the Laplace approximation becomes necessary for this setting.
\FloatBarrier

\subsection{Example 3: Thermo-mechanics example with inexact sampling}
In this example, $\bs{G}$ is the solution operator of the PDE described below; therefore, we applied an appropriate finite element approximation $\bs{G}_h$. The domain is given by $\cl{D}=[0,1]^2\setminus\cl{B}$, where $\cl{B}$ is a ball centered at $\bs{0}$ with a radius of 0.1. The problem is time-dependent, with observations occurring in $t\in[0,T]$. We solved the following problem for the unknown absolute temperature $\vartheta$ (not to be confused with the dummy variable presented in \eqref{eq:likelihood}) and unknown displacement $\bs{u}$ for fully coupled thermomechanical fields (for more information and a staggered algorithm, see \cite{Far91}).

A weak form of the heat equation is given as follows:
\begin{equation}\label{eq:weak.temperature}
    \int_{\cl{D}}\rho \vartheta_0\dot{s}\hat{\vartheta}\di{}\cl{D}-\int_{\cl{D}}\bs{q}\cdot\nabla\hat{\vartheta}\di{}\cl{D}=-\int_{\partial\cl{D}}\bs{q}\cdot \bs{n} \hat{\vartheta}\di{}S \quad \forall \hat{\vartheta}\in V_{\vartheta},
\end{equation}
where
\begin{equation}
\rho \vartheta_0\dot{s}=\rho \theta_2\dot{\vartheta} + \theta_1(3\lambda + 2\mu) \vartheta_0\tr(\dot{\bs{e}}),
\end{equation}
and
\begin{equation}
    \bs{q}=-\theta_3\nabla{\vartheta}.
\end{equation}

The parameters of interest are $\bs{\theta}=(\theta_1,\theta_2,\theta_3) \in\bs{\Theta}\subset\bb{R}^3$; $\theta_1$, the thermal expansion coefficient; $\theta_2$, the specific heat per unit volume at constant strain; and $\theta_3$, the thermal conductivity. Furthermore, the fixed parameters are the material density $\rho=2700\si{\kg}/\si{\m}^3$, entropy per unit of mass in the current $s$, Lam\'{e} constants $\lambda=7000/1.56$ and $\mu=35000/1.3$, and initial temperature $\vartheta_0=293\si{\kelvin}$. Finally, the strain tensor $\bs{e}=\nabla^{\rm{s}}\bs{u}$, the symmetric part of the gradient of $\bs{u}$,
\begin{equation}
\nabla^{\rm{s}}\bs{u}=\frac{1}{2}(\nabla\bs{u}+(\nabla\bs{u})^{\trans}),
\end{equation}
the unit outward normal $\bs{n}$, and the function space for the temperature field $V_{\vartheta}$ for all $t\in[0,T]$ were also assessed. For the time derivatives, the implicit Euler scheme with a log-spaced step size was applied. A weak form of the momentum equation is
\begin{equation}\label{eq:weak.mechanical}
    \int_{\cl{D}}(\lambda \tr(\bs{e})\bs{1} + 2\mu\bs{e} - \theta_1(3\lambda + 2\mu)(\vartheta-\vartheta_0)\bs{1})\colon\nabla_s \bs{\hat{u}}\di{}\cl{D}=W_{\rm{ext}}(\bs{\hat{u}}) \quad \forall \bs{\hat{u}}\in V_U,
\end{equation}
where $V_U$ denotes the function space for the displacement, $\bs{1}$ indicates the identity tensor, and $W_{\rm{ext}}$ represents the linear functional corresponding to the work of body external forces (neglecting inertia effects; thus, a quasi-static but transient model) and surface tractions on the boundary. The unknown absolute temperature $\vartheta$ was replaced with the temperature variation $\Delta=\vartheta-\vartheta_0$. For this experiment, a temperature increase of $\Delta=10\si{\degreeCelsius}$ was applied at the circular exclusion. Stress and flux-free conditions were applied at the remaining boundaries, and symmetry conditions were applied on the corresponding symmetry planes (for details on the FEM, see Appendix~\ref{ap:FEM}).

The discretization parameter $h$ affects the problem formulation as follows. The mesh-resolution parameter $h_{\rm{mesh}}$ of the FEM was set to $h_{\rm{mesh}}=\lceil 1/ h\rceil$, where $\lceil\cdot\rceil$ is the ceiling operator. Furthermore, the number of time steps was set to $h_{\rm{time}}=\lceil 1/ h^2\rceil$, because the solution to the heat equation \eqref{eq:weak.temperature} depends linearly on time and quadratically on space. Decoupling the parameters $h_{\rm{mesh}}$ and $h_{\rm{time}}$ would enable using multi-index MC techniques \cite{haji2016multi}. 
The actual observation times were selected on a log scale because the problem can be stiff initially. In addition, as the problem approaches a steady state, the observations become increasingly similar, leading to numerical issues unless the observation times are sufficiently spread apart.

\begin{figure}[ht]
	\subfloat[$t=1000$\label{subfig-1:comp.small}]{%
		\includegraphics[width=0.45\textwidth,trim={17cm 5cm 17cm 5cm},clip]{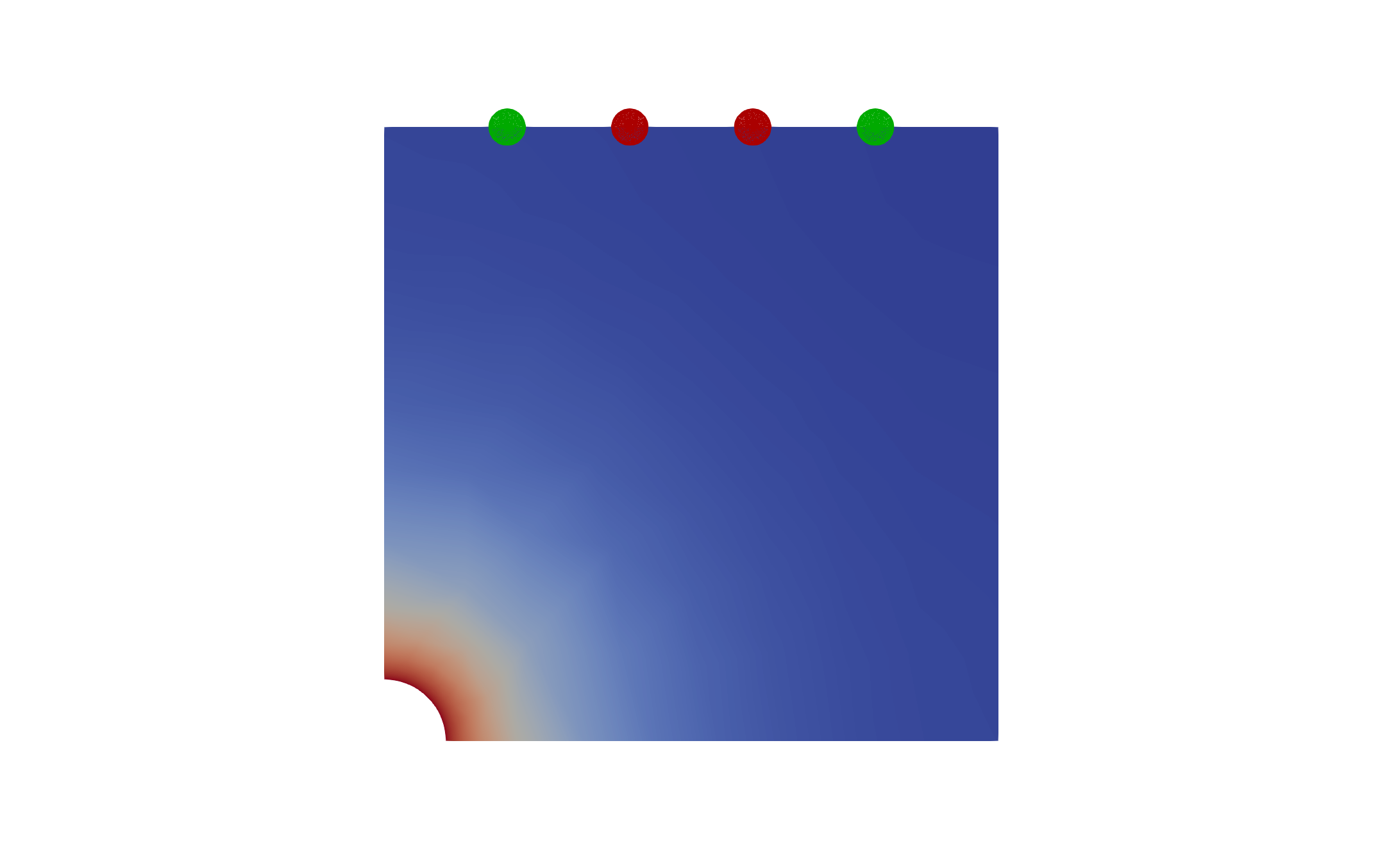}
	}
	\subfloat[$t=10000$\label{subfig-1:comp.int}]{%
		\includegraphics[width=0.45\textwidth,trim={17cm 5cm 17cm 5cm},clip]{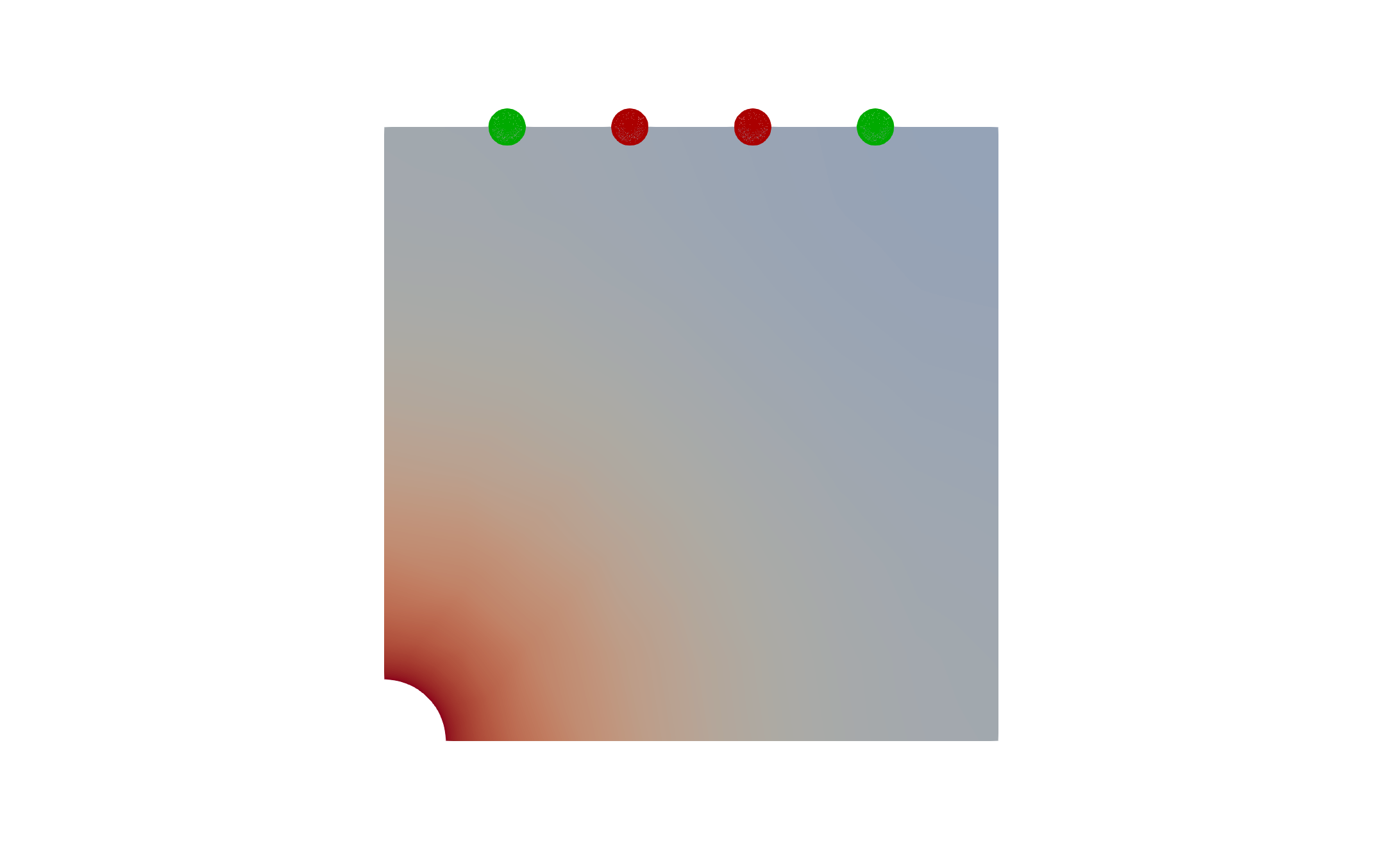}
	}
	\caption{Example 3 (thermo-mechanics example): temperature increases at observation times of the experiment for the rectangular domain with a circular exclusion. Point-sensor locations are displayed as circles at the top of the domain.}
	\label{fig:thermoelast}
\end{figure}

The following mutually independent prior distributions were assumed for the parameters of interest ($\bs{\theta}$):
\begin{equation}
	\theta_1\sim\cl{U}\left(1.81, 2.81\right)\times 10^{-5}, \quad \theta_2\sim\cl{U}\left(8.6,  9.6\right)\times 10^{-4}, \quad \theta_3\sim\cl{U}\left(1.87,  2.87\right)\times 10^{-4}.
\end{equation}
For the observation noise, we assumed that $\bs{\varepsilon}_i\sim\cl{N}(\bs{0},\bs{\Sigma}_{\bs{\varepsilon}})$, where $1\leq i\leq N_e$. Moreover, $\bs{\Sigma}_{\bs{\varepsilon}}$ represents a diagonal matrix in $\bb{R}^{d_y\times d_y}$, where the first $d_y/2$ diagonal entries corresponding to the strain measurements are $\sigma_{\bs{\varepsilon},1}^2=10^{-10}$, and the second $d_y/2$ diagonal entries corresponding to the measurements of the increase in temperature are denoted by $\sigma_{\bs{\varepsilon},2}^2=1/4$. The data model, as presented in \eqref{eq:data.model}, is as follows:
\begin{equation}\label{eq:ex.2}
    \bs{Y}_i(\bs{\xi})=\bs{G}_h(\bs{\theta}_t,\bs{\xi})+\bs{\varepsilon}_i, \quad 1\leq i\leq N_e.
\end{equation}
 We considered two measurements of the strain and two measurements of the temperature increase at three observation times each, yielding $d_y=12$ observations. For $N_e=1$ experiment, the outer integration in \eqref{eq:EIG} is $N_e\times d_y+d_{\theta}=15$ dimensional, and the inner integration is $d_{\theta}=3$ dimensional. 
 Figure~\ref{fig:thermoelast} presents the domain $\cl{D}$, and a circular exclusion occurs in the lower left corner. For the design parameter $\bs{\xi}=(\xi_1, \xi_2)$, $\xi_1$ represents the position of the sensors, and $\xi_2$ indicates the maximum observation time. The sensors for the strain are located at $(0.2,1.0)$ and $(0.8,1.0)$. The sensors for the temperature increase are located at $(0.4,1.0)$ and $(0.6,1.0)$. The temperature increase at the observation times $t=\{10^3, 10^4\}$ is depicted in Panel (A) and Panel (B), respectively.


\begin{figure}[ht]
	\subfloat[Outer variance]{%
		\includegraphics[width=0.45\textwidth]{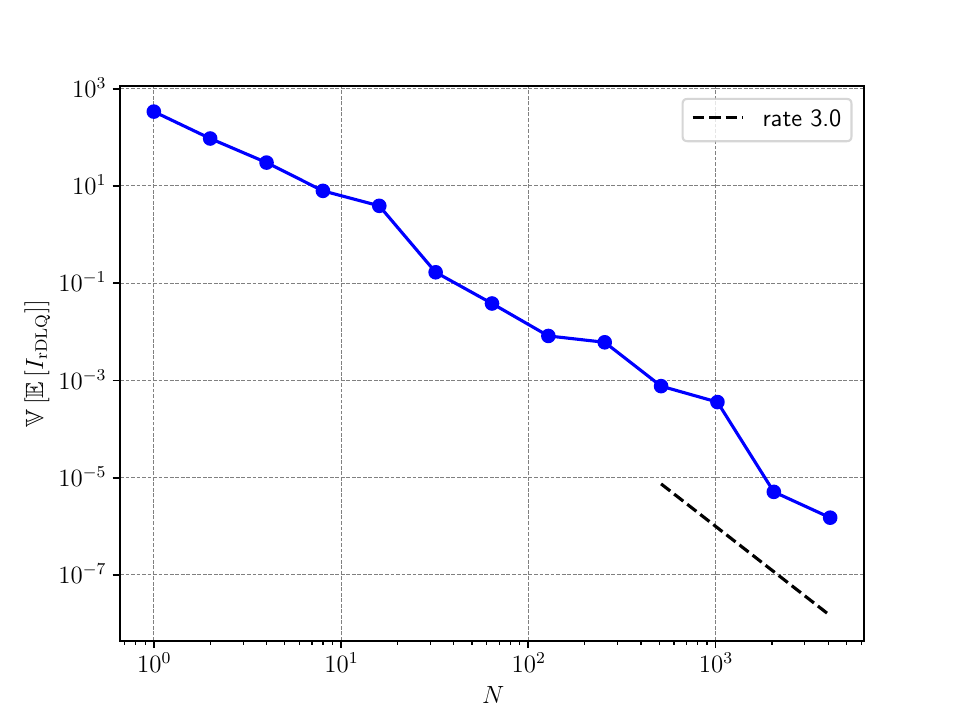}
	}
    \subfloat[Inner variance]{%
		\includegraphics[width=0.45\textwidth]{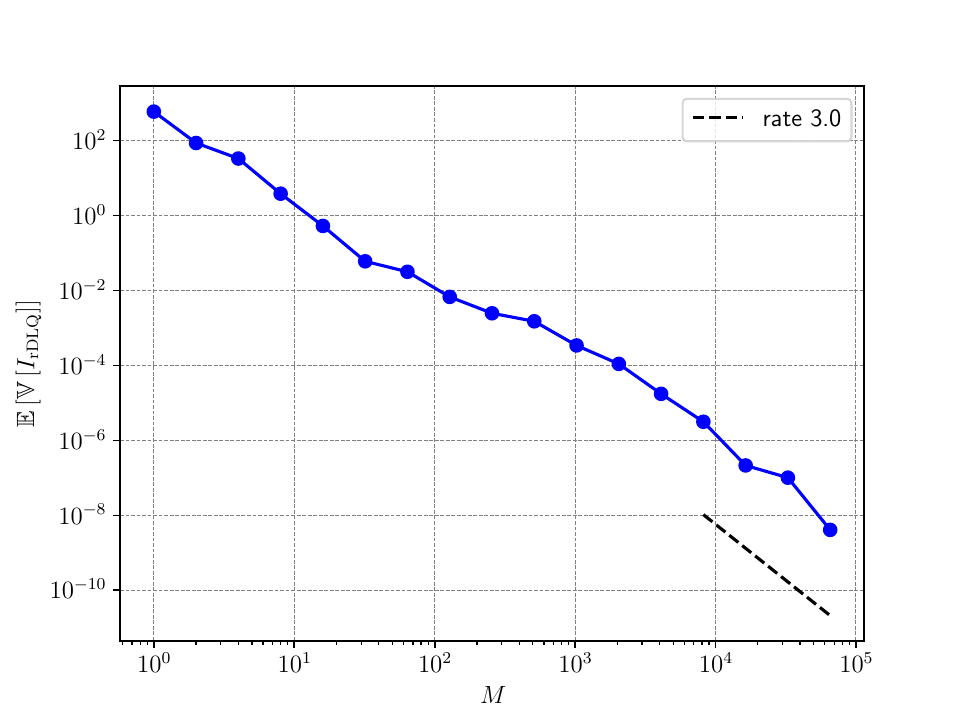}
	}\\
    \subfloat[Bias]{%
		\includegraphics[width=0.45\textwidth]{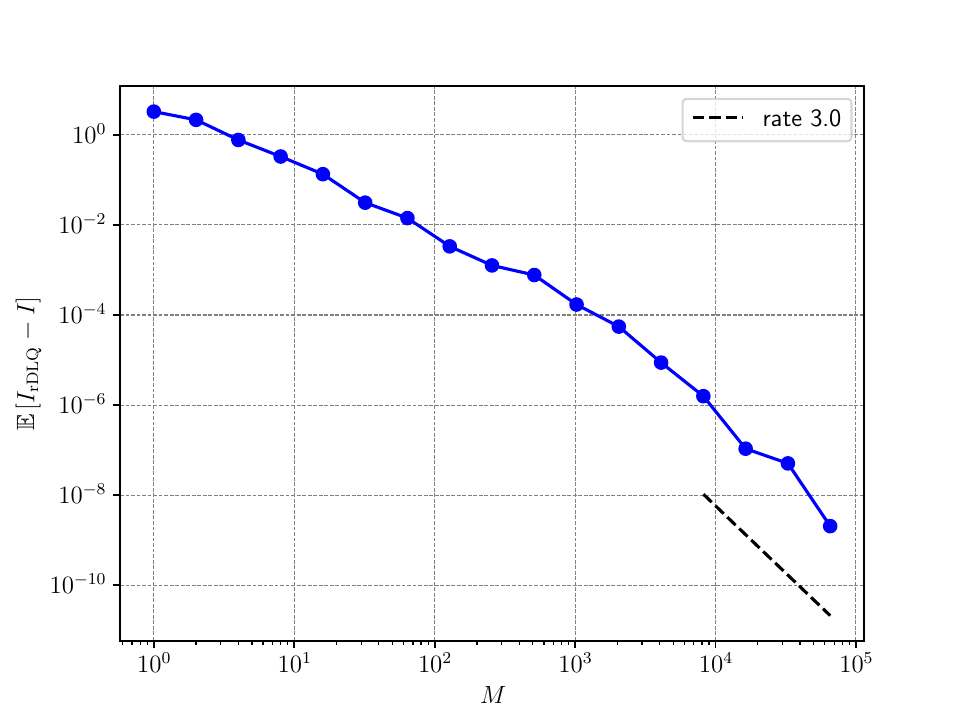}
	}
	\caption{Example 3 (thermo-mechanics example): results of the pilot runs to estimate the outer variance, inner variance, and bias. The bias is proportional to the inner variance. The same random seeds were used for all pilot runs.}
	\label{fig:ex2.pilot}
\end{figure}

Figure~\ref{fig:ex2.pilot} presents the results of the pilot runs for $S=2^5$, $N=2^{12}$, $R=1$, $M=1$ to estimate the outer variance (Panel (A)) and for $S=1$, $N=1$, $R=2^4$, $M=2^{16}$ for the inner variance (Panel (B)) and bias (Panel (C)). The outer variance decreased approximately at the expected rate three ($1+\beta\approx 3.119$), although a larger number of randomizations would be required to obtain clearer results. Due to budget constraints, only a relatively low number of randomizations was chosen for the pilots in this example. The inner variance and bias decreased at rate $1+\delta\approx2.984$. These results are significantly better than the rate one obtained for MC methods or the rate two for rQMC methods not based on Sobol' points with Owen's scrambling. 
The observed convergence rates, multiplicative constants, and total computational work including the finite element contribution are summarized in Table~\ref{table:TE}. The overall cost implied by these results for the rDLQMC estimator is $\cl{O}(TOL^{-(\max\left\{1,2/(1+\beta)\right\}+1/(1+\delta)+\gamma/\eta)})=\cl{O}(TOL^{-(1.335+1.078)})$ where the factor $\gamma/\eta\approx1.078$ comes from the finite element approximation. This rate is consistent with the theoretically optimal rate of $\cl{O}(TOL^{-(4/3+1.078)})$ within the accuracy of the pilot fit and thus lower than the cost of the DLMC estimator, that is, $\cl{O}(TOL^{-(3+1.078)})$ for this example. If other constructions of rDLQMC estimators were used, a rate of $\cl{O}(TOL^{-(3/2+1.078)})$ could be expected. Table~\ref{table:TE} presents the results of the pilot run for the thermo-mechanics example and the resulting optimal computational work compared to the optimal work for the DLMC estimator.

\begin{table}[ht]
\centering
\begin{tabular}{ cccccccc} 
 \toprule
$\beta$ & $\delta$ & $C_{\rm{Q},\epsilon}^{(1)}$ & $C_{\rm{Q},\epsilon}^{(2)}$ & $C_{\rm{Q},\epsilon}^{(3)}$ & $\gamma/\eta$ & $W^{\ast}_{\rm{rDLQ}}$ & $W^{\ast}_{\rm{DLMC}}$ \\ 
 \midrule
 $2.119$ & $1.984$ & $2.6\cdot 10^5$ & $1.4\cdot 10^6$ & $6.9\cdot10^5$ & 1.078 & $\cl{O}(TOL^{-(1.335+1.078)})$ & $\cl{O}(TOL^{-(3+1.078)})$\\ 
 \bottomrule
\end{tabular}
\caption{Example 3 (thermo-mechanics example): the observed rate of increase of the outer variance $1+\beta\approx3.119$ was slightly larger than the rate 3 predicted by the theory. The observed rate of increase of the inner variance and bias $1+\delta\approx2.984$ was slightly lower than the rate 3 predicted by the theory. These variations are likely a result of the low number of randomizations. The total computational work of the rDLQMC estimator, that is, $\cl{O}(TOL^{-(\max\left\{1,2/(1+\beta)\right\}+1/(1+\delta)+\gamma/\eta)})=\cl{O}(TOL^{-(1.335+1.078)})$ is significantly reduced compared to the DLMC estimator as a result of the improved inner and outer integral estimation.}
\label{table:TE}
\end{table}

\FloatBarrier

\section{Conclusion}\label{sec:conclusion}
This paper proposes estimating nested integrals using randomized quasi-Monte Carlo (rQMC) approximations for outer and inner integrals, delivering the double-loop randomized quasi-Monte Carlo estimator (rDLQMC). For this method, we derived asymptotic error bounds and attained the near-optimal number of samples required for each approximation under certain regularity assumptions as the primary contribution of this work. The convergence rates of the inner and outer integral approximations thus improved compared to Monte Carlo estimators. Particular constructions of rQMC estimators based on scrambled Sobol' sequences can achieve even faster error convergence in certain cases of highly regular integrands. The application to Bayesian optimal experimental design indicates that the proposed rDLQMC method is superior to using the rQMC method for only one or neither integral approximation. A truncation scheme for the observation noise made the general error bounds applicable to the estimation of the EIG of an experiment. Error bounds based on the Hardy--Krause variation (implying $L^1$ regularity of the mixed first-order derivatives) displayed only mild dependence on the size of the truncation region, whereas error bounds based on the generalized Hardy--Krause variation of order one (implying $L^2$ regularity) displayed dependence on the truncation region that is exponential in the dimension of the outer integrand. The inner integrand displayed sufficient $L^1$ and $L^2$ regularity, even without truncation. An improved convergence rate for the inner integral approximation error was also observed in our numerical experiments for the EIG application when using scrambled Sobol' sequences, whereas the outer integral approximation error was reduced by a multiplicative factor. A bound on the resulting truncation error was derived for additive Gaussian noise, and the resulting increase in computational cost was analyzed. For future work, the combination of Laplace-based importance sampling with rQMC seems promising. Combining the multilevel MC (or the even more advanced multi-index MC) method with the rQMC method for the outer integral approximation is another promising avenue for future research.

\section{Statements and Declarations}
\subsection{Conflict of interest}
The authors have no conflicts to disclose.

\subsection{Data availability}
The data that support the findings of this study are available from the corresponding author upon reasonable request.



\appendix

\section{Variance inequality}\label{app:variance}
We derived the following inequality for the variance of the sum of dependent random variables $X^{(j)}$, where $1\leq j\leq J$, possibly with different distributions, using the Cauchy--Schwarz inequality:
\begin{align}
    \bb{V}\left[\sum_{j=1}^JX^{(j)}\right]={}&\bb{E}\left[\left(\sum_{j=1}^JX^{(j)}-\bb{E}\left[\sum_{j=1}^JX^{(j)}\right]\right)^2\right],\nonumber\\
    ={}&\bb{E}\left[\left(\sum_{j=1}^JX^{(j)}-\sum_{j=1}^J\bb{E}\left[X^{(j)}\right]\right)^2\right],\nonumber\\
    ={}&\bb{E}\left[\left(\sum_{j=1}^J\left(X^{(j)}-\bb{E}\left[X^{(j)}\right]\right)\right)^2\right],\nonumber\\
    ={}&\bb{E}\left[\left(\sum_{j=1}^J1\cdot\left(X^{(j)}-\bb{E}\left[X^{(j)}\right]\right)\right)^2\right],\nonumber\\
    \leq{}&\bb{E}\left[\sum_{j=1}^J1^2\cdot\sum_{j=1}^J\left(X^{(j)}-\bb{E}\left[X^{(j)}\right]\right)^2\right],\nonumber\\
    ={}&J\sum_{j=1}^J\bb{E}\left[\left(X^{(j)}-\bb{E}\left[X^{(j)}\right]\right)^2\right],\nonumber\\
    ={}&J\sum_{j=1}^J\bb{V}\left[X^{(j)}\right].
\end{align}
If the variables $X^{(j)}$ are independent, then the sharper identity $\bb{V}[\sum_{j=1}^JX^{(j)}]=\sum_{j=1}^J\bb{V}[X^{(j)}]$ holds.

\section{Proof of Remark \ref{rmk:num}}\label{ap:Rmk}
\begin{proof}
First, we assumed that $N_e=d_y=1$. Next, from \eqref{eq:likelihood}, the numerator in the logarithm in \eqref{eq:EIG} is given by the following:
\begin{align}\label{eq:num}
    &\int_{\bb{R}}\log\left(\frac{1}{\sqrt{2\pi\sigma_{\varepsilon}^2}}\exp\left(-\frac{\varepsilon^2}{2\sigma_{\varepsilon}^2}\right)\right)\frac{1}{\sqrt{2\pi\sigma_{\varepsilon}^2}}\exp\left(-\frac{\varepsilon^2}{2\sigma_{\varepsilon}^2}\right)\di{}\varepsilon\nonumber\\
    =&-\frac{1}{2}\log\left(2\pi\sigma_{\varepsilon}^2\right)-\frac{1}{\sqrt{2\pi\sigma_{\varepsilon}^2}}\frac{1}{2\sigma_{\varepsilon}^2}\int_{\bb{R}}\varepsilon^2\exp\left(-\frac{\varepsilon^2}{2\sigma_{\varepsilon}^2}\right)\di{}\varepsilon.
\end{align}
We define $s\coloneqq1/(2\sigma_{\varepsilon}^2)$ and solve the integral in \eqref{eq:num} by taking the derivative with respect to $s$:
\begin{align}
    \int_{\bb{R}}\varepsilon^2e^{-s\varepsilon^2}\di{}\varepsilon&=-\int_{\bb{R}}\frac{\partial}{\partial s}e^{-s\varepsilon^2}\di{}\varepsilon=-\frac{\partial}{\partial s}\int_{\bb{R}}e^{-s\varepsilon^2}\di{}\varepsilon=-\sqrt{\pi}\frac{\partial}{\partial s}s^{-\frac{1}{2}}\nonumber\\
    &=\sqrt{\frac{\pi}{s}}\frac{1}{2s}=\sqrt{2\pi\sigma_{\varepsilon}^2}\sigma_{\varepsilon}^2.
\end{align}
Inserting this into \eqref{eq:num} changes the numerator in \eqref{eq:EIG} to $-(\log(2\pi\sigma_{\varepsilon}^2)+1)/2$. The case for $\bs{\varepsilon}_i\in\bb{R}^{d_y}$, where $1\leq i\leq N_e$, follows trivially under the hypothesis.
\end{proof}

\section{Verification of Assumption~\ref{eq:boundary.growth} for the EIG}\label{app:cond}
To verify Assumption~\ref{eq:boundary.growth} for the nested integrand in the EIG \eqref{eq:EIG}, we introduce the notation 
\begin{equation}
    \bs{\tilde{G}}(\cdot)\coloneqq \bs{G}( F_{\bs{\theta}}^{-1}(\cdot),\bs{\xi}),
\end{equation}
where $F_{\bs{\theta}}^{-1}$ is the inverse CDF of the parameters of interest $\bs{\theta}$ and the dependence on the design $\bs{\xi}$ has been omitted and
state the following:
\begin{asu}[Conditions on the experiment model]\label{asu:Lipschitz}
    Let $\bs{\tilde{G}}:[0,1]^{d_2}\to\bb{R}^{d_1-d_2}$ and assume that there exists $0\leq L<\infty$ such that    \begin{equation}\label{eq:mixed.derivatives}
    \left\lVert\left(\prod_{j\in u}\frac{\partial^u}{\partial z_j}\right)\bs{\tilde{G}}(\bs{z})\right\rVert_2\leq L
\end{equation}
for all $\bs{z}\in[0,1]^{d_{2}}$, and all $u\subseteq\{1,\ldots,d_2\}$. Moreover, let $\bs{\tilde{G}}$ be Lipschitz continuous, that is,
    \begin{equation}\label{eq:Lipschitz}
        \left\lVert \bs{\tilde{G}}(\bs{z})-\bs{\tilde{G}}(\bs{z}^{\prime})\right\rVert_{2}\leq L\left\lVert\bs{z}-\bs{z}^{\prime}\right\rVert_{2},
    \end{equation}
    for $\bs{z},\bs{z}^{\prime}\in[0,1]^{d_{2}}$.
\end{asu}
Lipschitz continuity \eqref{eq:Lipschitz} of the function $\bs{\tilde{G}}$ is a consequence of the differentiability condition \eqref{eq:mixed.derivatives}, and was stated for exposition.
\begin{lem}[Verification of Assumption~\ref{eq:boundary.growth} for the EIG]\label{lem:boundary.condition}
Let $f\equiv\log$ and $g:[0,1]^{d_{1}}\times [0,1]^{d_2}\to\bb{R}_{+}$, where
\begin{equation}\label{def:inner.integrand}
    g(\bs{y},\bs{x})= \frac{1}{\det(2\pi\bs{\Sigma}_{\bs{\varepsilon}})^{\frac{1}{2}}}e^{-\frac{1}{2}\left\lVert \bs{\tilde{G}}(\bs{y}_1)+\bs{\Sigma}_{\bs{\varepsilon}}^{\frac{1}{2}}\Phi^{-1}(\bs{y}_2)-\bs{\tilde{G}}(\bs{x})\right\rVert^2_{\bs{\Sigma}_{\bs{\varepsilon}}^{-1}}},
\end{equation}
where $d_1>d_2$, $\bs{y}=(\bs{y}_1,\bs{y}_2)\in[0,1]^{d_1}$, $\bs{y}_1$, $\bs{x}\in[0,1]^{d_2}$, $\bs{y}_2\in[0,1]^{d_1-d_2}$, and $\bs{\tilde{G}}:[0,1]^{d_2}\to\bb{R}^{d_1-d_2}$ satisfies Assumption~\ref{asu:Lipschitz}. Moreover, $\bs{\Sigma}_{\bs{\varepsilon}}\in\bb{R}^{(d_1-d_2)\times(d_1-d_2)}$ is a positive definite matrix and $\Phi^{-1}$ is the inverse CDF of the standard normal in dimension $d_1-d_2$ and define:
\begin{equation}
    \tilde{f}\coloneqq f\left(\int_{[0,1]^{d_2}}g(\bs{y},\bs{x})\di{}\bs{x}\right).
\end{equation}
    There exists $0<b<\infty$ such that
    \begin{equation}
    \left|\left(\prod_{j\in u}\frac{\partial^u}{\partial y_j}\right)\tilde{f}(\bs{y})\right|\leq b\prod_{i=1}^{d_1}\min (y_i,1-y_i)^{-A_i-\mathds{1}_{\{i\in u\}}},
\end{equation}
for any $A_i>0$, where $1\leq i\leq d_1$, for all $\bs{y}\in(0,1)^{d_1}$ and all $u\subseteq\{1,\ldots,d_1\}$.
    \end{lem}
    \begin{proof}
For illustration, we assumed that $d_{1}=2$, $d_2=1$, and $\bs{\Sigma}_{\bs{\varepsilon}}=1$, which results in
\begin{equation}\label{eq:g.1d}
    g(\bs{y},x)=\frac{1}{\sqrt{2\pi}}e^{-\frac{1}{2}\left(\tilde{G}(y_1)+\Phi^{-1}(y_2)-\tilde{G}(x)\right)^2}.
\end{equation}
We now demonstrate that
\begin{align}
    \left|\frac{\partial^2}{\partial y_1\partial y_2}\tilde{f}(y_1,y_2)\right|{}&\leq b\min(y_1,1-y_1)^{-A_1-1}\min(y_2,1-y_2)^{-A_2-1},\label{eq:f.3}\\
    \left|\frac{\partial}{\partial y_1}\tilde{f}(y_1,y_2)\right|{}&\leq b\min(y_1,1-y_1)^{-A_1-1}\min(y_2,1-y_2)^{-A_2},\label{eq:f.1}\\
    \left|\frac{\partial}{\partial y_2}\tilde{f}(y_1,y_2)\right|{}&\leq b\min(y_1,1-y_1)^{-A_1}\min(y_2,1-y_2)^{-A_2-1},\label{eq:f.2}\\
    \left|\tilde{f}(y_1,y_2)\right|{}&\leq b\min(y_1,1-y_1)^{-A_1}\min(y_2,1-y_2)^{-A_2},\label{eq:f.0}
\end{align}
where
\begin{equation}
    \tilde{f}(y_1,y_2)\equiv\log\left(\frac{1}{\sqrt{2\pi}}\int_{[0,1]}e^{-\frac{1}{2}\left(\tilde{G}(y_1)+\Phi^{-1}(y_2)-\tilde{G}(x)\right)^2}\di{}x\right).
\end{equation}
First, to verify Condition \eqref{eq:f.1}, it follows that
\begin{align}\label{eq:f.y1}
    \left|\frac{\partial}{\partial y_1}\tilde{f}(y_1,y_2)\right|{}&=\left|\frac{\int_{[0,1]}\left(\tilde{G}(y_1)+\Phi^{-1}(y_2)-\tilde{G}(x)\right)\frac{\di{}}{\di{} y_1}\tilde{G}(y_1)e^{-\frac{1}{2}\left(\tilde{G}(y_1)+\Phi^{-1}(y_2)-\tilde{G}(x)\right)^2}\di{}x}{\int_{[0,1]}e^{-\frac{1}{2}\left(\tilde{G}(y_1)+\Phi^{-1}(y_2)-\tilde{G}(x)\right)^2}\di{}x}\right|,\nonumber\\
    {}&\leq\max_{x\in[0,1]}\left\{\left|\left(\tilde{G}(y_1)+\Phi^{-1}(y_2)-\tilde{G}(x)\right)\frac{\di{}}{\di{} y_1}\tilde{G}(y_1)\right|\right\},\nonumber\\
    {}&\leq L\left(L+\left|\Phi^{-1}(y_2)\right|\right)
\end{align}
by Assumption~\ref{asu:Lipschitz}. For the case in which $y_2\to0^{+}$, the following approximation is employed \cite{Liu23, Owe06, Pat96}:
\begin{align}\label{eq:CDF.approx.0}
    \Phi^{-1}(y_2){}&=-\sqrt{-2\log(y_2)}+o(1),
\end{align}
for which it holds that
\begin{equation}\label{eq:inverse.CDF.ineq}
    \left|\Phi^{-1}(y_2)\right|\leq \sqrt{-2\log(y_2)}+P,
\end{equation}
for any $P>0$ and all $y_2\leq y_{2,0}$ for some $y_{2,0}\in(0,1)$ as $y_2\to 0^{+}$. Next, we apply the inequality
\begin{align}\label{eq:sqrt.A}
    \sqrt{2}\sqrt{\log(y_2^{-1})}{}&\leq A_2\log(y_2^{-1}),\nonumber\\
    {}&\leq y_2^{-A_2},
\end{align}
for any $A_2>0$ as $y_2\to 0^{+}$ to obtain
\begin{align}\label{eq:f.y1.2}
    \left|\frac{\partial}{\partial y_1}\tilde{f}(y_1,y_2)\right|{}&\leq L\left(L+\left|\Phi^{-1}(y_2)\right|\right),\nonumber\\
    {}&\leq L\left(L+y_2^{-A_2}+P\right).
\end{align}
The case in which $y_2\to 1^{-}$ follows by a similar derivation, employing the following approximation \cite{Liu23, Owe06, Pat96}:
\begin{equation}\label{eq:CDF.approx.1}
    \Phi^{-1}(y_2)=\sqrt{-2\log(1-y_2)}+o(1).
\end{equation}
Finally,
\begin{equation}
    \left|\frac{\partial}{\partial y_1}\tilde{f}(y_1,y_2)\right|\leq b_1\min(y_1,1-y_1)^{-A_1-1}\min(y_2,1-y_2)^{-A_2}
\end{equation}
for some $0<b_1<\infty$ and any $A_1,A_2>0$. 
For Condition \eqref{eq:f.2}, it holds that
\begin{align}\label{eq:f.y2}
    \left|\frac{\partial}{\partial y_2}\tilde{f}(y_1,y_2)\right|{}&=\left|\frac{\int_{[0,1]}\left(\tilde{G}(y_1)+\Phi^{-1}(y_2)-\tilde{G}(x)\right)\frac{\di{}}{\di{} y_2}\Phi^{-1}(y_2)e^{-\frac{1}{2}\left(\tilde{G}(y_1)+\Phi^{-1}(y_2)-\tilde{G}(x)\right)^2}\di{}x}{\int_{[0,1]}e^{-\frac{1}{2}\left(\tilde{G}(y_1)+\Phi^{-1}(y_2)-\tilde{G}(x)\right)^2}\di{}x}\right|,\nonumber\\
    {}&\leq\max_{x\in[0,1]} \left\{\left|\left(\tilde{G}(y_1)+\Phi^{-1}(y_2)-\tilde{G}(x)\right)\frac{\di{}}{\di{} y_2}\Phi^{-1}(y_2)\right|\right\},\nonumber\\
    {}&\leq \left|\frac{\di{}}{\di{} y_2}\Phi^{-1}(y_2)\right|\left(L+\left|\Phi^{-1}(y_2)\right|\right).
\end{align}
It follows that
\begin{align}\label{eq:derivative.inverse.CDF}
    \frac{\di{}}{\di{} y_2}\Phi^{-1}(y_2){}&=\left(\Phi^{\prime}(\Phi^{-1}(y_2))\right)^{-1},\nonumber\\
{}&=\sqrt{2\pi}e^{\frac{1}{2}(\Phi^{-1}(y_2))^2},
\end{align}
yielding
\begin{align}\label{eq:der.inverse.CDF.ineq}
    \frac{\di{}}{\di{} y_2}\Phi^{-1}(y_2){}&=\sqrt{2\pi}e^{\frac{1}{2}(\Phi^{-1}(y_2))^2},\nonumber\\
    {}&=\sqrt{2\pi}e^{\frac{1}{2}\left(-\sqrt{-2\log(y_2)}+o(1)\right)^2},\nonumber\\
    {}&=\sqrt{2\pi}e^{\left(\log(y_2^{-1})+o\left(\sqrt{\log(y_2^{-1})}\right)+o(1)\right)},\nonumber\\
    {}&\leq\sqrt{2\pi}e^{\left(\log(y_2^{-1})+P\sqrt{\log(y_2^{-1})}+P\right)},\nonumber\\
    {}&\leq\sqrt{2\pi}e^{\left(\log(y_2^{-1})+B_2\log(y_2^{-1})+P\right)},\nonumber\\
    {}&=\sqrt{2\pi}e^{P}y_2^{-B_2-1}.
\end{align}
for any $P,B_2>0$ and all $y_2\leq y_{2,0}$ for some $y_{2,0}\in(0,1)$ as $y_2\to 0^{+}$. 
Thus, it follows that
\begin{align}
    \left|\Phi^{-1}(y_2)\right|\left|\frac{\di{}}{\di{} y_2}\Phi^{-1}(y_2)\right|{}&\leq\left(\sqrt{-2\log(y_2)}+P\right)\sqrt{2\pi}e^{P}y_2^{-B_2-1},\nonumber\\
    {}&\leq\left(y_2^{-C_2}+P\right)\sqrt{2\pi}e^{P}y_2^{-B_2-1},
\end{align}
for any $C_2>0$ as $y_2\to 0^{+}$. The case in which $y_2\to 1^{-}$ again follows from the application of the Approximation in \eqref{eq:CDF.approx.1}. Finally,
\begin{equation}
    \left|\frac{\partial}{\partial y_2}\tilde{f}(y_1,y_2)\right|\leq b_2\min(y_1,1-y_1)^{-A_1}\min(y_2,1-y_2)^{-A_2-1}
\end{equation}
for some $0<b_2<\infty$, any $A_1>0$, and any $A_2\coloneqq B_2+C_2>0$.
Condition~\eqref{eq:f.3} follows by a similar derivation for some $0<b_3\,<\infty$ and any $A_1, A_2>0$. For Condition~\eqref{eq:f.0}, Assumption~\ref{asu:Lipschitz} was used to obtain the following:
\begin{align}
    {}&\left|\log\left(\int_{[0,1]}e^{-\frac{1}{2}\left(\tilde{G}(y_1)+\Phi^{-1}(y_2)-\tilde{G}(x)\right)^2}\di{}x\right)\right|\nonumber\\
    {}&=\left|\log\left(\int_{[0,1]}e^{-\frac{1}{2}\left(\left(\tilde{G}(y_1)-\tilde{G}(x)\right)^2+2\left(\tilde{G}(y_1)-\tilde{G}(x)\right)\Phi^{-1}(y_2)+\left(\Phi^{-1}(y_2)\right)^2\right)}\di{}x\right)\right|\nonumber\\
    {}&\leq\left|\log\left(\int_{[0,1]}e^{-\frac{1}{2}\left(\left(\tilde{G}(y_1)-\tilde{G}(x)\right)^2+2\left|\tilde{G}(y_1)-\tilde{G}(x)\right|\left|\Phi^{-1}(y_2)\right|+\left(\Phi^{-1}(y_2)\right)^2\right)}\di{}x\right)\right|\nonumber\\
    {}&\leq\left|\log\left(\int_{[0,1]}e^{-\frac{1}{2}\left(L^2\left(y_1-x\right)^2+2L\left|y_1-x\right|\left|\Phi^{-1}(y_2)\right|+\left(\Phi^{-1}(y_2)\right)^2\right)}\di{}x\right)\right|\nonumber\\
    {}&\leq\left|\log\left(e^{-\frac{1}{2}\left(L^2+2L\left|\Phi^{-1}(y_2)\right|+\left(\Phi^{-1}(y_2)\right)^2\right)}\right)\right|\nonumber\\
    {}&=\left|\frac{1}{2}L^2+L\left|\Phi^{-1}(y_2)\right|+\frac{1}{2}\left(\Phi^{-1}(y_2)\right)^2\right|.
\end{align}
The approximation \eqref{eq:CDF.approx.0} as $y_2\to 0^{+}$ was used to obtain:
\begin{align}
    \left|\frac{1}{2}L^2+L\left|\Phi^{-1}(y_2)\right|+\left(\Phi^{-1}(y_2)\right)^2\right|{}&=\left|\frac{1}{2}L^2+L\left|-\sqrt{-2\log(y_2)}+o(1)\right|+\left(-\sqrt{-2\log(y_2)}+o(1)\right)^2\right|,\nonumber\\
    {}&\leq\left|\frac{1}{2}L^2+L\sqrt{2}\sqrt{\log(y_2^{-1})}+LP+\log(y_2^{-1})+P\sqrt{2}\sqrt{\log(y_2^{-1})}+\frac{1}{2}P\right|,\nonumber\\
    {}&=\left|\frac{1}{2}L^2+(L+P)\sqrt{2}\sqrt{\log(y_2^{-1})}+LP+\log(y_2^{-1})+\frac{1}{2}P\right|,\nonumber\\
    {}&\leq\left|\frac{1}{2}L^2+(L+P)\sqrt{2}y_2^{-C_2}+LP+y_2^{-B_2}+\frac{1}{2}P\right|
\end{align}
for any $B_2,C_2>0$. The case in which $y_2\to 1^{-}$ follows from the application of the Approximation in \eqref{eq:CDF.approx.1}, providing

\begin{align}
    \left|\tilde{f}(y_1,y_2)\right|{}&\leq b_4\min(y_1,1-y_1)^{-A_1}\min(y_2,1-y_2)^{-A_2}
\end{align}
for some $0<b_4<\infty$, any $A_1$, and any $A_2\coloneqq \max\{B_2,C_2\}>0$. Setting $b\coloneqq\max_{i\in\{1,\ldots,4\}}b_i$ demonstrates the Conditions~\eqref{eq:f.3} to~\eqref{eq:f.0}. Moreover, the general case in higher dimensions can be demonstrated using the Fa\`{a} di Bruno formula under the auxiliary assumption that the smallest eigenvalue of the noise covariance matrix is bounded away from zero.
\end{proof}

\begin{rmk}[Dependence of the constant $b$ on the integrand dimension $d_1$]
    The analysis of the multiplicative constant $b$ in the above result using Fa\`{a} di Bruno's formula is beyond the scope of this work; however, we note that $b$ possibly exhibits a strong dependence on the integrand dimension $d_1$.
\end{rmk}

\section{Verification of Condition~(\ref{eq:reg:g}) for the EIG}\label{app:Condition.42}
To verify Condition~\eqref{eq:reg:g}, we assume the same setting as in Lemma~\ref{lem:boundary.condition} and set
\begin{align}\label{eq:truncation.split}
    I={}&\int_{[0,1]^2}\log\left(\frac{1}{\sqrt{2\pi}}\int_{[0,1]}e^{-\frac{1}{2}\left(\tilde{G}(y_1)+\Phi^{-1}(y_2)-\tilde{G}(x)\right)^2}\di{}x\right)\di{}y_2\di{}y_1,\nonumber\\
    ={}&\int_{[0,1]}\int_{\bb{R}}\log\left(\frac{1}{\sqrt{2\pi}}\int_{[0,1]}e^{-\frac{1}{2}\left(\tilde{G}(y_1)+\varepsilon-\tilde{G}(x)\right)^2}\di{}x\right)\frac{1}{\sqrt{2\pi}}e^{-\frac{1}{2}\varepsilon^2}\di{}\varepsilon\di{}y_1,\nonumber\\
    ={}&\int_{[0,1]}\int_{-c(TOL)}^{c(TOL)}\log\left(\frac{1}{\sqrt{2\pi}}\int_{[0,1]}e^{-\frac{1}{2}\left(\tilde{G}(y_1)+\varepsilon-\tilde{G}(x)\right)^2}\di{}x\right)\frac{1}{\sqrt{2\pi}}e^{-\frac{1}{2}\varepsilon^2}\di{}\varepsilon\di{}y_1,\nonumber\\
    {}&+\int_{[0,1]}\int_{-\infty}^{-c(TOL)}\log\left(\frac{1}{\sqrt{2\pi}}\int_{[0,1]}e^{-\frac{1}{2}\left(\tilde{G}(y_1)+\varepsilon-\tilde{G}(x)\right)^2}\di{}x\right)\frac{1}{\sqrt{2\pi}}e^{-\frac{1}{2}\varepsilon^2}\di{}\varepsilon\di{}y_1,\nonumber\\
    {}&+\int_{[0,1]}\int_{c(TOL)}^{\infty}\log\left(\frac{1}{\sqrt{2\pi}}\int_{[0,1]}e^{-\frac{1}{2}\left(\tilde{G}(y_1)+\varepsilon-\tilde{G}(x)\right)^2}\di{}x\right)\frac{1}{\sqrt{2\pi}}e^{-\frac{1}{2}\varepsilon^2}\di{}\varepsilon\di{}y_1,\nonumber\\
    ={}&I^{\rm{tr}}+\int_{[0,1]}\int_{-\infty}^{-c(TOL)}\log\left(\frac{1}{\sqrt{2\pi}}\int_{[0,1]}e^{-\frac{1}{2}\left(\tilde{G}(y_1)+\varepsilon-\tilde{G}(x)\right)^2}\di{}x\right)\frac{1}{\sqrt{2\pi}}e^{-\frac{1}{2}\varepsilon^2}\di{}\varepsilon\di{}y_1,\nonumber\\
    {}&+\int_{[0,1]}\int_{c(TOL)}^{\infty}\log\left(\frac{1}{\sqrt{2\pi}}\int_{[0,1]}e^{-\frac{1}{2}\left(\tilde{G}(y_1)+\varepsilon-\tilde{G}(x)\right)^2}\di{}x\right)\frac{1}{\sqrt{2\pi}}e^{-\frac{1}{2}\varepsilon^2}\di{}\varepsilon\di{}y_1,\nonumber\\
\end{align}
where 
\begin{equation}\label{eq:c.TOL}
    c(TOL)\coloneqq (2(1+p))^{\frac{1}{2}}\log(TOL^{-1})^{\frac{1}{2}},
\end{equation}
for any $p>0$, signifies a region of truncation for a tolerance $TOL>0$ and
\begin{align}\label{eq:I.tr}
    I^{\rm{tr}}{}&=\int_{[0,1]}\int_{-c(TOL)}^{c(TOL)}\log\left(\frac{1}{\sqrt{2\pi}}\int_{[0,1]}e^{-\frac{1}{2}\left(\tilde{G}(y_1)+\varepsilon-\tilde{G}(x)\right)^2}\di{}x\right)\frac{1}{\sqrt{2\pi}}e^{-\frac{1}{2}\varepsilon^2}\di{}\varepsilon\di{}y_1,\nonumber\\{}&=(2\Phi(c(TOL))-1)\int_{[0,1]^2}\log\left(\frac{1}{\sqrt{2\pi}}\int_{[0,1]}e^{-\frac{1}{2}\left(\tilde{G}(y_1)+F_{c(TOL)}^{-1}(y_2)-\tilde{G}(x)\right)^2}\di{}x\right)\di{}y_2\di{}y_1,
\end{align}
where 
\begin{equation}\label{eq:inverse.CDF}
F_{c(TOL)}^{-1}(y_2)\coloneqq\Phi^{-1}\Big((2\Phi(c(TOL))-1)y_2+(1-\Phi(c(TOL)))\Big)
\end{equation}
is the inverse CDF of the standard normal distribution truncated to the interval $[-c(TOL),c(TOL)]$.
 \begin{lem}[Verification of Condition~(39) for the EIG]\label{lem:cond}
 Given the assumptions of Lemma~\ref{lem:boundary.condition} for $F_{c(TOL)}^{-1}(\cdot)$ instead of $\Phi^{-1}(\cdot)$, for a fixed and sufficiently large value of $c(TOL)$ as in \eqref{eq:c.TOL}, where $TOL>0$, there exists $0\leq \tilde{L}<\infty$ such that 
 \begin{equation}
    \sup_{\bs{y}\in[0,1]^{2}} \left\lvert \frac{V_{\rm{HK}}(g_{h}(\bs{y},\cdot))}{ \bar{g}_{h}(\bs{y})} \right\rvert \leq \tilde{L}^{d_2}c(TOL)^{d_2}.
\end{equation}
 \end{lem}
 \begin{proof}
 The Hardy--Krause variation is as follows:
\begin{align}
    V_{\rm{HK}}(g_{h}(\bs{y},\cdot)){}&=\int_{[0,1]}\left|\frac{\partial}{\partial x}\frac{1}{\sqrt{2\pi}}e^{-\frac{1}{2}\left(\tilde{G}_{h}(y_1)+F_{c(TOL)}^{-1}(y_2)-\tilde{G}_{h}(x)\right)^2}\right|\di{}x,\nonumber\\
    {}&=\frac{1}{\sqrt{2\pi}}\int_{[0,1]}\left|\left(\tilde{G}_{h}(y_1)+F_{c(TOL)}^{-1}(y_2)-\tilde{G}_{h}(x)\right)\frac{\di{}}{\di{} x}\tilde{G}_{h}(x)e^{-\frac{1}{2}\left(\tilde{G}_{h}(y_1)+F_{c(TOL)}^{-1}(y_2)-\tilde{G}_{h}(x)\right)^2}\right|\di{}x,\nonumber\\
    {}&=\frac{1}{\sqrt{2\pi}}\int_{[0,1]}\left|\frac{\di{}}{\di{} x}\tilde{G}_{h}(x)\right|\left|\left(\tilde{G}_{h}(y_1)+F_{c(TOL)}^{-1}(y_2)-\tilde{G}_{h}(x)\right)e^{-\frac{1}{2}\left(\tilde{G}_{h}(y_1)+F_{c(TOL)}^{-1}(y_2)-\tilde{G}_{h}(x)\right)^2}\right|\di{}x,\nonumber\\
    {}&\leq\frac{1}{\sqrt{2\pi}}\int_{[0,1]}\left|\frac{\di{}}{\di{} x}\tilde{G}_{h}(x)\right|\di{}x,\nonumber\\
    {}&\leq\frac{L}{\sqrt{2\pi}},
\end{align}
For Condition~\eqref{eq:reg:g},
\begin{align}
    \sup_{\bs{y}\in[0,1]^{2}}{}&\left|\frac{\frac{1}{\sqrt{2\pi}}\int_{[0,1]}\left|\left(\tilde{G}_{h}(y_1)+F_{c(TOL)}^{-1}(y_2)-\tilde{G}_{h}(x)\right)\frac{\di{}}{\di{} x}\tilde{G}_{h}(x)e^{-\frac{1}{2}\left(\tilde{G}_{h}(y_1)+F_{c(TOL)}^{-1}(y_2)-\tilde{G}_{h}(x)\right)^2}\right|\di{}x}{\frac{1}{\sqrt{2\pi}}\int_{[0,1]}e^{-\frac{1}{2}\left(\tilde{G}_{h}(y_1)+F_{c(TOL)}^{-1}(y_2)-\tilde{G}_{h}(x)\right)^2}\di{}x}\right|\nonumber\\
    {}&\leq \sup_{\bs{y}\in[0,1]^{2}}\left|\max_{x\in[0,1]}\left\{\left(\tilde{G}_{h}(y_1)+F_{c(TOL)}^{-1}(y_2)-\tilde{G}_{h}(x)\right)\frac{\di{}}{\di{} x}\tilde{G}_{h}(x)\right\}\right|,\nonumber\\
    {}&\leq L\left(L+\sup_{\bs{y}\in[0,1]^{2}}\left|F_{c(TOL)}^{-1}(y_2)\right|\right).
\end{align}
The CDF $\Phi(\cdot)$ of the standard normal and its inverse are monotonic; hence, the supremum of $|F_{c(TOL)}^{-1}(y_2)|$ is attained at $y_2\in\{0,1\}$, yielding
\begin{align}
    \sup_{y_2\in[0,1]}|F_{c(TOL)}^{-1}(y_2)|{}&=c(TOL),
\end{align}
with $c(TOL)$ as in \eqref{eq:c.TOL}. It follows that
\begin{align}
    \sup_{\bs{y}\in[0,1]^{2}}{}&\left|\frac{\frac{1}{\sqrt{2\pi}}\int_{[0,1]}\left|\left(\tilde{G}_{h}(y_1)+F_{c(TOL)}^{-1}(y_2)-\tilde{G}_{h}(x)\right)\frac{\di{}}{\di{} x}\tilde{G}_{h}(x)e^{-\frac{1}{2}\left(\tilde{G}_{h}(y_1)+F_{c(TOL)}^{-1}(y_2)-\tilde{G}_{h}(x)\right)^2}\right|\di{}x}{\frac{1}{\sqrt{2\pi}}\int_{[0,1]}e^{-\frac{1}{2}\left(\tilde{G}_{h}(y_1)+F_{c(TOL)}^{-1}(y_2)-\tilde{G}_{h}(x)\right)^2}\di{}x}\right|\nonumber\\
    {}&\leq \tilde{L}c(TOL)
\end{align}
for $TOL$ sufficiently small, where $\tilde{L}>0$ depends on the Lipschitz constant $L$ of $G_{h}$.
The general case in higher dimensions can also be demonstrated using the Fa\`{a} di Bruno formula.
\end{proof}
The constant $k=\tilde{L}^{d_2}c(TOL)^{d_2}$ in Condition~\eqref{eq:reg:g} thus displays only weak dependence on the error tolerance ($TOL$). The above result was demonstrated for parameters of interest with a uniform prior distribution. For parameters of interest distributed according to a normal distribution, the present inverse-inequality argument does not apply directly because derivatives of $\Phi^{-1}$ introduce additional boundary growth. However, numerical results indicate that the bias and inner variance still display behavior consistent with the results of Proposition~\ref{prop:B.DLQ} and~\ref{prop:V.DLQ}; see Figure~\ref{fig:ex0.5.pilot} Panel~(B) and Panel~(C). The analysis of this special case is thus left for future work.
\section{Bound on the truncation error of the nested integral in the EIG}\label{app:Truncation.EIG}
\begin{lem}[Bound on the truncation error of the nested integral in the EIG]\label{lem:trunc}
Given the assumptions of Lemma~\ref{lem:cond}, the truncation error is bounded as
\begin{align}\label{eq:bound.truncation}
    |I-I^{\rm{tr}}|=o(TOL),
\end{align}
where $I$ as in \eqref{eq:truncation.split} and $I^{\rm{tr}}$ as in \eqref{eq:I.tr},
as $TOL\to 0$.
\end{lem}

\begin{proof}
By the triangle inequality,
\begin{align}\label{eq:tr.error}
    |I-I^{\rm{tr}}|={}&\left|\int_{[0,1]}\int^{-c(TOL)}_{-\infty}\log\left(\frac{1}{\sqrt{2\pi}}\int_{[0,1]}e^{-\frac{1}{2}\left(\tilde{G}(y_1)+\varepsilon-\tilde{G}(x)\right)^2}\di{}x\right)\frac{1}{\sqrt{2\pi}}e^{-\frac{1}{2}\varepsilon^2}\di{}\varepsilon\di{}y_1\right.\nonumber\\
    {}&\left.+\int_{[0,1]}\int_{c(TOL)}^{\infty}\log\left(\frac{1}{\sqrt{2\pi}}\int_{[0,1]}e^{-\frac{1}{2}\left(\tilde{G}(y_1)+\varepsilon-\tilde{G}(x)\right)^2}\di{}x\right)\frac{1}{\sqrt{2\pi}}e^{-\frac{1}{2}\varepsilon^2}\di{}\varepsilon\di{}y_1\right|,\nonumber\\
    \leq{}&\left|\int_{[0,1]}\int^{-c(TOL)}_{-\infty}\log\left(\frac{1}{\sqrt{2\pi}}\int_{[0,1]}e^{-\frac{1}{2}\left(\tilde{G}(y_1)+\varepsilon-\tilde{G}(x)\right)^2}\di{}x\right)\frac{1}{\sqrt{2\pi}}e^{-\frac{1}{2}\varepsilon^2}\di{}\varepsilon\di{}y_1\right|\nonumber\\
    {}&+\left|\int_{[0,1]}\int_{c(TOL)}^{\infty}\log\left(\frac{1}{\sqrt{2\pi}}\int_{[0,1]}e^{-\frac{1}{2}\left(\tilde{G}(y_1)+\varepsilon-\tilde{G}(x)\right)^2}\di{}x\right)\frac{1}{\sqrt{2\pi}}e^{-\frac{1}{2}\varepsilon^2}\di{}\varepsilon\di{}y_1\right|.
\end{align}
By Jensen's inequality it follows that
\begin{align}
    \frac{1}{2}\left(\tilde{G}(y_1)+\varepsilon-\tilde{G}(x)\right)^2{}&=2\left(\frac{(\tilde{G}(y_1)-\tilde{G}(x))}{2}+\frac{\varepsilon}{2}\right)^2,\nonumber\\
    {}&\leq2\left(\frac{(\tilde{G}(y_1)-\tilde{G}(x))^2}{2}+\frac{\varepsilon^2}{2}\right),\nonumber\\
    {}&=(\tilde{G}(y_1)-\tilde{G}(x))^2+\varepsilon^2,
\end{align}
and thus for the last term in \eqref{eq:tr.error} that
\begin{align}
{}&\left|\int_{[0,1]}\int_{c(TOL)}^{\infty}\log\left(\frac{1}{\sqrt{2\pi}}\int_{[0,1]}e^{-\frac{1}{2}\left(\tilde{G}(y_1)+\varepsilon-\tilde{G}(x)\right)^2}\di{}x\right)\frac{1}{\sqrt{2\pi}}e^{-\frac{1}{2}\varepsilon^2}\di{}\varepsilon\di{}y_1\right|\nonumber\\
    {}&\leq\int_{[0,1]}\int_{c(TOL)}^{\infty}\left|\log\left(\frac{1}{\sqrt{2\pi}}\int_{[0,1]}e^{-\frac{1}{2}\left(\tilde{G}(y_1)+\varepsilon-\tilde{G}(x)\right)^2}\di{}x\right)\right|\frac{1}{\sqrt{2\pi}}e^{-\frac{1}{2}\varepsilon^2}\di{}\varepsilon\di{}y_1,\nonumber\\
    {}&\leq\int_{[0,1]}\int_{c(TOL)}^{\infty}\left|\log\left(\frac{1}{\sqrt{2\pi}}\int_{[0,1]}e^{-\left(\tilde{G}(y_1)-\tilde{G}(x)\right)^2-\varepsilon^2}\di{}x\right)\right|\frac{1}{\sqrt{2\pi}}e^{-\frac{1}{2}\varepsilon^2}\di{}\varepsilon\di{}y_1{}.
\end{align}
Next, from Assumption~\ref{asu:Lipschitz} it follows that
\begin{align}\label{eq:truncation.bound}
    {}&\int_{[0,1]}\int_{c(TOL)}^{\infty}\left|\log\left(\frac{1}{\sqrt{2\pi}}\int_{[0,1]}e^{-\left(\tilde{G}(y_1)-\tilde{G}(x)\right)^2-\varepsilon^2}\di{}x\right)\right|\frac{1}{\sqrt{2\pi}}e^{-\frac{1}{2}\varepsilon^2}\di{}\varepsilon\di{}y_1{},\nonumber\\
    {}&\leq\int_{[0,1]}\int_{c(TOL)}^{\infty}\left|\log\left(\frac{1}{\sqrt{2\pi}}\int_{[0,1]}e^{-L^2\left(y_1-x\right)^2-\varepsilon^2}\di{}x\right)\right|\frac{1}{\sqrt{2\pi}}e^{-\frac{1}{2}\varepsilon^2}\di{}\varepsilon\di{}y_1,\nonumber\\
    {}&<\int_{c(TOL)}^{\infty}\left|\log\left(\frac{1}{\sqrt{2\pi}}e^{-L^2-\varepsilon^2}\right)\right|\frac{1}{\sqrt{2\pi}}e^{-\frac{1}{2}\varepsilon^2}\di{}\varepsilon,\nonumber\\
    {}&=\int_{c(TOL)}^{\infty}\left|T-\varepsilon^2\right|\frac{1}{\sqrt{2\pi}}e^{-\frac{1}{2}\varepsilon^2}\di{}\varepsilon,\nonumber\\
    {}&\leq\int_{c(TOL)}^{\infty}\left|T\right|\frac{1}{\sqrt{2\pi}}e^{-\frac{1}{2}\varepsilon^2}\di{}\varepsilon+\int_{c(TOL)}^{\infty}\varepsilon^2\frac{1}{\sqrt{2\pi}}e^{-\frac{1}{2}\varepsilon^2}\di{}\varepsilon,
\end{align}
where $T\coloneqq-L^2-\log(\sqrt{2\pi})$. The first term in \eqref{eq:truncation.bound} is bounded via integration by parts,
\begin{align}\label{eq:T.truncated}
    \int_{c(TOL)}^{\infty}\left|T\right|\frac{1}{\sqrt{2\pi}}e^{-\frac{1}{2}\varepsilon^2}\di{}\varepsilon{}&=|T|\frac{1}{\sqrt{2\pi}}\int_{c(TOL)}^{\infty}\frac{1}{\varepsilon}\cdot\varepsilon e^{-\frac{1}{2}\varepsilon^2}\di{}\varepsilon,\nonumber\\
    {}&=|T|\frac{1}{\sqrt{2\pi}}\left(\left[-\frac{1}{\varepsilon} e^{-\frac{1}{2}\varepsilon^2}\right]_{c(TOL)}^{\infty}-\int_{c(TOL)}^{\infty}\frac{1}{\varepsilon^2}e^{-\frac{1}{2}\varepsilon^2}\di{}\varepsilon\right),\nonumber\\
    {}&=|T|\frac{1}{\sqrt{2\pi}}\left(\frac{1}{c(TOL)} e^{-\frac{1}{2}(c(TOL))^2}-\int_{c(TOL)}^{\infty}\frac{1}{\varepsilon^2}e^{-\frac{1}{2}\varepsilon^2}\di{}\varepsilon\right),\nonumber\\
    {}&<|T|\frac{1}{\sqrt{2\pi}}\frac{1}{c(TOL)} e^{-\frac{1}{2}(c(TOL))^2},
\end{align} 
where the last term in the second to last line is strictly positive. Substituting \eqref{eq:c.TOL} into \eqref{eq:T.truncated}, we find that
\begin{align}
    \int_{c(TOL)}^{\infty}\left|T\right|\frac{1}{\sqrt{2\pi}}e^{-\frac{1}{2}\varepsilon^2}\di{}\varepsilon{}&<|T|\frac{1}{\sqrt{2\pi}}(2(1+p))^{-\frac{1}{2}}\log(TOL^{-1})^{-\frac{1}{2}} TOL^{1+p},\nonumber\\
    {}&=o(TOL)
\end{align}
as $TOL\to 0$. The second term in \eqref{eq:truncation.bound} is bounded via integration by parts as
\begin{align}
    \int_{c(TOL)}^{\infty}\varepsilon^2\frac{1}{\sqrt{2\pi}}e^{-\frac{1}{2}\varepsilon^2}\di{}\varepsilon{}&=\frac{1}{\sqrt{2\pi}}\int_{c(TOL)}^{\infty}\varepsilon\cdot \varepsilon e^{-\frac{1}{2}\varepsilon^2}\di{}\varepsilon,\nonumber\\
    {}&=\frac{1}{\sqrt{2\pi}}\left[-\varepsilon e^{-\frac{1}{2}\varepsilon^2}\right]_{c(TOL)}^{\infty}+\frac{1}{\sqrt{2\pi}}\int_{c(TOL)}^{\infty} e^{-\frac{1}{2}\varepsilon^2}\di{}\varepsilon,\nonumber\\
    {}&=\frac{1}{\sqrt{2\pi}}c(TOL)e^{-\frac{1}{2}(c(TOL))^2}+\frac{1}{\sqrt{2\pi}}\int_{c(TOL)}^{\infty} e^{-\frac{1}{2}\varepsilon^2}\di{}\varepsilon,\nonumber\\
    {}&<\frac{1}{\sqrt{2\pi}}c(TOL)e^{-\frac{1}{2}(c(TOL))^2}+\frac{1}{\sqrt{2\pi}}\frac{1}{c(TOL)} e^{-\frac{1}{2}(c(TOL))^2},\nonumber\\
    {}&<\frac{1}{\sqrt{2\pi}}(2(1+p))^{\frac{1}{2}}\log(TOL^{-1})^{\frac{1}{2}} TOL^{1+p}+\frac{1}{\sqrt{2\pi}}(2(1+p))^{-\frac{1}{2}}\log(TOL^{-1})^{-\frac{1}{2}} TOL^{1+p},\nonumber\\
    {}&=o(TOL)
\end{align}
as $TOL\to 0$ by a similar derivation as in \eqref{eq:T.truncated}. The bound on the first term in \eqref{eq:tr.error} follows analogously. Finally, it holds that
\begin{align}
    |I-I^{\rm{tr}}|{}&<\frac{2}{\sqrt{2\pi}}\left((2(1+p))^{\frac{1}{2}}\log(TOL^{-1})^{\frac{1}{2}} +(1+|T|)(2(1+p))^{-\frac{1}{2}}\log(TOL^{-1})^{-\frac{1}{2}}\right)TOL^{1+p},\nonumber\\
    {}&=o(TOL).
\end{align}
If we assume a standard normal distribution for the parameters of interest, the EIG has the following shape instead:
\begin{equation}
    I=\int_{\bb{R}^2}\log\left(\frac{1}{\sqrt{2\pi}}\int_{\bb{R}}e^{-\frac{1}{2}\left(G(\theta)+\varepsilon-G(\vartheta)\right)^2}\frac{1}{\sqrt{2\pi}}e^{-\frac{1}{2}\vartheta^2}\di{}\vartheta\right)\frac{1}{\sqrt{2\pi}}e^{-\frac{1}{2}\varepsilon^2}\di{}\varepsilon\frac{1}{\sqrt{2\pi}}e^{-\frac{1}{2}\theta^2}\di{}\theta.
\end{equation}
We introduce the following notations:
\begin{equation}
    \tilde{I}^{\rm{tr}}\coloneqq \int_{-c(TOL)}^{c(TOL)}\int_{-c(TOL)}^{c(TOL)}\log\left(\frac{1}{\sqrt{2\pi}}\int_{\bb{R}}e^{-\frac{1}{2}\left(G(\theta)+\varepsilon-G(\vartheta)\right)^2}\pi(\vartheta)\di{}\vartheta\right)\pi(\varepsilon)\di{}\varepsilon\pi(\theta)\di{}\theta,
\end{equation}
and
\begin{equation}
    I^{\rm{tr}}\coloneqq \int_{-c(TOL)}^{c(TOL)}\int_{-c(TOL)}^{c(TOL)}\log\left(\frac{1}{\sqrt{2\pi}}\int_{-c(TOL)}^{c(TOL)}e^{-\frac{1}{2}\left(G(\theta)+\varepsilon-G(\vartheta)\right)^2}\pi(\vartheta)\di{}\vartheta\right)\pi(\varepsilon)\di{}\varepsilon\pi(\theta)\di{}\theta.
\end{equation}
It follows that
\begin{equation}\label{eq:trunc.split}
    |I-I^{\rm{tr}}|\leq |I-\tilde{I}^{\rm{tr}}|+|\tilde{I}^{\rm{tr}}-I^{\rm{tr}}|.
\end{equation}
To bound the first term in~\eqref{eq:trunc.split}, we split the outer region of truncation into four parts as follows:
\begin{align}\label{eq:regions}
    \bb{R}^2\setminus[-c(TOL),c(TOL)]^2={}&[c(TOL),\infty]\times[-c(TOL),\infty]\nonumber\\
    {}&\cup[-c(TOL),\infty]\times[-\infty,-c(TOL)],\nonumber\\
    {}&\cup[-\infty,-c(TOL)]\times[-\infty,c(TOL)]\nonumber\\
    {}&\cup[-\infty,c(TOL)]\times[c(TOL),\infty].
\end{align}
From the triangle inequality, it follows that
\begin{align}
    |I-\tilde{I}^{\rm{tr}}|\leq{}&\int_{c(TOL)}^{\infty}\int_{-c(TOL)}^{\infty}\left|\log\left(\frac{1}{\sqrt{2\pi}}\int_{\bb{R}}e^{-\frac{1}{2}\left(G(\theta)+\varepsilon-G(\vartheta)\right)^2}\pi(\vartheta)\di{}\vartheta\right)\right|\pi(\varepsilon)\di{}\varepsilon\pi(\theta)\di{}\theta\nonumber\\
    {}&+\int_{-c(TOL)}^{\infty}\int_{-\infty}^{-c(TOL)}\left|\log\left(\frac{1}{\sqrt{2\pi}}\int_{\bb{R}}e^{-\frac{1}{2}\left(G(\theta)+\varepsilon-G(\vartheta)\right)^2}\pi(\vartheta)\di{}\vartheta\right)\right|\pi(\varepsilon)\di{}\varepsilon\pi(\theta)\di{}\theta\nonumber\\
    {}&+\int_{-\infty}^{-c(TOL)}\int_{-\infty}^{c(TOL)}\left|\log\left(\frac{1}{\sqrt{2\pi}}\int_{\bb{R}}e^{-\frac{1}{2}\left(G(\theta)+\varepsilon-G(\vartheta)\right)^2}\pi(\vartheta)\di{}\vartheta\right)\right|\pi(\varepsilon)\di{}\varepsilon\pi(\theta)\di{}\theta\nonumber\\
    {}&+\int_{-\infty}^{c(TOL)}\int_{c(TOL)}^{\infty}\left|\log\left(\frac{1}{\sqrt{2\pi}}\int_{\bb{R}}e^{-\frac{1}{2}\left(G(\theta)+\varepsilon-G(\vartheta)\right)^2}\pi(\vartheta)\di{}\vartheta\right)\right|\pi(\varepsilon)\di{}\varepsilon\pi(\theta)\di{}\theta.
\end{align}
From the Lipschitz property of $G$, it follows that
\begin{align}
    |G(z)|{}&=|G(z)-G(0)+G(0)|,\nonumber\\
    {}&\leq |G(z)-G(0)|+|G(0)|,\nonumber\\
    {}&\leq L|z|+|G(0)|,
\end{align}
and thus also that
\begin{align}
    \frac{1}{2}\left(G(\theta)+\varepsilon-G(\vartheta)\right)^2{}&\leq \left(G(\theta)+\varepsilon\right)^2+\left(G(\vartheta)\right)^2,\nonumber\\
    {}&\leq \left(G(\theta)+\varepsilon\right)^2+2\left(G(0)\right)^2+2L^2\vartheta^2.
\end{align}
We begin with the first region in~\eqref{eq:regions}, and note that bounds on the remaining regions follow from similar arguments. It holds that
\begin{align}
    {}&\int_{c(TOL)}^{\infty}\int_{-c(TOL)}^{\infty}\left|\log\left(\frac{1}{\sqrt{2\pi}}\int_{\bb{R}}e^{-\frac{1}{2}\left(G(\theta)+\varepsilon-G(\vartheta)\right)^2}\pi(\vartheta)\di{}\vartheta\right)\right|\pi(\varepsilon)\di{}\varepsilon\pi(\theta)\di{}\theta\nonumber\\
    {}&\leq\int_{c(TOL)}^{\infty}\int_{-c(TOL)}^{\infty}\left|\log\left(\frac{1}{\sqrt{2\pi}}\int_{\bb{R}}e^{-\left(G(\theta)+\varepsilon\right)^2-2\left(G(0)\right)^2-2L^2\vartheta^2}\pi(\vartheta)\di{}\vartheta\right)\right|\pi(\varepsilon)\di{}\varepsilon\pi(\theta)\di{}\theta,\nonumber\\
    {}&=\int_{c(TOL)}^{\infty}\int_{-c(TOL)}^{\infty}\left|R-\left(G(\theta)+\varepsilon\right)^2\right|\pi(\varepsilon)\di{}\varepsilon\pi(\theta)\di{}\theta,
\end{align}
where
\begin{align}
    R{}&\coloneqq \log\left(\frac{1}{\sqrt{2\pi}}\right)-2(G(0))^2+\log\left(\int_{\bb{R}}\pi\left(\sqrt{4L^2+1}\vartheta\right)\di{}\vartheta\right),\nonumber\\
    {}&=\log\left(\sqrt{\frac{4\pi(L^2+1)}{2\pi}}\right)-2(G(0))^2.
\end{align}
It follows that
\begin{align}
    \int_{c(TOL)}^{\infty}\int_{-c(TOL)}^{\infty}\left|R-\left(G(\theta)+\varepsilon\right)^2\right|\pi(\varepsilon)\di{}\varepsilon\pi(\theta)\di{}\theta{}&\leq\int_{c(TOL)}^{\infty}\int_{\bb{R}}\left|R-\left(G(\theta)+\varepsilon\right)^2\right|\pi(\varepsilon)\di{}\varepsilon\pi(\theta)\di{}\theta,\nonumber\\
    {}&\leq \int_{c(TOL)}^{\infty}\int_{\bb{R}}\left(|R|+2(G(\theta))^2+2\varepsilon^2\right)\pi(\varepsilon)\di{}\varepsilon\pi(\theta)\di{}\theta,\nonumber\\
    {}&= \int_{c(TOL)}^{\infty}\left(|R|+(G(\theta))^2\right)\pi(\theta)\di{}\theta+2\int_{c(TOL)}^{\infty}\pi(\theta)\di{}\theta,\nonumber\\
    {}&=o(TOL),
\end{align}
where the last line follows from a similar argument as in~\eqref{eq:T.truncated}. The bound on the other regions in~\eqref{eq:regions} follows analogously. Moreover, it follows for the second term in~\eqref{eq:trunc.split} that
\begin{align}\label{eq:log.bound}
    |\tilde{I}^{\rm{tr}}-I^{\rm{tr}}|{}&=\left|\int_{-c(TOL)}^{c(TOL)}\int_{-c(TOL)}^{c(TOL)}\log\left(\frac{\int_{\bb{R}}e^{-\frac{1}{2}\left(G(\theta)+\varepsilon-G(\vartheta)\right)^2}\pi(\vartheta)\di{}\vartheta}{\int_{-c(TOL)}^{c(TOL)}e^{-\frac{1}{2}\left(G(\theta)+\varepsilon-G(\vartheta)\right)^2}\pi(\vartheta)\di{}\vartheta}\right)\pi(\varepsilon)\di{}\varepsilon\pi(\theta)\di{}\theta\right|,\nonumber\\
    {}&=\int_{-c(TOL)}^{c(TOL)}\int_{-c(TOL)}^{c(TOL)}\log\left(\frac{\int_{\bb{R}}e^{-\frac{1}{2}\left(G(\theta)+\varepsilon-G(\vartheta)\right)^2}\pi(\vartheta)\di{}\vartheta}{\int_{-c(TOL)}^{c(TOL)}e^{-\frac{1}{2}\left(G(\theta)+\varepsilon-G(\vartheta)\right)^2}\pi(\vartheta)\di{}\vartheta}\right)\pi(\varepsilon)\di{}\varepsilon\pi(\theta)\di{}\theta,
\end{align}
where the term inside the logarithm is always greater or equal to one. Next, we apply the following bound:
\begin{align}\label{eq:num.denom}
    \log\left(\frac{\int_{\bb{R}}e^{-\frac{1}{2}\left(G(\theta)+\varepsilon-G(\vartheta)\right)^2}\pi(\vartheta)\di{}\vartheta}{\int_{-c(TOL)}^{c(TOL)}e^{-\frac{1}{2}\left(G(\theta)+\varepsilon-G(\vartheta)\right)^2}\pi(\vartheta)\di{}\vartheta}\right){}&=\log\left(1+\frac{\int_{|\vartheta|>c(TOL)}e^{-\frac{1}{2}\left(G(\theta)+\varepsilon-G(\vartheta)\right)^2}\pi(\vartheta)\di{}\vartheta}{\int_{-c(TOL)}^{c(TOL)}e^{-\frac{1}{2}\left(G(\theta)+\varepsilon-G(\vartheta)\right)^2}\pi(\vartheta)\di{}\vartheta}\right),\nonumber\\
    {}&\leq\frac{\int_{|\vartheta|>c(TOL)}e^{-\frac{1}{2}\left(G(\theta)+\varepsilon-G(\vartheta)\right)^2}\pi(\vartheta)\di{}\vartheta}{\int_{-c(TOL)}^{c(TOL)}e^{-\frac{1}{2}\left(G(\theta)+\varepsilon-G(\vartheta)\right)^2}\pi(\vartheta)\di{}\vartheta}.
\end{align}
We now derive an upper bound for the numerator and a lower bound for the denominator in~\eqref{eq:num.denom}. First, it holds that
\begin{equation}
    e^{-\frac{1}{2}\left(G(\theta)+\varepsilon-G(\vartheta)\right)^2}\leq 1,
\end{equation}
and thus the numerator is bounded using a similar argument as in~\eqref{eq:T.truncated} as follows:
\begin{align}
    \int_{|\vartheta|>c(TOL)}e^{-\frac{1}{2}\left(G(\theta)+\varepsilon-G(\vartheta)\right)^2}\pi(\vartheta)\di{}\vartheta{}&\leq \int_{|\vartheta|>c(TOL)}\pi(\vartheta)\di{}\vartheta,\nonumber\\
    {}&\leq\frac{2}{\sqrt{2\pi}}\frac{1}{c(TOL)} e^{-\frac{1}{2}(c(TOL))^2}.
\end{align}
To bound the denominator in~\eqref{eq:num.denom}, we note that $\theta,\vartheta\in[-c(TOL),c(TOL)]$ and thus
\begin{equation}
    \int_{-c(TOL)}^{c(TOL)}e^{-\frac{1}{2}\left(G(\theta)+\varepsilon-G(\vartheta)\right)^2}\pi(\vartheta)\di{}\vartheta\geq \int_{J_{\theta}}e^{-\frac{1}{2}\left(G(\theta)+\varepsilon-G(\vartheta)\right)^2}\pi(\vartheta)\di{}\vartheta
\end{equation}
for any interval $J_{\theta}\subseteq [-c(TOL),c(TOL)]$. To conveniently apply the Lipschitz bound on $G$ in the following step, we select $J_{\theta}$ such that it has length at most one. If $c(TOL)\leq 1/2$, then we choose 
\begin{equation}
    J_{\theta}\coloneqq[-c(TOL),c(TOL)];
\end{equation}
thus, $|J_{\theta}|\leq 1$. Otherwise, we choose
\begin{equation}
    J_{\theta}\coloneqq[\theta-j_{\theta},\theta-j_{\theta}+1],
\end{equation}
where
\begin{equation}
    j_{\theta}\coloneqq\min\{1,|\theta+c(TOL)|\}. 
\end{equation}
Thus, $\theta\in J_{\theta}$ and $|J_{\theta}|=1$. Note that other constructions for $J_{\theta}$ are possible. Then, it follows that
\begin{align}
    \int_{J_{\theta}}e^{-\frac{1}{2}\left(G(\theta)+\varepsilon-G(\vartheta)\right)^2}\pi(\vartheta)\di{}\vartheta{}&\geq \int_{J_{\theta}}e^{-\frac{1}{2}\left(|\varepsilon|+|G(\theta)-G(\vartheta)|\right)^2}\pi(\vartheta)\di{}\vartheta,\nonumber\\
    {}&\geq \int_{J_{\theta}}e^{-\frac{1}{2}\left(|\varepsilon|+L|\theta-\vartheta|\right)^2}\pi(\vartheta)\di{}\vartheta,\nonumber\\
    {}&\geq \int_{J_{\theta}}e^{-\frac{1}{2}\left(|\varepsilon|+L\right)^2}\pi(\vartheta)\di{}\vartheta.
\end{align}
Moreover, it holds that $|\vartheta|\leq |\theta|+1$ for $\vartheta\in J_{\theta}$; thus,
\begin{align}
    \int_{J_{\theta}}e^{-\frac{1}{2}\left(|\varepsilon|+L\right)^2}\pi(\vartheta)\di{}\vartheta{}&\geq \int_{J_{\theta}}e^{-\frac{1}{2}\left(|\varepsilon|+L\right)^2}\pi(|\theta|+1)\di{}\vartheta,\nonumber\\
    {}&=|J_{\theta}| e^{-\frac{1}{2}\left(|\varepsilon|+L\right)^2}\pi(|\theta|+1).
\end{align}
Combining the previous results and substituting into~\eqref{eq:log.bound}, we arrive at
\begin{align}
    {}&\int_{-c(TOL)}^{c(TOL)}\int_{-c(TOL)}^{c(TOL)}\log\left(\frac{\int_{\bb{R}}e^{-\frac{1}{2}\left(G(\theta)+\varepsilon-G(\vartheta)\right)^2}\pi(\vartheta)\di{}\vartheta}{\int_{-c(TOL)}^{c(TOL)}e^{-\frac{1}{2}\left(G(\theta)+\varepsilon-G(\vartheta)\right)^2}\pi(\vartheta)\di{}\vartheta}\right)\pi(\varepsilon)\di{}\varepsilon\pi(\theta)\di{}\theta\nonumber\\
    {}&\leq \int_{-c(TOL)}^{c(TOL)}\int_{-c(TOL)}^{c(TOL)}\frac{\frac{2}{\sqrt{2\pi}}\frac{1}{c(TOL)} e^{-\frac{1}{2}(c(TOL))^2}}{|J_{\theta}|e^{-\frac{1}{2}\left(|\varepsilon|+L\right)^2}\pi(|\theta|+1)}\pi(\varepsilon)\di{}\varepsilon\pi(\theta)\di{}\theta,\nonumber\\
    {}&=\frac{2}{\sqrt{2\pi}}\frac{e^{-\frac{1}{2}(c(TOL))^2}}{|J_{\theta}|c(TOL)}\int_{-c(TOL)}^{c(TOL)}\int_{-c(TOL)}^{c(TOL)}\frac{1}{\sqrt{2\pi}}e^{\frac{1}{2}\left(|\varepsilon|+L\right)^2+\frac{1}{2}\left(|\theta|+1\right)^2-\frac{1}{2}\varepsilon^2-\frac{1}{2}\theta^2}\di{}\varepsilon\di{}\theta,\nonumber\\
    {}&=\frac{1}{\pi}\frac{1}{|J_{\theta}|c(TOL)} e^{-\frac{1}{2}(c(TOL))^2}\int_{-c(TOL)}^{c(TOL)}\int_{-c(TOL)}^{c(TOL)}e^{|\varepsilon|L+\frac{1}{2}L^2+|\theta|+\frac{1}{2}}\di{}\varepsilon\di{}\theta.
\end{align}
Finally, we note that
\begin{align}
    \int_{-c(TOL)}^{c(TOL)}e^{|\theta|}\di{}\theta{}&=2\int_{0}^{c(TOL)}e^{\theta}\di{}\theta,\nonumber\\
    {}&=2\left(e^{c(TOL)}-1\right)
\end{align}
and
\begin{equation}
    \int_{-c(TOL)}^{c(TOL)}e^{|\varepsilon|L}\di{}\varepsilon=\frac{2}{L}\left(e^{c(TOL)}-1\right).
\end{equation}
Thus, it follows that
\begin{align}
    |\tilde{I}^{\rm{tr}}-I^{\rm{tr}}|{}&\leq \frac{4}{L\pi}\frac{1}{|J_{\theta}|c(TOL)} e^{-\frac{1}{2}(c(TOL))^2}e^{\frac{L^2+1}{2}}\left(e^{c(TOL)}-1\right)^2,\nonumber\\
    {}&=o(TOL).
\end{align}
Here, $|J_{\theta}|=\min\{1,2c(TOL)\}$. We are interested in the case where $TOL\to 0$ and thus $c(TOL)\to\infty$ and $|J_{\theta}|=1$ for $c(TOL)\geq 1/2$.

\end{proof}
\begin{rmk}[Extensions to general dimensions]
   The case for dimension $d>1$ follows immediately for a diagonal covariance matrix $\bs{\Sigma}_{\bs{\varepsilon}}=\bs{I}_{d\times d}$, yielding only weak dependence of the truncation error on the integrand dimension.
\end{rmk}
\begin{rmk}[Parameters of interest with Gaussian prior]\label{rmk:Lipschitz.Gaussian.prior}
    Assumption~\ref{asu:Lipschitz} generally does not hold for parameters of interest distributed according to a Gaussian distribution. However, a truncation scheme similar to that for the observation noise may be applied to $\bs{\tilde{G}}(\cdot)=\bs{G}(F_{\bs{\theta}}^{-1}(\cdot),\bs{\xi})$ if $\bs{G}$ is smooth in the first argument.
\end{rmk}
\section{Bounds on the generalized Hardy--Krause variation for the truncated integrands}\label{app:GHK}
We begin by introducing some definitions. For continuously differentiable integrands $\varphi$, the generalized Vitali variation of order one in dimension $d$ is defined as follows (see \cite[Chapter 13]{Dic10}):
    \begin{equation}\label{def:GVV}
    V^{(d)}(\varphi)\coloneqq \left(\int_{[0,1]^d}\left|\frac{\partial^d\varphi}{\partial z_1\ldots \partial z_d}(\bs{z})\right|^2\di{}\bs{z}\right)^{\frac{1}{2}}.
    \end{equation}
    Moreover, let $\cl{I}_{d}\coloneqq\{1,\ldots,d\}$, and $\mathfrak{u}$ be the set of all subsets of $\cl{I}_d$. We denote the cardinality of $\mathfrak{u}$ by $|\mathfrak{u}|$ and define
    \begin{equation}
        \varphi_{\mathfrak{u}}(\bs{z}_{\mathfrak{u}})\coloneqq \int_{[0,1]^{d-|\mathfrak{u}|}}\varphi(\bs{z})\di{}\bs{z}_{\cl{I}_d\setminus\mathfrak{u}},
    \end{equation}
    where $\bs{z}_{\mathfrak{u}}$ contains only those entries of $\bs{z}$ that are in $\mathfrak{u}$ and $\bs{z}_{\cl{I}_d\setminus\mathfrak{u}}$ contains only those entries of $\bs{z}$ that are not in $\mathfrak{u}$. Then, the generalized Hardy--Krause variation of order one of $\varphi$ on $[0,1]^d$ is defined as follows:
    \begin{equation}
        V_{\rm{GHK}}(\varphi)\coloneqq \left(\left|\int_{[0,1]^d}\varphi(\bs{z})\di{}\bs{z}\right|^2+\sum_{\emptyset\neq \mathfrak{u}\subseteq\cl{I}_{d}}\left(V^{(|\mathfrak{u}|)}(\varphi_{\mathfrak{u}})\right)^2\right)^{\frac{1}{2}}.
    \end{equation}
    Next, we demonstrate bounds on the generalized Hardy--Krause variation of order one for the inner and outer integrands in the EIG application, where the outer integrand is the logarithm and the inner integrand is the likelihood.
\begin{lem}\label{lem:GHK.inner}
    Let $g$ be as in~\eqref{def:inner.integrand} and define
    \begin{equation}
        g_{\bs{y}}(\bs{x})\coloneqq g(\bs{y},\bs{x})|_{\bs{y}},
    \end{equation}
where $g_{\bs{y}}(\bs{x}):[0,1]^{d_{2}}\to\bb{R}_{+}$ is the likelihood as a function of the second variable with the first variable fixed. The generalized Hardy--Krause variation of $g_{\bs{y}}$ is uniformly bounded for all $\bs{y}\in[0,1]^{d_1}$. Moreover, for a fixed and sufficiently large value of $c(TOL)$ as in \eqref{eq:c.TOL}, where $TOL>0$, there exists $0\leq a_{\epsilon}<\infty$ such that 
 \begin{equation}\label{eq:inv.ineq.GHK}
    \sup_{\bs{y}\in[0,1]^{2}} \left\lvert \frac{V_{\rm{GHK}}(g(\bs{y},\cdot))}{ \bar{g}(\bs{y})} \right\rvert \leq a_{\epsilon}TOL^{-\epsilon}
\end{equation}
for any $\epsilon>0$, where $a_{\epsilon}\to\infty$ as $\epsilon\to 0$.
\end{lem}
\begin{proof}
    Let $d_1=2$ $d_2=1$. It follows from Assumption~\ref{asu:Lipschitz} that
    \begin{align}
        \left(V^{(1)}(g_{\bs{y}})_{\{1\}}\right)^2{}&=\frac{1}{2\pi}\int_{[0,1]}\left|\frac{\di{}}{\di{}x}e^{-\frac{1}{2}\left(\tilde{G}(y_1)+\Phi^{-1}(y_2)-\tilde{G}(x)\right)^2}\right|^2\di{}x,\nonumber\\
        {}&=\frac{1}{2\pi}\int_{[0,1]}\left|\left(\tilde{G}(y_1)+\Phi^{-1}(y_2)-\tilde{G}(x)\right)\frac{\di{}}{\di{}x}\tilde{G}(x)e^{-\frac{1}{2}\left(\tilde{G}(y_1)+\Phi^{-1}(y_2)-\tilde{G}(x)\right)^2}\right|^2\di{}x,\nonumber\\
        {}&\leq \frac{L^2}{2\pi}.
    \end{align}
    The other term in the generalized Hardy--Krause variation is bounded by $1/2\pi$ by a similar derivation. To demonstrate the bound in~\eqref{eq:inv.ineq.GHK}, we must show that there exists $0<a_{\epsilon}<\infty$ such that
    \begin{align}
    \sup_{\bs{y}\in[0,1]^{2}}\left|\frac{\left(\left|\int_{[0,1]}e^{-\frac{1}{2}\left(z(\bs{y},x)\right)^2}\di{}x\right|^2+\int_{[0,1]}\left|\left(\frac{\di{}}{\di{} x}\tilde{G}_{h}(x)\right)z(\bs{y},x)e^{-\frac{1}{2}\left(z(\bs{y},x)\right)^2}\right|^2\di{}x\right)^{\frac{1}{2}}}{\int_{[0,1]}e^{-\frac{1}{2}\left(z(\bs{y},x)\right)^2}\di{}x}\right|\leq a_{\epsilon}TOL^{-\epsilon},
\end{align}
where
    \begin{equation}
        z(\bs{y},x)\coloneqq \tilde{G}_{h}(y_1)+F_{c(TOL)}^{-1}(y_2)-\tilde{G}_{h}(x),
    \end{equation}
for any $\epsilon>0$, where $a_{\epsilon}\to\infty$ as $TOL\to 0$. First, we use
\begin{equation}
    \left(\frac{\di{}}{\di{} x}\tilde{G}_{h}(x)\right)z(\bs{y},x)\leq L\left(L+c(TOL)\right)
\end{equation}
via Assumption~\ref{asu:Lipschitz} to find that
    \begin{align}\label{eq:sup.z}
    {}&\sup_{\bs{y}\in[0,1]^{2}}\left|\frac{\left(\left|\int_{[0,1]}e^{-\frac{1}{2}\left(z(\bs{y},x)\right)^2}\di{}x\right|^2+\int_{[0,1]}\left|\left(\frac{\di{}}{\di{} x}\tilde{G}_{h}(x)\right)z(\bs{y},x)e^{-\frac{1}{2}\left(z(\bs{y},x)\right)^2}\right|^2\di{}x\right)^{\frac{1}{2}}}{\int_{[0,1]}e^{-\frac{1}{2}\left(z(\bs{y},x)\right)^2}\di{}x}\right|\nonumber\\
    {}&= \sup_{\bs{y}\in[0,1]^{2}}\left|\left(1+\frac{\int_{[0,1]}\left|\left(\frac{\di{}}{\di{} x}\tilde{G}_{h}(x)\right)z(\bs{y},x)\right|^2e^{-\left(z(\bs{y},x)\right)^2}\di{}x}{\left(\int_{[0,1]}e^{-\frac{1}{2}\left(z(\bs{y},x)\right)^2}\di{}x\right)^2}\right)^{\frac{1}{2}}\right|,\nonumber\\
    {}&\leq \sup_{\bs{y}\in[0,1]^{2}}\left|\left(1+L^2\left(L+c(TOL)\right)^2\frac{\int_{[0,1]}e^{-\left(z(\bs{y},x)\right)^2}\di{}x}{\left(\int_{[0,1]}e^{-\frac{1}{2}\left(z(\bs{y},x)\right)^2}\di{}x\right)^2}\right)^{\frac{1}{2}}\right|.
\end{align}
We look at the numerator first and find that
\begin{align}
    \int_{[0,1]}e^{-\left(z(\bs{y},x)\right)^2}\di{}x{}&=\int_{[0,1]}e^{-\left(\tilde{G}_{h}(y_1)-\tilde{G}_{h}(x)\right)^2-2\left(\tilde{G}_{h}(y_1)-\tilde{G}_{h}(x)\right)F_{c(TOL)}^{-1}(y_2)-\left(F_{c(TOL)}^{-1}(y_2)\right)^2}\di{}x,\nonumber\\
    {}&\leq e^{-\left(F_{c(TOL)}^{-1}(y_2)\right)^2}\int_{[0,1]}e^{-\left(\tilde{G}_{h}(y_1)-\tilde{G}_{h}(x)\right)^2}e^{2\left|\left(\tilde{G}_{h}(y_1)-\tilde{G}_{h}(x)\right)\right|\left|F_{c(TOL)}^{-1}(y_2)\right|}\di{}x.
\end{align}
The first term inside the integral is bounded by one and from Assumption~\ref{asu:Lipschitz}, it follows that
\begin{align}\label{eq:ubn}
    \int_{[0,1]}e^{-\left(z(\bs{y},x)\right)^2}\di{}x{}&\leq e^{-\left(F_{c(TOL)}^{-1}(y_2)\right)^2}e^{2L\left|F_{c(TOL)}^{-1}(y_2)\right|}.
\end{align}
Next, we consider the denominator in~\eqref{eq:sup.z} to see that
\begin{align}
    \int_{[0,1]}e^{-\frac{1}{2}\left(z(\bs{y},x)\right)^2}\di{}x\geq e^{-\frac{1}{2}\left(F_{c(TOL)}^{-1}(y_2)\right)^2}e^{-\frac{L^2}{2}}e^{-L\left|F_{c(TOL)}^{-1}(y_2)\right|}.
\end{align}
Thus, it follows that
\begin{equation}\label{eq:lbd}
    \left(\int_{[0,1]}e^{-\frac{1}{2}\left(z(\bs{y},x)\right)^2}\di{}x\right)^2\geq e^{-\left(F_{c(TOL)}^{-1}(y_2)\right)^2}e^{-L^2}e^{-2L\left|F_{c(TOL)}^{-1}(y_2)\right|}.
\end{equation}
Combining the upper bound on the numerator in~\eqref{eq:ubn} and the lower bound on the denominator in~\eqref{eq:lbd}, we find that
\begin{align}
    \frac{\int_{[0,1]}e^{-\left(z(\bs{y},x)\right)^2}\di{}x}{\left(\int_{[0,1]}e^{-\frac{1}{2}\left(z(\bs{y},x)\right)^2}\di{}x\right)^2}{}&\leq e^{L^2}e^{4L\left|F_{c(TOL)}^{-1}(y_2)\right|},\nonumber\\
    {}&\leq e^{L^2}e^{4Lc(TOL)}.
\end{align}
We recall that $c(TOL)=(2(1+p))^{\frac{1}{2}}\log(TOL^{-1})^{\frac{1}{2}}$ and find that
\begin{align}
    e^{4Lc(TOL)}{}&=\left(e^{\frac{1}{2}(c(TOL))^2}\right)^{\frac{8L}{c(TOL)}},\nonumber\\
    {}&=\left(TOL^{-(1+p)}\right)^{\frac{8L}{c(TOL)}},\nonumber\\
    {}&=TOL^{-\frac{8L(1+p)}{c(TOL)}}.
\end{align}
As $TOL$ approaches zero, $c(TOL)$ approaches infinity and the exponent in the above expression approaches zero, that is,
\begin{equation}
    \frac{8L(1+p)}{c(TOL)}\to 0 \quad \text{as } TOL\to 0.
\end{equation}
Thus, there exists a constant $0<b_{\tilde{\epsilon}}<\infty$ such that
\begin{equation}
    \frac{\int_{[0,1]}e^{-\left(z(\bs{y},x)\right)^2}\di{}x}{\left(\int_{[0,1]}e^{-\frac{1}{2}\left(z(\bs{y},x)\right)^2}\di{}x\right)^2}\leq b_{\tilde{\epsilon}}TOL^{-\tilde{\epsilon}},
\end{equation}
for any $\tilde{\epsilon}>0$ where $b_{\tilde{\epsilon}}\to\infty$ as $\tilde{\epsilon}\to 0$. Finally, it follows that
\begin{align}
    \sup_{\bs{y}\in[0,1]^{2}}\left|\left(1+L^2\left(L+c(TOL)\right)^2\frac{\int_{[0,1]}e^{-\left(z(\bs{y},x)\right)^2}\di{}x}{\left(\int_{[0,1]}e^{-\frac{1}{2}\left(z(\bs{y},x)\right)^2}\di{}x\right)^2}\right)^{\frac{1}{2}}\right|\leq \left(1+L^2\left(L+c(TOL)\right)^2b_{\epsilon}TOL^{-\epsilon}\right)^{\frac{1}{2}},
\end{align}
and from this it follows that there exists a constant $0<a_{\epsilon}<\infty$ such that
\begin{align}
    \left(1+L^2\left(L+c(TOL)\right)^2b_{\epsilon}TOL^{-\epsilon}\right)^{\frac{1}{2}}{}&=\left(1+L^2\left(L+(2(1+p))^{\frac{1}{2}}\log(TOL^{-1})^{\frac{1}{2}}\right)^2b_{\epsilon}TOL^{-\epsilon}\right)^{\frac{1}{2}},\nonumber\\
    {}&\leq a_{\epsilon}TOL^{-\epsilon},
\end{align}
for any $\epsilon>0$, where $a_{\epsilon}\to\infty$ as $\epsilon\to 0$. The derivation for higher dimensions follows from the Fa\`{a} di Bruno formula for any positive definite covariance matrix of the observation noise.
\end{proof}
\begin{lem}\label{lem:GHK}
    Let $f\equiv\log$ and $\bar{g}$ be as in~\eqref{eq:bar.g}. The generalized Hardy--Krause variation of $f(\bar{g})$ is unbounded for $y_2\in[0,1]$. Moreover, given Assumption~\ref{asu:Lipschitz}, there exists $0<a_{\epsilon}<\infty$ such that the following bound holds for $y_2\in[\Phi(-c(TOL)),\Phi(c(TOL))]$:
    \begin{equation}
        V_{\rm{GHK}}(f(\bar{g}))\leq a_{\epsilon} TOL^{-\frac{1}{2}(1+p+\epsilon)},
    \end{equation}
    for any $\epsilon,p>0$ as $TOL\to 0$, where $a_{\epsilon}\to\infty$ as $\epsilon\to 0$. Here, $\Phi$ is the CDF of the standard normal and $c(TOL)$ is as in~\eqref{eq:c.TOL}. Furthermore, under the same assumptions, there exists $0<b_{\epsilon}<\infty$ such that the Hardy--Krause variation~\ref{eq:VHK} is bounded as follows:
    \begin{equation}
        V_{\rm{HK}}(f(\bar{g}))=b_{\epsilon}TOL^{-\epsilon},
    \end{equation}
    for any $\epsilon>0$ as $TOL\to 0$, where $b_{\epsilon}\to\infty$ as $\epsilon\to 0$.
\end{lem}
The above lemma implies that although the generalized Hardy--Krause variation of the truncated integrand is bounded for any $TOL>0$, it increases at a rate of essentially $TOL^{-\frac{1}{2}}$. The error bound on the outer variance in Corollary~\ref{cor:optimal.work.s} contains the squared generalized Hardy--Krause variation. Thus, the faster convergence observed for the scrambled Sobol' sequence again reduces to the rate predicted by the Koksma--Hlawka inequality as well as by Proposition~\ref{cor:he} based on Owen's boundary growth condition for the untruncated integrand in the EIG application. The inner variance still converges at a faster rate as this integrand has bounded generalized Hardy--Krause variation independent of $TOL$ with or without truncation.
\begin{proof}[Proof of Lemma\ref{lem:GHK}]
It follows immediately from the definition in~\eqref{def:GVV} that the generalized Vitali variation and subsequently also the generalized Hardy--Krause variation are unbounded for $\varphi\equiv f(\bar{g}(\bs{y}))$, where $f\equiv\log$ and $\bs{y}\in[0,1]^d$ due to the boundary singularity. Next, let $d_1=2$ and $d_2=1$. Moreover, let
\begin{equation}
    \bar{g}^{\rm{tr}}(\bs{y})\coloneqq \frac{1}{\sqrt{2\pi}}\int_{[0,1]}e^{-\frac{1}{2}\left(G(y_1)+F_{c(TOL)}^{-1}(y_2)-G(x)\right)^2}\di{} x,
\end{equation}
where $F_{c(TOL)}^{-1}$ is as in~\eqref{eq:inverse.CDF}. For this case, the generalized Hardy--Krause variation consists of four terms. We begin by bounding
\begin{align}
   &{}\left(V^{(1)}(f(\bar{g}^{\rm{tr}})_{\{2\}})\right)^2\nonumber\\
   {}&=(2\Phi(c(TOL))-1)\int_{[0,1]}\left|\frac{\partial}{\partial y_2}\int_{[0,1]}\log\left(\frac{1}{\sqrt{2\pi}}\int_{[0,1]}e^{-\frac{1}{2}\left(G(y_1)+F_{c(TOL)}^{-1}(y_2)-G(x)\right)^2}\di{} x\right)\di{}y_1\right|^2\di{}y_2,\nonumber\\
   {}&=\int_{J_{c(TOL)}}\left|\frac{\partial}{\partial y_2}\int_{[0,1]}\log\left(\frac{1}{\sqrt{2\pi}}\int_{[0,1]}e^{-\frac{1}{2}\left(G(y_1)+\Phi^{-1}(y_2)-G(x)\right)^2}\di{} x\right)\di{}y_1\right|^2\di{}y_2,\nonumber\\
\end{align}
where $J_{c(TOL)}\coloneqq [\Phi(-c(TOL)),\Phi(c(TOL))]$. It follows that
\begin{align}
   &{}\left(V^{(1)}(f(\bar{g}^{\rm{tr}})_{\{2\}})\right)^2\nonumber\\{}&=\int_{J_{c(TOL)}}\left|\int_{[0,1]}\frac{\frac{\partial}{\partial y_2}\left(\frac{1}{\sqrt{2\pi}}\int_{[0,1]}e^{-\frac{1}{2}\left(G(y_1)+\Phi^{-1}(y_2)-G(x)\right)^2}\di{} x\right)}{\frac{1}{\sqrt{2\pi}}\int_{[0,1]}e^{-\frac{1}{2}\left(G(y_1)+\Phi^{-1}(y_2)-G(x)\right)^2}\di{} x}\di{}y_1\right|^2\di{}y_2,\nonumber\\
   {}&=\int_{J_{c(TOL)}}\left|\int_{[0,1]}\frac{\frac{\di{}}{\di{}{y_2}}\Phi^{-1}(y_2)\int_{[0,1]}\left(G(y_1)+\Phi^{-1}(y_2)-G(x)\right)e^{-\frac{1}{2}\left(G(y_1)+\Phi^{-1}(y_2)-G(x)\right)^2}\di{} x}{\int_{[0,1]}e^{-\frac{1}{2}\left(G(y_1)+\Phi^{-1}(y_2)-G(x)\right)^2}\di{} x}\di{}y_1\right|^2\di{}y_2,\nonumber\\
   {}&\leq\int_{J_{c(TOL)}}\left|\int_{[0,1]}\frac{\left|\frac{\di{}}{\di{}{y_2}}\Phi^{-1}(y_2)\right|\int_{[0,1]}\left|G(y_1)+\Phi^{-1}(y_2)-G(x)\right|e^{-\frac{1}{2}\left(G(y_1)+\Phi^{-1}(y_2)-G(x)\right)^2}\di{} x}{\int_{[0,1]}e^{-\frac{1}{2}\left(G(y_1)+\Phi^{-1}(y_2)-G(x)\right)^2}\di{} x}\di{}y_1\right|^2\di{}y_2.
\end{align}
From Assumption~\ref{asu:Lipschitz}, it follows that
\begin{align}\label{eq:GHK.HK}
   &{}\left(V^{(1)}(f(\bar{g}^{\rm{tr}})_{\{2\}})\right)^2\leq\int_{J_{c(TOL)}}\left(\left|\frac{\di{}}{\di{}{y_2}}\Phi^{-1}(y_2)\right|\left(L+\left|\Phi^{-1}(y_2)\right|\right)\right)^2\di{}y_2.
\end{align}
We use the expression for derivative of the inverse CDF~\ref{eq:derivative.inverse.CDF} to obtain
\begin{align}\label{eq:GHK.transformation}
  \int_{J_{c(TOL)}}\left(\left|\frac{\di{}}{\di{}{y_2}}\Phi^{-1}(y_2)\right|\left(L+\left|\Phi^{-1}(y_2)\right|\right)\right)^2\di{}y_2{}&=\int_{-c(TOL)}^{c(TOL)}\left(
  \sqrt{2\pi}e^{\frac{1}{2}\varepsilon^2}\left(L+|\varepsilon|\right)\right)^2\frac{1}{\sqrt{2\pi}}e^{-\frac{1}{2}\varepsilon^2}\di{}\varepsilon\nonumber\\
  {}&=\int_{-c(TOL)}^{c(TOL)}
  \sqrt{2\pi}e^{\frac{1}{2}\varepsilon^2}\left(L+|\varepsilon|\right)^2\di{}\varepsilon\nonumber\\
  {}&=2\int_{0}^{c(TOL)}
  \sqrt{2\pi}e^{\frac{1}{2}\varepsilon^2}\left(L+\varepsilon\right)^2\di{}\varepsilon\nonumber\\
  {}&=2\sqrt{2\pi}\int_{0}^{c(TOL)}
  e^{\frac{1}{2}\varepsilon^2}\left(L^2+2L\varepsilon+\varepsilon^2\right)\di{}\varepsilon.
\end{align}
We first consider the term
\begin{align}
    2L\int_{0}^{c(TOL)}
  \varepsilon e^{\frac{1}{2}\varepsilon^2}\di{}\varepsilon{}&=2L\left(e^{\frac{1}{2}(c(TOL))^2}-1\right).
\end{align}
Next, we use integration by parts to obtain
\begin{align}
   \int_{0}^{c(TOL)}
  \varepsilon^2 e^{\frac{1}{2}\varepsilon^2}\di{}\varepsilon{}&= \int_{0}^{c(TOL)}
  \varepsilon\cdot\varepsilon e^{\frac{1}{2}\varepsilon^2}\di{}\varepsilon,\nonumber\\
  {}&=\left(c(TOL)e^{\frac{1}{2}(c(TOL))^2}-\int_{0}^{c(TOL)}e^{\frac{1}{2}\varepsilon^2}\di{}\varepsilon\right),\nonumber\\
  {}&\leq c(TOL)e^{\frac{1}{2}(c(TOL))^2}.
\end{align}
Combining these terms, it follows that
\begin{equation}
    \left(V^{(1)}(f(\bar{g}^{\rm{tr}})_{\{2\}})\right)^2\leq 2\sqrt{2\pi}\left(2L\left(e^{\frac{1}{2}(c(TOL))^2}-1\right) + c(TOL)e^{\frac{1}{2}(c(TOL))^2}+L^2\int_{0}^{c(TOL)}e^{\frac{1}{2}\varepsilon^2}\di{}\varepsilon\right).
\end{equation}
The remaining integral can be bounded as follows:
\begin{equation}
    \int_{0}^{c(TOL)}e^{\frac{1}{2}\varepsilon^2}\di{}\varepsilon\leq c(TOL)e^{\frac{1}{2}\left(c(TOL)\right)^2},
\end{equation}
resulting in
\begin{align}
    \left(V^{(1)}(f(\bar{g}^{\rm{tr}})_{\{2\}})\right)^2{}&\leq e^{\frac{1}{2}(c(TOL))^2}\left(2L+c(TOL) + L^2c(TOL)\right)2\sqrt{2\pi},\nonumber\\
    {}&=TOL^{-(1+p)}\left(L^2+\left(L^2+1\right)(2(1+p))^{\frac{1}{2}}\log(TOL^{-1})^{\frac{1}{2}}\right)2\sqrt{2\pi}.
\end{align}
Finally, there exists $0<a_{\epsilon}<\infty$ such that
\begin{equation}\label{eq:GHK.bound.TOL}
    TOL^{-(1+p)}\left(L^2+\left(L^2+1\right)(2(1+p))^{\frac{1}{2}}\log(TOL^{-1})^{\frac{1}{2}}\right)2\sqrt{2\pi}\leq a_{\epsilon}TOL^{-(1+p+\epsilon)}
\end{equation}
for any $\epsilon>0$ as $TOL\to 0$, where $a_{\epsilon}\to\infty$ as $\epsilon\to 0$. The derivation of the bounds for the other terms in the generalized Hardy--Krause variations  follows along similar lines. If we retrace the steps in the above proof up to~\eqref{eq:GHK.HK} for the Hardy--Krause variation rather than the stricter generalized Hardy--Krause variation of order one, we end up with a similar expression, except that the integrand is not squared. It follows directly from the transformations in~\eqref{eq:GHK.transformation} that the resulting bound contains only terms that are logarithmic in $TOL$; thus, there exists $0<b_{\epsilon}<\infty$ such that
\begin{align}\label{eq:HK.transformation}
  \int_{J_{c(TOL)}}\left|\frac{\di{}}{\di{}{y_2}}\Phi^{-1}(y_2)\right|\left(L+\left|\Phi^{-1}(y_2)\right|\right)\di{}y_2{}&=\int_{-c(TOL)}^{c(TOL)}\left(L+|\varepsilon|\right)\di{}\varepsilon\nonumber\\
  {}&=2\int_{0}^{c(TOL)}\left(L+\varepsilon\right)\di{}\varepsilon\nonumber\\
  {}&=2Lc(TOL)+\left(c(TOL)\right)^2,\nonumber\\
  {}&=2L(2(1+p))^{\frac{1}{2}}\log(TOL^{-1})^{\frac{1}{2}}+2(1+p)\log(TOL^{-1}),\nonumber\\
  {}&\leq b_{\epsilon}TOL^{-\epsilon},
\end{align}
for any $\epsilon>0$ as $TOL\to 0$, where $b_{\epsilon}\to\infty$ as $\epsilon\to 0$.
\end{proof}
This result demonstrates that although the error bound based on the generalized Hardy--Krause variation of order one converges faster by one order of magnitude compared to the bound in the Koksma--Hlawka inequality based on the Hardy--Krause variation, we lose that one order of magnitude when estimating the EIG with truncation. This again indicates comparable behavior as predicted from the bound on the untruncated EIG based on Owen's boundary growth condition in Proposition~\ref{cor:he}. The derivation of the bounds for higher dimensions follows from the Fa\`{a} di Bruno formula, and we also expect a similar behavior, up to multiplicative logarithmic factors. On the other hand, from the above results, we expect the bound on the generalized Hardy--Krause variation of order one in~\eqref{eq:GHK.bound.TOL} for the multidimensional noise case to behave as $\cl{O}\left(TOL^{-d(1+p+\epsilon)}\right)$ for $\bs{y}\in[0,1]^d$. This would imply an error bound that is significantly worse than the ones implied by the Hardy--Krause variation, Proposition~\ref{cor:he}, or even the CLT for higher dimensions. Similar explosion results were derived in~\cite[Theorem 6.2]{Kaa24} for another regularity criterion also implying $L^2$ regularity of the mixed derivatives. We did not observe such a behavior in our numerical experiments, where it appears that the multiplicative error term $C_{\rm{Q}}^{(1)}$ only increases approximately at rate one even in dimension $d=4$ (see Fig.~\ref{fig:ex0.5.pilot.truncation.c}).

\section{Total error of the EIG estimation via the rDLQMC estimator}\label{app:Tota.error}
\begin{cor}[Total error of the EIG estimation via the rDLQMC estimator]\label{cor:total.error}
    Under the assumptions stated in Proposition~\ref{prop:B.DLQ}, Proposition~\ref{prop:V.DLQ}, Lemma~\ref{lem:cond}, and Lemma~\ref{lem:trunc}, there exists $0<\delta_{K}<1$ such that
    \begin{align}
    |I-I^{\rm{tr}}|={}&o(TOL),\\
        \left|\bb{E}[I^{\rm{tr}}_{\rm{rDLQ}}]-I^{\rm{tr}}\right|\leq{}& C_{\mathrm{disc}}h^{\eta}+ \frac{\bb{E}[|\bar{g}_{h}|^2\left|f''(\bar{g}_{h})\right|]\tilde{L}^{2d_2}(c(TOL))^{2d_2}C_{\epsilon,d_2}^2}{2M^{2-2\epsilon}}+\frac{\bb{E}[\left|\bar{g}_{h}\right|^3\left|f'''(K(TOL)\bar{g}_{h})\right|]\tilde{L}^{3d_2}(c(TOL))^{3d_2}C_{\epsilon,d_2}^3}{6M^{3-3\epsilon}},\label{eq:eig.bias}\\
        \bb{V}[I^{\rm{tr}}]\leq{}&\frac{2b^2B_A^2C_{\epsilon, d_1}^2}{N^{2-2\epsilon-2\max_{i}A_i}} +\frac{2\mathbb{E}[|\bar{g}_h|^2\left|f'(\bar{g}_h)\right|^2]\tilde{L}^{2d_2}(c(TOL))^{2d_2}C_{\epsilon, d_2}^2}{NM^{2-2\epsilon}}\nonumber\\
        {}&+\frac{\mathbb{E}[|\bar g_h|^4\left|f''(K(TOL)\bar{g}_h)\right|^2](\Gamma_{d_1}+1)\tilde{L}^{4d_2}(c(TOL))^{4d_2}C_{\epsilon, d_2}^4}{2NM^{4-4\epsilon}}\label{eq:eig.variance},
    \end{align}
    for all $K(TOL)\in(1-\delta_{K},1+\delta_{K})$ for a specified error tolerance $TOL>0$, where $I^{\rm{tr}}$ is the truncated integral introduced in \eqref{eq:truncation.split} with $c(TOL)$ as in \eqref{eq:c.TOL}. Moreover, this result holds for any $\epsilon>0$, where $C_{\epsilon,d_1},C_{\epsilon,d_2}\to\infty$ as $\epsilon\to 0$, $B_A\to\infty$ as $\min_{i}A_i\to 0$, and $\Gamma_{d_1}\to\infty$ as $d_1\to\infty$.
\end{cor}
\begin{proof}
    The result follows directly from Proposition~\ref{prop:B.DLQ}, Proposition~\ref{prop:V.DLQ}, Lemma~\ref{lem:cond}, and Lemma~\ref{lem:trunc}. For the logarithm, it follows that $|f'(\bar{g}_{h})|=|\bar{g}_{h}^{-1}|$, $|f''(\bar{g}_{h})|=|\bar{g}_{h}^{-2}|$, and $|f'''(\bar{g}_{h})|=|2\bar{g}_{h}^{-3}|$ and thus the terms on the right-hand sides of \eqref{eq:eig.bias} and \eqref{eq:eig.variance} are bounded. From the introduction of $K$ in \eqref{eq:monotonicity}, it immediately follows that the condition that $K(TOL)\in(1-\delta_{K},1+\delta_{K})$ for $0<\delta_{K}<1$ is equivalent to
    \begin{equation}
        C_{\epsilon,d_2}M^{-1+\epsilon}\tilde{L}^{d_2}c(TOL)^{d_2}<1,
    \end{equation}
     which is verified for $M>(C_{\epsilon,d_2}\tilde{L}^{d_2}c(TOL)^{d_2})^{\frac{1}{1-\epsilon}}$,
     imposing a weaker requirement than the condition in \eqref{eq:Bias.Constraint}.
\end{proof}
\begin{rmk}[Total complexity of the EIG estimation via the rDLQMC estimator]
    The total computational complexity to estimate the EIG for a prescribed tolerance ($TOL$) via the rDLQMC estimator differs from the result stated in Proposition~\ref{prop:W.DLQ} by a log factor introduced by the restriction to truncated observation noise. This immediately follows from Proposition~\ref{prop:W.DLQ} for the error bounds presented in Corollary~\ref{cor:total.error}.
\end{rmk}
\begin{rmk}[Total complexity of the EIG estimation for scrambled Sobol' sequence]
    When using rDLQMC estimators based on scrambled Sobol' sequences for EIG estimation with truncation, improved inner sampling bias and variance errors can be obtained compared to the bounds in Corollary~\ref{cor:total.error} (see Corollary~\ref{cor:optimal.work.s} and Lemma~\ref{lem:GHK.inner}). In particular, the inner sampling bias error can be significantly reduced due to the smoothness of the inner integrand.
\end{rmk}
\clearpage
\section{Algorithm for the pilot to estimate convergence rates and multiplicative constants}\label{app:pilot}
\begin{algorithm}
\caption{Pilot to estimate convergence rates and multiplicative constants}
\SetKwInOut{Input}{Input}
\SetKwInOut{Output}{Output}
\Input{Functions $f: \mathbb{R} \to \mathbb{R}$, $g: [0,1]^{d_1+d_2} \to \mathbb{R}$}
\Output{Constants and rates $C_{{\mathrm{Q}}}^{(1)}$, $C_{{\mathrm{Q}}}^{(2)}$, $C_{{\mathrm{Q}}}^{(3)}$, $\beta$, $\delta$}
\For{$s = 1$ to $S_{\mathrm{pilot}}$}
{Generate outer Sobol' points $\{\bs{y}^{(s,n)}\}_{n=1}^{N_{\mathrm{pilot}}}$ using scrambling $\bs{\tau}^{(s)}$\\
\For{$n = 1$ to $N_{\mathrm{pilot}}$}
{
\For{$r = 1$ to $R_{\mathrm{pilot}}$}
{
Generate inner Sobol' points $\{\bs{x}^{(s,n,r,m)}\}_{m=1}^{M_{\mathrm{pilot}}}$ using scrambling $\bs{\rho}^{(s,n,r)}$\\
\For{$m = 1$ to $M_{\mathrm{pilot}}$}
{
Compute $g(\bs{y}^{(s,n)}, \bs{x}^{(s,n,r,m)})$
}
}
}
}
Compute outer sample variance to estimate $\mathbb{V}[\mathbb{E}[\frac{1}{N}\sum_{n=1}^Nf(\hat{g}_M(\bs{y}^{(n)}))|\{\bs{y}^{(n)}\}_{n=1}^N]]$\\
Fit outer multiplicative constant $C_{{\mathrm{Q}}}^{(1)}$ and convergence rate $1+\beta$
\\Compute inner sample variance to estimate $\mathbb{E}[\mathbb{V}[\frac{1}{N}\sum_{n=1}^Nf(\hat{g}_M(\bs{y}^{(n)}))|\{\bs{y}^{(n)}\}_{n=1}^N]]$\\
Fit inner multiplicative constant $C_{{\mathrm{Q}}}^{(2)}$ and convergence rate $1+\delta$
\\Compute sample mean to estimate $\frac{1}{2}|\bb{E}[f''(\bar{g}(\bs{y}))\left(\Delta g_M(\bs{y})\right)^2]|$\\
Fit inner multiplicative constant $C_{{\mathrm{Q}}}^{(3)}$ and convergence rate $1+\delta$
\label{alg.pil}
\end{algorithm}

\section{Derivation of the finite element formulation}\label{ap:FEM}
We let $(\Omega,\cl{F},\bb{P})$ be a complete probability space with outcomes $\Omega$, $\sigma$-field $\cl{F}$, and probability measure $\bb{P}$. Moreover, $\cl{H}\coloneqq H^1(\cl{D})\oplus H^1(\cl{D})$ is defined as the space of the solution for the coupled thermomechanical fields $(\vartheta(\omega), \bs{u}(\omega))$ for $\omega\in\Omega$, where $H^1(\cl{D})$ represents the standard Sobolev space $W^{1,2}(\cl{D})$ with its corresponding Sobolev norm. Furthermore, the Bochner space is defined as follows:
\begin{equation}
   V_{\vartheta}\oplus V_U\equiv L^2_{\bb{P}}(\Omega;\cl{H})\coloneqq \left\{(\vartheta,\bs{u})\colon\Omega\to\cl{H}\quad \text{s.t.}\,\left(\int_{\Omega}\lVert (\vartheta(\omega),\bs{u}(\omega))\rVert_{\cl{H}}^2\di{}\bb{P}(\omega)\right)^{1/2}<\infty\right\}.
\end{equation}
We aimed to determine $(\vartheta,\bs{u})\in L^2_{\bb{P}}(\Omega;\cl{H})$ such that the weak formulations, \eqref{eq:weak.temperature} and \eqref{eq:weak.mechanical}, are fulfilled for all $(\hat{\vartheta},\bs{\hat{u}})\in L^2_{\bb{P}}(\Omega;\cl{H})$.



\footnotesize

\begin{thebibliography}{1}

\bibitem{Rya03}
Ryan KJ.
\newblock {Estimating expected information gains for experimental designs with
  application to the random fatigue-limit model}.
\newblock {\em {Journal of Computational and Graphical Statistics}},
  12:585--603, 2003.

\bibitem{Bec18}
Beck J, Dia BM, Espath LF, Long Q, and Tempone R.
\newblock Fast {B}ayesian experimental design: {L}aplace-based importance
  sampling for the expected information gain.
\newblock {\em Computer Methods in Applied Mechanics and Engineering},
  334:523--553, 2018.

\bibitem{Bec20}
Beck J, Dia BM, Espath LF, and Tempone R.
\newblock Multilevel double loop {M}onte {C}arlo and stochastic collocation
  methods with importance sampling for {B}ayesian optimal experimental design.
\newblock {\em International Journal for Numerical Methods in Engineering},
  121:3482--3503, 2020.

\bibitem{Car20}
Carlon AG, Dia BM, Espath LF, Lopez RH, and Tempone R.
\newblock Nesterov-aided stochastic gradient methods using {L}aplace
  approximation for {B}ayesian design optimization.
\newblock {\em Computer Methods in Applied Mechanics and Engineering},
  363:112909, 2020.

\bibitem{Sha48}
Shannon CE.
\newblock {A mathematical theory of communication}.
\newblock {\em {Bell System Technical Journal}}, 27:379--423, 1948.

\bibitem{Kul59}
Kullback S.
\newblock {\em Information theory and statistics}.
\newblock Wiley, 1959.

\bibitem{Lon13}
Long Q, Scavino M, Tempone R, and Wang S.
\newblock Fast estimation of expected information gains for {B}ayesian
  experimental designs based on {L}aplace approximations.
\newblock {\em Computer Methods in Applied Mechanics and Engineering},
  259:24--39, 2013.

\bibitem{Lon21}
Long Q.
\newblock Multimodal information gain in {B}ayesian design of experiments.
\newblock {\em Computational Statistics}, 37:865--885, 2022.

\bibitem{Sch20}
Schillings C, Sprungk B, and Wacker P.
\newblock On the convergence of the {L}aplace approximation and noise-level-robustness of {L}aplace-based {M}onte {C}arlo methods for {B}ayesian inverse problems.
\newblock {\em Numerische Mathematik},
  145:915--971, 2020.

\bibitem{Wac17}
Wacker P.
\newblock Laplace's method in {B}ayesian inverse problems.
\newblock {\em arXiv preprint arXiv:1701.07989 [math.PR]}, 2017.

\bibitem{Cha95}
Chaloner K and Verdinelli I.
\newblock Bayesian experimental design: a review.
\newblock {\em Statistical Science},
  10:273--304, 1995.

\bibitem{Kul51}
Kullback S and Leibler RA.
\newblock On information and sufficiency.
\newblock {\em Annals of Mathematical Statistics},
  22:79--86, 1951.

\bibitem{Sti86}
Stigler SM.
\newblock Laplace's 1774 memoir on inverse probability.
\newblock {\em Statistical Science},
  1:359--378, 1986.

\bibitem{Lin56}
Lindley DV.
\newblock On a measure of information provided by an experiment.
\newblock {\em Annals of Mathematical Statistics},
  27:986--1005, 1956.

\bibitem{Tie86}
Tierney L and Kadane JB.
\newblock Accurate approximations for posterior moments and
marginal densities.
\newblock {\em Journal of the American Statistical Association},
  81:82--86, 1986.

\bibitem{Tie89}
Tierney L, Kass RE, and Kadane JB.
\newblock Fully exponential {L}aplace approximations to expectations and variances of nonpositive functions.
\newblock {\em Journal of the American Statistical Association},
  84:710--716, 1989.

\bibitem{Kas90}
Kass RE, Tierney L, and Kadane JB.
\newblock The validity of posterior expansions based on {L}aplace's method.
\newblock {\em in: Geisser S, Hodges JS, Press SJ, and Zellner A (Eds.), Essays
in Honor of George Barnard},
  473--488, 1990.

\bibitem{Lon15}
Long Q, Scavino M, Tempone R, and Wang S.
\newblock A {L}aplace method for under-determined {B}ayesian optimal experimental designs.
\newblock {\em Computer Methods in Applied Mechanics and Engineering},
  285:849--876, 2015.

\bibitem{Gil08}
Giles MB.
\newblock Multilevel {M}onte {C}arlo path simulation.
\newblock {\em Operations Research},
56(3):607--617, 2008.

\bibitem{Tsi17}
Tsilifis P, Ghanem RG, and Hajali P.
\newblock Efficient {B}ayesian experimentation using an expected information gain lower bound.
\newblock {\em SIAM/ASA Journal on Uncertainty Quantification},
5:30--62, 2017.

\bibitem{God18}
Goda T, Murakami D, Tanaka K, and Sato K.
\newblock Decision-theoretic sensitivity analysis for reservoir development under uncertainty using multilevel quasi-{M}onte {C}arlo methods.
\newblock {\em Computational Geosciences},
22:1009--1020, 2018.

\bibitem{Xu20}
Xu Z, He Z, and Wang X.
\newblock Efficient risk estimation via nested multilevel quasi-{M}onte {C}arlo simulation.
\newblock {\em Journal of Computational and Applied Mathematics}, 
443:115745, 2024.

\bibitem{Fan22}
Fang W, Wang Z, Giles MB, Jackson CH, Welton NJ, Andrieu C, and Thom H.
\newblock
Multilevel and quasi {M}onte {C}arlo methods for the calculation of the expected value of partial perfect information.
\newblock {\em Medical Decision Making},
42(2):168--181, 2022.

\bibitem{Caf98}
Caflisch RE.
\newblock
Monte {C}arlo and quasi-{M}onte {C}arlo methods.
\newblock {\em Acta Numerica},
7:1--49, 1998.

\bibitem{Hic98}
Hickernell FJ.
\newblock
Lattice rules: {H}ow well do they measure up?
\newblock {\em In: Hellekalek P and Larcher G (Eds.). Random and Quasi-Random Point Sets. Lecture Notes in Statistics 138},
\newblock {New York, Springer},
109--166, 1998.

\bibitem{Nie92}
Niederreiter H.
\newblock
Random number generation and quasi-{M}onte {C}arlo methods.
\newblock {\em SIAM},
\newblock {Philadelphia},
1992.

\bibitem{Owe03}
Owen AB.
\newblock
Quasi-{M}onte {C}arlo sampling.
\newblock {\em In: Monte Carlo Ray Tracing: SIGGRAPH},
69--88, 2003.

\bibitem{Tuf04}
Tuffin B.
\newblock
Randomization of quasi-{M}onte {C}arlo methods for error estimation: survey and normal approximation.
\newblock {\em In: Monte Carlo Methods and Applications},
\newblock {De Gruyter},
10(3-4):617--628, 2004.


\bibitem{Dur19}
Durrett R.
\newblock
Probability: {T}heory and examples.
\newblock {\em Cambridge University Press},
2019.

\bibitem{Dic10}
Dick J and Pillichshammer F.
\newblock
Digital nets and sequences: discrepancy theory and quasi-{M}onte {C}arlo integration.
\newblock {\em Cambridge University Press},
2010.


\bibitem{Dic13}
Dick J, Kuo FY, and Sloan IH.
\newblock
High-dimensional integration: the quasi-{M}onte {C}arlo way.
\newblock {\em Acta Numerica},
22:133--288, 2013.

\bibitem{Sob67}
Sobol' IM.
\newblock
Distribution of points in a cube and approximate evaluation of integrals.
\newblock {\em \v{Z}. Vy\v{c}isl. Mat. i Mat. Fiz. (in Russian)},
7:784--802, 1967.

\bibitem{Dro18}
Drovandi CC and Tran M-N.
\newblock
Improving the efficiency of fully {B}ayesian optimal design of experiments using randomised {Q}uasi-{M}onte {C}arlo.
\newblock {\em Bayesian Analysis},
13(1):139--162, 2018.

\bibitem{Lec10}
L'Ecuyer P, Munger D, and Tuffin B.
\newblock
On the distribution of integration error by randomly-shifted lattice rules.
\newblock {\em Electronic Journal of Statistics},
4:950--993, 2010.

\bibitem{Gob22}
Gobet E, Lerasle M, and M\'{e}tivier D.
\newblock
Mean estimation for randomized {Q}uasi {M}onte {C}arlo method.
\newblock {\em Hal preprint hal-03631879v2},
2022.

\bibitem{Hel22}
Helin T and Kretschmann R.
\newblock Non-asymptotic error estimates for the {L}aplace approximation in {B}ayesian inverse problems.
\newblock {\em Numerische Mathematik},
  150:521--549, 2022.

\bibitem{Lec18}
L'Ecuyer P.
\newblock Randomized {Q}uasi-{M}onte {C}arlo: {A}n introduction for practitioners.
\newblock {\em In: Owen A and Glynn P (Eds.). Monte Carlo and Quasi-Monte Carlo Methods. MCQMC 2016. Springer Proceedings in Mathematics \& Statistics},
\newblock {Cham, Springer}
  241:29--52, 2018.

\bibitem{Hla61}
Hlawka E.
\newblock Funktionen von beschr\"{a}nkter Variation in der Theorie der Gleichverteilung.
\newblock {\em Annali di Matematica Pura ed Applicata},
  54:325--333, 1961.

\bibitem{Owe08}
Owen AB.
\newblock Local antithetic sampling with scrambled nets.
\newblock {\em The Annals of Statistics},
  36(5):2319--2343, 2008.

\bibitem{Far91}
Farhat C,  Park KC, and Dubois-Pelerin Y.
\newblock An unconditionally stable staggered algorithm for transient finite element analysis of coupled thermoelastic problems.
\newblock {\em Computer Methods in Applied Mechanics and
Engineering},
  85(3):349--365, 1991.

\bibitem{haji2016multi}
Haji-Ali A-L,  Nobile F, Tamellini, L, and Tempone R.
\newblock Multi-index stochastic collocation convergence rates for random {PDE}s with parametric regularity.
\newblock {\em Foundations of Computational Mathematics},
  16:1555--1605, 2016.
  
\bibitem{Owe06}
Owen AB.
\newblock Halton sequences avoid the origin.
\newblock {\em SIAM Review},
  48(3):487--503, 2006.

\bibitem{Liu23}
Liu Y and Tempone R.
\newblock Error analysis of randomized quasi-{M}onte {C}arlo: {N}on-asymptotic error bound, importance sampling and application to linear elliptic {PDE}s with lognormal coefficients.
\newblock {\em Journal of Computational and Applied Mathematics}, 482:117310, 2026.

\bibitem{Gra15}
Graham IG, Kuo FY, Nichols JA, Scheichl R, Schwab C, and Sloan IH.
\newblock Quasi-{M}onte {C}arlo finite element methods for elliptic {PDE}s with lognormal random coefficients.
\newblock {\em Numerische Mathematik},
  131:329--368, 2015.

\bibitem{He23}
He Z, Zheng Z, and Wang X.
\newblock On the error rate of importance sampling with randomized quasi-{M}onte {C}arlo.
\newblock {\em SIAM Journal on Numerical Analysis},
  61(2):515--538, 2023.

\bibitem{Ouy23}
Ouyang D, Wang X, and He Z.
\newblock Achieving high convergence rates by quasi-{M}onte {C}arlo and importance sampling for unbounded integrands.
\newblock {\em SIAM Journal on Numerical Analysis}, 62(5):2393--2414, 2024.

\bibitem{Owe95}
Owen  AB.
\newblock Randomly permuted (t,m,s)-nets and (t,s)-sequences.
\newblock {\em in: Niederreiter H and Shiue PJS (Eds.), Monte Carlo and Quasi-Monte Carlo methods in Scientific Computing. Lecture Notes in Statistics 106},
\newblock {New York, Springer},
1995.

\bibitem{Rai18}
Rainforth T, Cornish R, Yang H, Warrington A, and Wood F.
\newblock On nesting Monte Carlo estimators.
\newblock {\em Proceedings of the $35^{\text{th}}$ International Conference on Machine Learning}, 80:4267--4276,
2018.

\bibitem{Hon09}
Hong LJ, and Juneja S.
\newblock Estimating the mean of a non-linear function of conditional expectation.
\newblock {\em in: Rossetti MD, Hill RR, Johansson B, Dunkin A and Ingalls RG (Eds.), Proceedings of the 2009 Winter Simulation Conference (WSC)}
\newblock {Austin},
1223--1236, 2009.



\bibitem{Lec23}
L'Ecuyer P, Nakayama MK, Owen AB, and Tuffin B.
\newblock Confidence intervals for randomized quasi-{M}onte {C}arlo estimators.
\newblock {\em Winter Simulation Conference},
445--456, 2023.

\bibitem{Gob17}
Fort G, Gobet E, and Moulines E.
\newblock MCMC design-based non-parametric regression for rare event. Application to nested risk computations.
\newblock {\em Monte Carlo Methods and Applications},
23(1):21--42, 2017.

\bibitem{Fos17}
Foster A, Jankowiak M, Bingham E, Horsfall P, Teh YW, Rainforth T, and Goodman N.
\newblock Variational {B}ayesian optimal experimental design.
\newblock {\em Advances in Neural Information Processing Systems 32},
2019.

\bibitem{God20}
Goda T, Hironaka T, and Iwamoto T.
\newblock Multilevel {M}onte {C}arlo estimation of expected information gains.
\newblock {\em Stochastic Analysis and Applications},
38(4):581--600, 2020.

\bibitem{God22}
Goda T, Hironaka T, Kitade W, and Foster A.
\newblock Unbiased MLMC stochastic gradient-based optimization of {B}ayesian experimental designs.
\newblock {\em SIAM Journal on Scientific Computing},
44(1):A286--A311, 2022.

\bibitem{Kuo12}
Kuo FY, Schwab C, and Sloan IH.
\newblock Quasi-Monte {C}arlo finite element methods for a class of elliptic partial differential equations with random coefficients.
\newblock {\em SIAM Journal on Numerical Analysis},
50(6):3351--3374, 2012.

\bibitem{Kuo15}
Kuo FY, Schwab C, and Sloan IH.
\newblock Multi-level quasi-{M}onte {C}arlo finite element methods for a class of elliptic PDEs with random coefficients.
\newblock {\em Foundations of Computational Mathematics},
15:411--449, 2015.



\bibitem{Owe97}
Owen AB.
\newblock Scrambled net variance for integrals of smooth functions.
\newblock {\em The Annals of Statistics},
25(4):1541--1562, 1997.

\bibitem{Cra76}
Cranley R, and Patterson TNL.
\newblock Randomization of number theoretic methods for multiple integration.
\newblock {\em SIAM Journal on Numerical Analysis},
13(6):904--914, 1976.

\bibitem{Pan23}
Pan Z, and Owen AB.
\newblock The nonzero gain coefficients of {S}obol’s sequences are always powers of two.
\newblock {\em Journal of Complexity},
75:101700, 2023.

\bibitem{Bor11}
Bornkamp B.
\newblock Approximating probability densities by iterated {L}aplace approximations.
\newblock {\em Journal of Computational and Graphical Statistics},
20(3):656--669, 2011.

\bibitem{Kaa24}
Kaarnioja V, and Schillings C.
\newblock Quasi-{M}onte {C}arlo for {B}ayesian design of
experiment problems governed by parametric {PDE}s.
\newblock {\em arXiv preprint arXiv:2405.03529v3 [math.NA]}, 2026.

\bibitem{Pat96}
Patel JK, and Read CB.
\newblock Handbook of the normal distribution.
\newblock {\em 2nd ed.},
\newblock {New York, Marcel Dekker},
1996.

\bibitem{Fen19}
Feng C, and Marzouk YM.
\newblock A layered multiple importance sampling scheme for focused optimal
  {B}ayesian experimental design.
\newblock {\em arXiv preprint arXiv:1903.11187 [stat.CO]}, 2019.

\bibitem{Bar25}
Bartuska A.
\newblock Hierarchical methods for {B}ayesian
optimal experimental design.
\newblock {\em Ph.D. thesis, RWTH Aachen University}, 2025.

\bibitem{Col14}
Collier N, Haji-Ali, AL, Nobile F, von Schwerin E, and Tempone
R.
\newblock A continuation multilevel {M}onte {C}arlo algorithm.
\newblock {\em BIT Numerical Mathematics}, 55(2):399–432, 2015.
\end{thebibliography}
\end{document}